\documentclass[runningheads,envcountsame]{llncs}
\usepackage[T1]{fontenc}
\usepackage{graphicx}
\usepackage{physics}
\usepackage{amssymb}
\usepackage{mathtools}
\usepackage{stmaryrd}
\usepackage{tikz-cd}
\usepackage{adjustbox}
\usepackage[bookmarks,colorlinks]{hyperref}
\hypersetup{
    colorlinks=true,
    linkcolor=black,
    filecolor=magenta,
    urlcolor=black,
    citecolor=blue
    }
\usepackage[capitalise]{cleveref}
\usepackage{cmll}
\usepackage{xspace}
\usepackage{thmtools}

\newcommand{\secref}[1]{\S \ref{#1}}

\newcommand{\Schrod}{Schr\"odinger}%
\newcommand{\independent}{\rotatebox[origin=c]{90}{$\models$}}

\DeclareMathOperator{\hagertimes}{\ensuremath{\stackrel{\scriptscriptstyle h}{\otimes}}}%
\DeclareMathOperator{\projtens}{\ensuremath{\widehat{\otimes}}}%
\DeclareMathOperator{\ptimes}{\projtens}%
\DeclareMathOperator{\pBantimes}{\mathbin{\ensuremath{\stackrel{\scriptscriptstyle \gamma}{\otimes}}}}%
\DeclareMathOperator{\iBantimes}{\mathbin{\ensuremath{\stackrel{\scriptscriptstyle \epsilon}{\otimes}}}}%
\DeclareMathOperator{\injtens}{\ensuremath{\widecheck{\otimes}}}%
\DeclareMathOperator{\itimes}{\injtens}%
\DeclareMathOperator{\alphatens}{\ensuremath{\stackrel{\scriptscriptstyle \alpha}{\otimes}}}
\DeclareMathOperator{\aconj}{\&}%

\makeatletter
\DeclareRobustCommand\widecheck[1]{{\mathpalette\@widecheck{#1}}}
\def\@widecheck#1#2{%
    \setbox\z@\hbox{\m@th$#1#2$}%
    \setbox\tw@\hbox{\m@th$#1%
       \widehat{%
          \vrule\@width\z@\@height\ht\z@
          \vrule\@height\z@\@width\wd\z@}$}%
    \dp\tw@-\ht\z@
    \@tempdima\ht\z@ \advance\@tempdima2\ht\tw@ \divide\@tempdima\thr@@
    \setbox\tw@\hbox{%
       \raise\@tempdima\hbox{\scalebox{1}[-1]{\lower\@tempdima\box
\tw@}}}%
    {\ooalign{\box\tw@ \cr \box\z@}}}
\makeatother

\newcommand{\fdOS}{\ensuremath{\mathbf{FdOS}}}%
\newcommand{\FdOS}{\fdOS}%
\newcommand{\FdOScb}{\fdOS_{\mathrm{cb}}}%
\newcommand{\FdBan}{\ensuremath{\mathbf{FdBan}}}%
\newcommand{\QCoh}{\ensuremath{\mathbf{Q}}}%
\newcommand{\vNalg}{\ensuremath{\mathbf{vNAlg_{fd}}}}%
\newcommand{\vNAlg}{\vNalg}%
\newcommand{\vNcoalg}{\ensuremath{\mathbf{vNCoalg_{fd}}}}%
\newcommand{\vNCoalg}{\vNcoalg}%
\newcommand{\FdVect}{\ensuremath{\mathbf{FdVect}}}%
\newcommand{\HH}{\ensuremath{\mathbf{H}}}%
\newcommand{\SSS}{\ensuremath{\mathbf{S}}}%

\newcommand{\loneplus}{\mathbin{\ensuremath{\stackrel{\scriptscriptstyle 1}{\oplus}}}}%
\newcommand{\linftyplus}{\mathbin{\ensuremath{\stackrel{\scriptscriptstyle{\infty}}{\oplus}}}}%
\newcommand{\CB}{\ensuremath{\mathrm{CB}}}%
\newcommand{\MM}{\ensuremath{\mathbb{M}}}%

\newcommand{\Ob}{\ensuremath{\mathrm{Ob}}}%
\newcommand{\QQ}{\ensuremath{\mathbf{Q}}}%
\newcommand{\polar}{\ensuremath{\circ}}%
\newcommand{\bipolar}{{\ensuremath{\circ\circ}}}%
\newcommand{\trnorm}[1]{\ensuremath{\norm{#1}_{\mathrm{tr}}}}%%

\newcommand{\id}{\ensuremath{\mathrm{\id}}}
\newcommand{\Ball}{\ensuremath{\mathrm{Ball}}}

\newcommand{\defeq}{\triangleq}
\newcommand{\eqdef}{\defeq}
\newcommand{\fullsubright}{\ensuremath{\xhookrightarrow{\mathrm{full}}}}%

\newcommand{\qsw}{\ensuremath{\mathrm{qsw}}}%

\begin{document}
\title{Quantum Coherence Spaces Revisited:\\ A von Neumann (Co)Algebraic Approach}
\titlerunning{Quantum Coherence Spaces Revisited}
\author{Thea Li \and Vladimir Zamdzhiev}
\institute{Université Paris-Saclay, CNRS, ENS Paris-Saclay, Inria, Laboratoire Méthodes Formelles,
91190, Gif-sur-Yvette, France
\\
\email{thea.li@inria.fr, vladimir.zamdzhiev@inria.fr}
}
\authorrunning{Thea Li \and Vladimir Zamdzhiev}
\maketitle

\begin{abstract}
   We describe a categorical model of MALL (Multiplicative Additive Linear
   Logic) inspired by the Heisenberg-\Schrod{} duality of finite-dimensional
   quantum theory. Proofs of formulas with positive logical polarity correspond
   to CPTP (completely positive trace-preserving) maps in our model, i.e. the
   quantum operations in the \Schrod{} picture, whereas proofs of formulas with
   negative logical polarity correspond to CPU (completely positive unital)
   maps, i.e. the quantum operations in the Heisenberg picture. The
   mathematical development is based on noncommutative geometry and
   finite-dimensional von Neumann (co)algebras, which can be defined as special
   kinds of (co)monoid objects internal to the category of finite-dimensional
   operator spaces.
\keywords{Quantum Theory, Linear Logic, Categorical Semantics}
\end{abstract}

\section{Introduction}

Linear Logic \cite{linear-logic} was discovered by Girard while studying the
mathematical model of \emph{coherence spaces} \cite{systemF,linear-logic}.
Coherence spaces ``played an essential role in the discovery of linear logic''
\cite[p. 967]{pcoh-ll} and clearly they have had a major impact on the study of
(semantics of) Linear Logic (LL). Later, \emph{Probabilistic Coherence Spaces}
(PCSs) \cite{pcoh-ll} were described and they give a model of LL which is
suitable for modeling discrete probability. The study of PCSs has resulted in
impressive results, such as full abstraction for probabilistic PCF
\cite{ppcf-full-abstraction} and also full abstraction for a more complicated
probabilistic programming language with recursive types
\cite{pcbv-full-abstraction}.

A natural next step in the development of these models is to consider
\emph{quantum coherence spaces} (QCSs) with the intention of taking the
semantic study of LL to quantum theory. Girard already proposed models of QCSs
\cite{qcs1,qcs2,qcs3} and they were later studied in \cite{qcs-baratella}, but
this approach has already been criticised by Selinger \cite{qcs-selinger} as
being inappropriate for quantum theory. One of the main problems with this
proposal is that the morphisms between QCSs do not appear to
correspond to \emph{completely positive maps} on relevant homsets (see also \cite[p.
2]{qcs-baratella}). This is problematic from the point of view of quantum
theory, because \emph{quantum operations} (also known as \emph{quantum
channels}) are modelled mathematically as certain kinds of completely positive
maps.

Inspired by the success of (probabilistic) coherence spaces, in this work we
set out to paint a new picture of what a quantum coherence space ought to be.
We consider finite-dimensional quantum theory and we describe a model of
MALL based on the following natural ideas: (1) proofs $P \vdash R$ of formulas with
\emph{positive} logical polarities should admit an interpretation as CPTP
(completely positive trace-preserving) maps, i.e. as the quantum operations in
the \emph{\Schrod{} picture} of quantum theory; (2) proofs $M \vdash N$ of
formulas with \emph{negative} logical polarities should admit an interpretation
as CPU (completely positive unital) maps, i.e. as the quantum operations in the
\emph{Heisenberg picture} of quantum theory; (3) the LL duality should coincide
with the Heisenberg-\Schrod{} duality of quantum theory on polarised formulas.
These desiderata naturally lead to the theory of operator spaces
\cite{er2000operator,blecher-merdy,Pisier_2003} which has excellent
categorical properties \cite{category-os} and which can be used to construct a
model of (full) LL whose duality is compatible with the Heisenberg-\Schrod{}
duality on the level of objects/formulas \cite{os-lics}. However, in
\cite{os-lics}, the morphisms/proofs do \emph{not} correspond to the quantum
operations in either picture. We address this problem (which is
suggested for future work in \cite{os-lics}) for \emph{finite-dimensional} (f.d.)
quantum theory and MALL.

In \secref{sec:background}, we provide relevant background on operator spaces.
We begin our contributions in \secref{sec:fdos}, where we show that $\FdOS$,
the category of f.d. operator spaces, is a model of MALL, and we prove other
relevant categorical properties. The category has a very important monoidal
structure $(\FdOS, \mathbb C, \hagertimes)$ that is outside the scope of LL --
it is given by the Haagerup tensor \cite[\S 9]{er2000operator},
\cite[pp. 30 -- 34]{blecher-merdy}, \cite[\S 5]{Pisier_2003}, which is one of the highlights of
operator space theory. This monoidal structure allows us to give a new, but
equivalent, definition of a f.d. \emph{von Neumann algebra} (vN-algebra) in
\secref{sec:vNAlg} as a certain kind of internal monoid object in $(\FdOS,
\mathbb C, \hagertimes)$.
It is well-known that CPU maps may be defined on
vN-algebras and this allows us to identify the subcategory of $\FdOS$
corresponding to negative logical polarity and the Heisenberg picture. The
novel categorical definition of f.d. vN-algebras, and the self-duality of the
Haagerup tensor, allow us to easily dualise this definition and we introduce
f.d. \emph{von Neumann coalgebras} (vN-coalgebras) in \secref{sec:vNCoalg} as
certain kinds of internal comonoid objects in $(\FdOS, \mathbb C,
\hagertimes).$ We show that CPTP maps are a natural notion of morphism between
vN-coalgebras and so we identify the subcategory of $\FdOS$ corresponding to
positive logical polarity and the \Schrod{} picture. In \secref{sec:hs} we
prove that there is a categorical duality between the    subcategories, which
gives a categorical formulation of the Heisenberg-\Schrod{} duality. Finally,
in \secref{sec:qcs}, we apply the semantic technique of \emph{gluing and orthogonality}
of Hyland and Schalk \cite{HS} to $\FdOS$ which results in a model of MALL
that satisfies the desiderata (1) -- (3), see Figure \ref{fig:the-two-tables}.

\section{Background on Operator Spaces and Banach Spaces}
\label{sec:background}
Operator spaces are widely seen as the ``non-commutative'' or ``quantised''
generalisation of Banach spaces. In this section we recall some basic
background material on these topics from several books on the matter \cite{er2000operator},
\cite{blecher-merdy}, \cite{Pisier_2003}.

%%%%%%%%%%%%%%%%%%%%%%%%%%%%%%%%%%%%%%
\subsection{Normed Spaces and Banach Spaces}
\label{sub:banach}
%%%%%%%%%%%%%%%%%%%%%%%%%%%%%%%%%%%%%%

A \emph{normed space} is a pair $(X, \norm{\cdot}_X)$ where $X$ is a complex vector
space and $\norm{\cdot}_X$ a norm on it. A \emph{Banach space} is a normed
space that is complete with respect to the induced norm topology. In this paper
we work with finite-dimensional normed spaces and every such space is
necessarily a Banach space. We usually simply write $X$ to refer to a Banach
space and we just write $\norm{\cdot}$ for its norm, which should be
clear from context.
Recall that, if $X$ and $Y$ are two Banach spaces and $\varphi \colon X \to Y$ is a linear map
between them, then its \emph{operator norm} is defined by the assignment
$ \norm{\varphi} \eqdef \sup \{ \norm{\varphi(x)} \ :\ x \in X, \norm{x} \leq 1 \} . $
Then $\varphi$ is continuous iff $\varphi$ is bounded iff $\norm{\varphi} < \infty.$
A linear map $f \colon X \to Y$ between two Banach spaces is called an \emph{isometry}
if $\norm{f(x)} = \norm{x}$ for each $x \in X$ and it is called a
\emph{contraction} whenever $\norm{f} \leq 1$, equivalently when $\norm{f(x)}
\leq \norm{x}$ for every $x \in X.$
We write $\FdBan$ for the category whose
objects are finite-dimensional Banach spaces and whose morphisms are the linear
contractions between them. In order to gain some intuition for the
constructions that follow, it is useful to recall some of the categorical
properties of $\FdBan$. It has a zero object given by the zero-dimensional Banach space $0.$
The binary categorical product of $X$ and $Y$ is given by
the Banach space $X \linftyplus Y \eqdef (X \oplus Y, \norm{\cdot}_{\ell^\infty} ) $,
where $\norm{(x,y)}_{\ell^\infty} \eqdef \mathrm{max}(\norm{x}, \norm{y})$ is simply
the $\ell^\infty$-norm and $X \oplus Y$ is the usual direct sum of vector spaces.
The binary categorical coproduct of $X$ and $Y$ is given by
the Banach space $X \loneplus Y \eqdef (X \oplus Y, \norm{\cdot}_{\ell^1} ) $,
where $\norm{(x,y)}_{\ell^1} \eqdef \norm{x} + \norm{y}$ is simply
the $\ell^1$-norm.
The category $\FdBan$ is symmetric monoidal closed with respect to the
\emph{projective} tensor product $\pBantimes$ (details omitted) with internal
hom given by $B(X,Y)$, the space of bounded linear maps between $X$ and $Y$.
In fact, $\FdBan$ is a
model of MALL \cite[\S 2.6]{blanco-phd}: the additive conjunction (disjunction)
is modelled by the categorical (co)products; the multiplicative conjunction
(disjunction) corresponds to $\pBantimes$ ($\iBantimes$), where $\iBantimes$
stands for the \emph{injective} tensor product; negation corresponds to the
dual Banach space $X^* \eqdef B(X, \mathbb C).$

\subsection{Operator Spaces and their Morphisms}
An (abstract) operator space may be seen as a vector space $X$ equipped with
additional norms on matrices with entries in $X$ that satisfy some conditions.
\begin{definition}[Matrix Space]
  Given a complex vector space $X$, we write $\mathbb{M}_{n,m}(X)$ for the
  vector space consisting of the \emph{$n\times m$} matrices with entries in
  $X$. The vector space structure is defined in the obvious way, i.e.
  componentwise. We also use the shorthand notation $\mathbb{M}_n(X)
  \eqdef \mathbb M_{n,n}(X)$ and $\mathbb M_n \eqdef \mathbb M_n(\mathbb C).$
\end{definition}
\begin{definition}[{Matrix Norm \cite[p. 20]{er2000operator}}]
  A \emph{matrix norm} on a vector space $X$ is an assignment of a norm
  $\|\cdot \|_n$ on $\mathbb{M}_n(X)$ for each $n\in \mathbb{N}$.
  Given such a norm, we write $M_n(X)$ for the normed space $(\mathbb{M}_n(X), \| \cdot \|_n)$.
\end{definition}
\begin{example}
  \label{ex:complex-numbers}
  We write $M_{n,m}$ for the normed space $(\mathbb M_{n,m}, \norm{\cdot}_{\mathrm{op}})$, where
  $\norm{\cdot}_{\mathrm{op}}$ is the usual operator norm on $\mathbb M_{n,m}$.
  The normed spaces $M_n \defeq M_{n,n}$ for $n \in \mathbb N$ determine a
  matrix norm on $\mathbb C$ that we use throughout the paper.
\end{example}
\begin{definition}[{Operator Space \cite[pp. 34--35]{Pisier_2003}}]
  \label{def:operator-space}
  An \emph{(abstract) operator space} is a complex vector space $X$ equipped with a matrix norm
  that satisfies:
  \begin{enumerate}
    \item[(B)] The pair $(\MM_1(X), \norm{\cdot}_1)$ is a Banach space;
    \item[(M1)] $\|x\oplus y\|_{m+n}=\max\{\|x\|_m,\|y\|_n\}$
    \item[(M2)] $\|\alpha x\beta\|_m\leq\|\alpha\|\|x\|_m\|\beta\|$
  \end{enumerate}
  for $n,m \in \mathbb N$, $x\in M_m(X)$, $y \in M_n(X)$, $\alpha,\beta\in M_{m}$.
  We write
  $
    x\oplus y\eqdef
    \begin{bmatrix}
      x & 0\\
      0 & y
    \end{bmatrix}
    \in M_{m+n}(X)
  $
  for the indicated matrix and $\alpha x \beta \in M_m(X)$ for the matrix that results from the obvious generalisation of matrix
  multiplication.
\end{definition}
Note that $\alpha, \beta$ in the above definition are scalar (i.e. complex) matrices, but $x$ need not be one.
A matrix norm that satisfies the above criteria is called an \emph{operator space structure} (o.s.s.) on $X$.
The obvious isomorphism $M_1(X) \cong X$ determines a norm on $X$ for which we write $\norm{\cdot}_X$ or simply
$\norm{\cdot}$ in the sequel.
If $x \in M_n(X)$ we often simply write $\norm{x}$ for $\norm{x}_n.$
Standard textbook results show that for each operator space $X$, both
$M_n(X)$ and $X$ are Banach spaces.
In this paper we work with \emph{finite-dimensional} operator spaces, i.e.
operator spaces whose underlying vector space is finite-dimensional. The
condition (B) is then automatically satisfied.
\begin{example}
  \label{ex:Mn}
  The complex numbers $\mathbb C$ is an operator space when equipped with the matrix norm from Example \ref{ex:complex-numbers}.
  More generally, each vector space $\mathbb M_{n,m}$ can be equipped with
  an o.s.s by defining a norm on $\mathbb M_k(\mathbb M_{n,m})$ via the linear isomorphism $\mathbb M_k(\mathbb M_{n,m}) \cong M_{nk,mk}$
  and the norm on the latter space (see Example \ref{ex:complex-numbers}). We write $M_{n,m}$ for this operator space. Note that
  this notation is compatible with Example \ref{ex:complex-numbers}, i.e. the norm $\norm{\cdot}_{M_{n,m}}$ on $M_{n,m}$ is exactly the
  operator norm. The operator spaces $M_n \eqdef M_{n,n}$ are fundamental for the theory of operator spaces and also for the results that we present in this paper.
\end{example}
\begin{example}
  For an operator space $X$, each of the spaces $M_n(X)$ has a canonical
  o.s.s that is determined by the linear isomorphism
  $\MM_m(M_n(X)) \cong M_{mn}(X)$ and the already given norm on the latter space.
  More generally, each matrix space $\MM_{m,n}(X)$ also has a canonical o.s.s that can be
  defined via the linear embedding $\MM_{m,n}(X) \hookrightarrow M_p(X)$,
  where $p = \mathrm{max}(m,n)$, and the o.s.s of the latter space \cite[\S 2.1]{er2000operator}
  and we write $M_{m,n}(X)$ for this operator space.
\end{example}

Next, we introduce important classes of morphisms between operator spaces.

\begin{definition}[{\cite[\S 2.2]{er2000operator}}]
  Let $X$ and $Y$ be two vector spaces and $\varphi \colon X \to Y$
  a linear map between them. The \emph{$n$-th amplification} of $\varphi$ is the linear function
  $ \varphi_n : \mathbb{M}_n(X)  \to \mathbb{M}_n(Y) :: [x_{ij}] \mapsto [\varphi(x_{ij})] , $
  i.e. the componentwise application of $\varphi$ to the matrix $[x_{ij}]$.
\end{definition}

\begin{definition}[{\cite[\S 2.2]{er2000operator}}]
  Let $X$ and $Y$ be operator spaces and $\varphi \colon X \to Y$ a linear map.
  Consider the assignment
  $  \| \varphi \|_\text{cb} \eqdef
    \sup \left\{ \norm{\varphi_n} \ : \ n \in \mathbb{N} \right\} ,
  $
  where $\norm{\varphi_n}$ is the operator norm of $\varphi_n \colon M_n(X) \to M_n(Y).$
  We say that $\varphi$ is
  \begin{itemize}
    \item a \emph{completely bounded} map (c.b) if $\| \varphi \|_\text{cb} < \infty$;
    \item a \emph{complete contraction} (c.c) if $\| \varphi \|_\text{cb} \leq 1$, i.e. each $\varphi_n$ is a contraction;
    \item a \emph{complete isometry} (c.i) if each $\varphi_n$ is an isometry;
    \item a \emph{complete quotient} map (c.q) if each $\varphi_n$ is a quotient map, i.e.  $\varphi_n$ maps the open unit ball of $M_n(X)$ \emph{onto} the open unit ball of $M_n(Y)$.
    \item a \emph{completely isometric isomorphism} (c.i.i) if $\varphi$ is a surjective c.i.
  \end{itemize}
\end{definition}

Every linear map $\varphi \colon X \to Y$ between two finite-dimensional
operator spaces is (completely) bounded \cite[Corollary 2.2.4]{er2000operator}
and therefore $\varphi$ and all of its
amplifications $\varphi_n$ are continuous (and bounded). The morphisms that are
central to our development are the complete contractions. We write $\mathbf{(Fd)OS}$
for the category whose objects are the (f.d.) operator spaces with
morphisms the complete contractions between them. We write $\FdVect$ for the
category of finite-dimensional complex vector spaces with linear maps as
morphisms and we write $U \colon \FdOS \to \FdBan$ and
$V \colon \FdOS \to \FdVect$ for the forgetful functors.

The next proposition is useful because it simplifies some norm computations.

\begin{proposition}[{\cite[\S 2.2, \S 3.3]{er2000operator}}]
  \label{prop:cb-op}
  Let $\varphi \colon X \to Y$ be a c.b. map. If $X = \mathbb C$ or $Y = \mathbb C$, then $\norm{\varphi}_{cb} = \norm{\varphi},$
  i.e. the cb-norm and operator norm coincide.
\end{proposition}

The vector space $\CB(X,Y)$, consisting of the completely bounded maps between $X$ and $Y$,
becomes a Banach space when equipped with the cb-norm $\norm{\cdot}_{cb}$ and it becomes
an operator space when equipped with the matrix norm given by the linear isomorphism
$\MM_n( \CB(X,Y) )\cong \CB(X,M_n(Y))$ and the Banach space structure of the
latter space \cite[\S 3.2]{er2000operator}. Note that in $\FdOS$, the
underlying vector space of $\CB(X,Y)$ is exactly $L(X,Y),$ the vector space of linear maps between $X$ and $Y$, i.e. $V(\CB(X,Y)) = L(X,Y).$

If $X$ is an operator space, then its (operator space) \emph{dual} is the operator space $X^* \eqdef \CB(X, \mathbb{C}).$
This is consistent with the notation for the Banach space dual, i.e. we have a Banach space equality $U(X^*) = (UX)^*$ \cite[\S 3.2]{er2000operator}.
It is also consistent with the vector space dual for \emph{finite-dimensional}
spaces. That is, for $Y \in \FdVect$, we define $Y^* \eqdef L(Y, \mathbb C)$
and then for $X \in \FdOS$ we also have $V(X^*) = (VX)^*$. In the sequel, we
show that $\CB(X,Y)$ is the internal hom of $\FdOS$ (just like it is for $\mathbf{OS}$ \cite{os-lics})
and $X^*$ is the dual in the sense of $*$-autonomy.

\begin{example}
  \label{ex:Tn}
  The vector space $\MM_n$ can be equipped with another important o.s.s. via the linear isomorphism
  $\MM_n \cong M_n^* :: a \mapsto (b \mapsto \trace(ab))$
  and the o.s.s. of the latter space. We write $T_n$ for the resulting operator space.
  Here, $\trace$ stands for the usual trace functional and $ab$ for matrix multiplication.
  The induced norm on $T_n$ coincides with the \emph{trace norm} and we have a completely isometric isomorphism
  $T_n \cong M_n^*.$ This is known as the \emph{trace class} o.s.s. \cite[1.4.5]{blecher-merdy}.
  It follows that $\tr : T_n \to \mathbb C$ is a complete contraction, but $\tr : M_n \to \mathbb C$ is \emph{not} one for $n > 1.$
\end{example}

\begin{remark}
  \label{rem:ci-dual-cq}
  Complete isometries and complete quotient maps are dual notions. A
  c.b. map $\varphi \colon X \to Y$ is a complete isometry (complete quotient map) iff $\varphi^* \colon Y^* \to X^*$ is a complete quotient map (complete isometry) \cite[1.4.3]{blecher-merdy}.
  The complete isometries coincide with the strong monomorphisms and the complete quotient maps coincide with the strong epimorphisms in $\mathbf{FdOS}$ \cite[\S 4.3]{category-os}.
\end{remark}

Since $V(M_n) = \MM_n = V(T_n)$, the usual definition of a completely positive map makes sense for linear maps $M_n \to M_m$ and $T_n \to T_m.$
The next example shows that the quantum operations in the Heisenberg and \Schrod{} pictures provide two important classes of completely contractive maps.

\begin{example}
  Every CPTP map $\varphi \colon T_n \to T_m$ is a complete contraction and every CPU map $\psi \colon M_m \to M_n$ is also a complete contraction
  with $\norm{\varphi}_{cb} = 1 = \norm{\psi}_{cb}.$
\end{example}

%%%%%%%%%%%%%%%%%%%%%%%%%%%%%%%%%%%%%%
\subsection{Direct Sums and (Co)products}
\label{sub:direct-sums}
%%%%%%%%%%%%%%%%%%%%%%%%%%%%%%%%%%%%%%

The category $\FdOS$ has finite (co)products which follows immediately from
\cite[\S 4.2]{category-os}. We recall these constructions, which are actually
completely analogous and compatible with the (co)products in $\FdBan.$ The
category $\FdOS$ has a zero object given by the zero-dimensional operator space
$0$. If $X,Y \in \FdOS$, then the vector space $X \oplus Y$ can be equipped
with a matrix norm via the linear isomorphism $\MM_n(X \oplus Y) \cong M_n(X) \linftyplus M_n(Y)$
and the Banach space structure of the latter space (see \secref{sub:banach}).
This gives an o.s.s. on $X \oplus Y$ and we write $X \linftyplus Y$ for the
resulting operator space, which is the categorical product of $X$ and $Y.$ The
notation is consistent with the one in $\FdBan$, i.e. $U(X \linftyplus Y) = UX
\linftyplus UY$ as Banach spaces. The categorical coproduct of $X$ and $Y$ in $\FdOS$, written $X \loneplus Y$,
is the vector space $X \oplus Y$ together with the o.s.s inherited via the linear isomorphism
$X \oplus Y \cong (X^* \linftyplus Y^*)^*$ and the o.s.s. of the latter space. This is also consistent with
the notation in $\FdBan$, i.e. $U(X \loneplus Y) = UX \loneplus UY$ as Banach spaces.
Note that $V(X \loneplus Y) = VX \oplus VY = V(X \linftyplus Y),$ i.e. $V$ sends
(co)products of $\FdOS$ to biproducts in $\FdVect.$ The functorial action of $\loneplus$ and $\linftyplus$
coincides with that of $\oplus$ on morphisms, i.e. $f \loneplus g = f \oplus g = f \linftyplus g.$

\subsection{Tensor Products of Operator Spaces}
We recall the construction of three operator space tensor products that we use in the sequel.
The first two are the completely projective tensor $\ptimes$ and the completely injective one $\itimes.$
They are analogous to their non-complete Banach space counterparts $\pBantimes$ and $\iBantimes$ in that
$\ptimes$ ($\itimes$) is the largest (smallest) tensor o.s.s that can be assigned on the
vector space tensor product $X \otimes Y$, but we elide the details of what this means.
The third tensor, called the Haagerup tensor, is rather unique to operator space theory and it allows us
to define von Neumann (co)algebras as certain kinds of internal (co)monoid objects in the sequel.
We define all three tensors for finite-dimensional operator spaces for simplicity.

\begin{remark}
  The completely projective tensor is usually referred to as the ``projective
  tensor'' of operator spaces and the completely injective one as the
  ``minimal'' or ``injective'' tensor of operator spaces. We have chosen the
  alternative names in order to avoid confusion with their Banach space
  counterparts.
\end{remark}

\begin{definition}[{\cite[(1.5.1)]{blecher-merdy}}]
  \label{def:itimes}
  Let $X,Y \in \FdOS.$ The \emph{completely injective tensor product},
  written $X \itimes Y$, is the vector space $X \otimes Y$
  with the o.s.s inherited via the linear isomorphism
  $X \otimes Y \cong \CB(Y^*,X)$ and the o.s.s of the latter space.
\end{definition}

We show in the sequel that the above tensor product serves as the
multiplicative disjunction and that the next one serves as the multiplicative
conjunction.

\begin{definition}
  \label{def:ptimes-simple}
  Let $X,Y \in \FdOS.$ The \emph{completely projective tensor product},
  written $X \ptimes Y$, is the vector space $X \otimes Y$
  with the o.s.s inherited via the linear isomorphism
  $X \otimes Y \cong (X^* \itimes Y^*)^*$ and the o.s.s of the latter space.
\end{definition}

The above definition is not the standard one in the literature, but we prove in Appendix \ref{app:fdos} that
for \emph{finite-dimensional} operator spaces the two are equivalent.
Before we may recall the third tensor, we need an auxiliary definition.

\begin{definition}[{\cite[\S 9.1]{er2000operator}}]
  Given $x\in \mathbb{M}_{n, r}(X)$ and $y\in \mathbb{M}_{r,m}(Y)$ the \emph{matrix inner product}
  of $x$ and $y$ is the matrix $x \odot y \in \MM_{n,m}(X \otimes Y)$ defined as
  \begin{equation*}
    (x \odot y)_{i,j} \eqdef \sum_{k=1}^r x_{i,k}\otimes y_{k,j} .
  \end{equation*}
\end{definition}

\begin{definition}[{\cite[\S 9.2]{er2000operator}}]
  Let $X,Y \in \FdOS.$ The \emph{Haagerup tensor product}, written $X \hagertimes Y$,
  is the vector space $X \otimes Y$ together with the o.s.s defined for
  $v \in \MM_n(X \otimes Y)$ by:
  \[ \|v\|_h \eqdef \inf\{\|x\|\|y\|\colon v=x\odot y, x\in M_{n,r}(X),y\in M_{r,n}(Y), r\in\mathbb N\}.\]
\end{definition}

An interesting property of the Haagerup tensor is that it is self-dual in $\FdOS$.

\begin{proposition}[{\cite[Cor. 9.4.8]{er2000operator}}]
  \label{prop:hager-self-dual}
  Given $X, Y \in \FdOS,$ the linear isomorphism
  $ \theta : X^* \hagertimes Y^* \cong (X \hagertimes Y)^* :: f\otimes g \mapsto (x\otimes y \mapsto f(x)g(y)) $
  is a c.i.i.
\end{proposition}

Another fact that distinguishes it (again) from the other two tensors is that the Haagerup tensor is \emph{not} symmetric in $\FdOS$.
What they do have in common is that
$ V(X \ptimes Y) = V(X \hagertimes Y) = V(X \itimes Y) = VX \otimes VY , $
i.e. the underlying vector space of all three tensors coincides with the vector
space tensor product.
Furthermore, the functorial action of all three tensors on morphisms coincides with $\otimes$, i.e.
$(f \ptimes g) = (f \hagertimes g) = (f \itimes g) = f \otimes g$.
The identity map (on $X \otimes Y$) gives complete contractions
$X \ptimes Y \to X \hagertimes Y \to X \itimes Y$, because the respective norms on the three spaces
satisfy $\norm{\cdot}_{\land} \geq \norm{\cdot}_{h} \geq \norm{\cdot}_{\lor}.$
The next proposition (special case of \cite[Theorem 6.1]{er-shuffle}) shows even more.

\begin{proposition}
  \label{prop:shuffle}
  Given operator spaces $A, B, C, D \in \FdOS$,
  the permutation
  \[ u_{ABCD} \colon (A \otimes B) \otimes (C \otimes D) \to (A \otimes C) \otimes (B \otimes D) :: (a \otimes b) \otimes (c \otimes d) \mapsto (a \otimes c) \otimes (b \otimes d) \]
  gives the two complete contractions $v_{ABCD} = u_{ABCD} = w_{ABCD}$ with types
  \begin{align*}
    w_{ABCD} &\colon (A \hagertimes B) \ptimes (C \hagertimes D) \to (A \ptimes C) \hagertimes (B \ptimes D) \\
    v_{ABCD} &\colon (A \itimes B) \hagertimes (C \itimes D) \to (A \hagertimes C) \itimes (B \hagertimes D)
  \end{align*}
\end{proposition}

In operator space theory, maps of the above form are known as \emph{shuffle
maps}, whereas in category theory they are known as \emph{interchange laws}.

\subsection{Conjugate Operator Space and Opposite Operator Space}

In the sequel, we need the operation $(\cdot)^*$ of taking the
conjugate transpose of a matrix. This is not a linear map, nor a complete
contraction $M_n \to M_n$ (when $n > 1$). We deal with this by recalling two involutive
constructions.

An anti-linear map
$f \colon V \to W$ between two complex vector spaces may be equivalently
described as a linear map $f \colon V_c \to W$, where $V_c$ is the
\emph{conjugate} vector space of $V$. Recall that $V_c$ is a complex vector space with the same underlying set as $V$
and same additive structure as $V,$ but whose scalar multiplication is defined by
$\lambda * v \eqdef \overline{\lambda}v$. This makes sense for operator spaces
as well.

\begin{definition}[{\cite[Sec. 2.9]{Pisier_2003}}]
  Let $X$ be an operator space. The conjugate vector space $X_c$
  can be equipped with an o.s.s. defined by
  $\| x \|_{M_n(X_c)} \eqdef \| x \|_{M_n(X)}$ for $x \in \MM_n(X_c).$
  We say that $X_c$ is the \emph{conjugate operator space} of $X.$
\end{definition}
\begin{definition}[{\cite[Sec. 2.10]{Pisier_2003}}]
  Let $X$ be an operator space. The \emph{opposite operator space} of $X$, written $X_o$, has the
  same underlying vector space but its o.s.s. is defined by
  $\| [x_{ji}] \|_{M_n(X_o)} \eqdef \| [x_{ij}] \|_{M_n(X)}$, for
  $x \in \MM_n(X_o)$.
\end{definition}

For convenience, we also sometimes write $X^o$ or $X^c$ for the opposite/conjugate space, e.g. we write $M_n^o$ instead of $(M_n)_o.$

\begin{example}
  The conjugate transpose is a completely isometric isomorphism $(\cdot)^* \colon M_n^c \cong M_n^o.$
  This is just a special case of a result for C*-algebras \cite[p. 65]{Pisier_2003}.
\end{example}

%%%%%%%%%%%%%%%%%%%%%%%%%%%%%%%%%%%%%%%%%%%%%%%%%%%%%%%%%%%%%%%%%%%%%%%%%%%%%%
\section{Categorical Properties of $\fdOS$}
\label{sec:fdos}
%%%%%%%%%%%%%%%%%%%%%%%%%%%%%%%%%%%%%%%%%%%%%%%%%%%%%%%%%%%%%%%%%%%%%%%%%%%%%%

We begin our contributions by proving relevant categorical properties of
$\FdOS.$ Our first proposition follows easily from existing results (see
Appendix \ref{app:fdos}). Our first theorem, given after it, shows the
relevance of $\FdOS$ to Linear Logic.

\begin{restatable}{proposition}{fdosMonoidal}
  \label{prop:fdos-monoidal}
  The category $\FdOS$ is a monoidal category when equipped with any of the tensors $\ptimes / \itimes / \hagertimes$ as monoidal product and
  $\mathbb C$ as monoidal unit. The first two monoidal structures are symmetric and in all cases the (symmetric) monoidal natural isomorphisms
  coincide with those in $\FdVect,$ i.e. the forgetful functor $V \colon \FdOS
  \to \FdVect$ is a \emph{strict} (symmetric) monoidal functor.
\end{restatable}
\begin{table}[t]
  \begin{center}
  \begin{tabular}{ |c|c|c|c| }
    \hline
    MALL / BV & $\FdOS$ & $\FdBan$ & \FdVect  \\
   \hline
    $X \multimap Y$ & $\CB(X,Y)$ & $B(X,Y)$ & $L(X,Y)$ \\
   \hline
    $X \aconj Y$ & $X \linftyplus Y$ & $X \linftyplus Y$ & $X \oplus Y$ \\
   \hline
    $X \oplus Y$ & $X \loneplus Y$ & $X \loneplus Y$ & $X \oplus Y$ \\
   \hline
    $X \otimes Y$ & $X \ptimes Y$ & $X \pBantimes Y$ & $X \otimes Y$ \\
   \hline
    $X \parr Y$ & $X \itimes Y$ & $X \iBantimes Y$ & $X \otimes Y$ \\
   \hline
    $X \triangleleft Y$ & $X \hagertimes Y$ &  & $X \otimes Y$ \\
   \hline
  \end{tabular}
  \end{center}
  \caption{MALL and BV categorical structure.}
  \label{tab:fdos-mall}
\end{table}
\begin{restatable}{theorem}{fdosMall}
  \label{thm:fdos-mall}
  The category $\fdOS$ has finite (co)products, it is $*$-autonomous with relevant data given in Table \ref{tab:fdos-mall}, and is therefore a model of MALL.
\end{restatable}

The Haagerup tensor is outside the scope of Linear Logic, but $\FdOS$ together
with this tensor and the maps from Proposition \ref{prop:shuffle} constitute a
\emph{BV-category (with negation)} in the sense of \cite{bv-category} (see Appendix
\ref{app:fdos}). The binary connective $\triangleleft$
from Table \ref{tab:fdos-mall} stands for the ``seq'' tensor of BV-categories,
which corresponds to $\hagertimes$ in $\FdOS.$
In the sequel, we show that $\hagertimes$ has even more categorical properties.

Next, we consider some of the categorical properties of the assignments
$(\cdot)_c$ and $(\cdot)_o$ of taking conjugate/opposite operator spaces. The
second equality from our next proposition follows from {\cite[Lemma
2.1]{choi20}} where it is shown that the identity map is a c.i.i (and therefore
an equality). The first equality is easy.

\begin{restatable}{proposition}{opConjugateCB}
  \label{prop:op-conjugate-cb}
  Given operator spaces $X$ and $Y$, we have operator space equalities
  $\CB(X,Y)_c =  \CB(X_c, Y_c)$ and $\CB(X,Y)_o = \CB(X_o, Y_o)$.
\end{restatable}

Since $\mathbb C = \mathbb C_o$ and $\mathbb C_c \cong \mathbb C$, with conjugation serving as the completely isometric isomorphism,
we have as a special case that $(X_c)^* \cong (X^*)_c$ and
$(X_o)^* = (X^*)_o$. Moreover, we can see the assignments
$(\cdot)_c \colon \FdOS \to \FdOS$ and $(\cdot)_o \colon \FdOS \to \FdOS$
as (covariant) functors by simply defining $f_o \eqdef f$ and $f_c \eqdef f.$
It is easy to see that
$ (\cdot)_c \circ (\cdot)_c = \mathrm{Id} = (\cdot)_o \circ (\cdot)_o \text{ and }
   (\cdot)_c \circ (\cdot)_o = (\cdot)_o \circ (\cdot)_c  $
and therefore we have a pair of commuting functorial involutions on $\FdOS$
which are (set-theoretic) identities on morphisms, but not on objects.
Even more, these functors behave well with other functorial constructions we have seen.

\begin{restatable}{proposition}{functorialInvolutions}
  Given operator spaces $X, Y \in \FdOS$, we have:
  \begin{multicols}{3}
    \begin{itemize}
      \item $(X \ptimes Y)_c = X_c \ptimes Y_c$
      \item $(X \itimes Y)_c = X_c \itimes Y_c$
      \item $(X \ptimes Y)_o = X_o \ptimes Y_o$
      \item $(X \itimes Y)_o = X_o \itimes Y_o$
    \end{itemize}

    \begin{itemize}
      \item $(X \loneplus Y)_c = X_c \loneplus Y_c$
      \item $(X \linftyplus Y)_c = X_c \linftyplus Y_c$
      \item $(X \loneplus Y)_o = X_o \loneplus Y_o$
      \item $(X \linftyplus Y)_o = X_o \linftyplus Y_o$
    \end{itemize}

    \begin{itemize}
      \item $(X \hagertimes Y)_c = X_c \hagertimes Y_c$
      \item a c.i.i., called $\gamma$, $(X \hagertimes Y)_o \cong Y_o \hagertimes X_o$
        $\gamma(x \otimes y) \eqdef (y \otimes x)$.
    \end{itemize}
  \end{multicols}
\end{restatable}
\begin{proof}
  Some of these facts are already known and the rest are in Appendix \ref{app:fdos}.
\end{proof}

Therefore, even though the Haagerup tensor is \emph{not} symmetric in $\FdOS$,
we can still talk about the usual symmetry $\gamma$, provided we use the involution $(\cdot)_o.$
We also use $\gamma$ to denote its own inverse, i.e. $\gamma: Y_o \hagertimes X_o \cong (X \hagertimes Y)_o$,
since both $\gamma$ and $\gamma^{-1}$ are defined in the same way, namely, via swapping.

We also briefly use the category $\FdOScb$ which has f.d. operator spaces as
objects, but whose morphisms are the c.b. maps between them. All the
categorical properties from this section hold verbatim for $\FdOScb$ and the
subcategory inclusion $\FdOS \hookrightarrow \FdOScb$ preserves this data.
However, as we already explained, the linear and c.b. maps coincide for f.d.
operator spaces, so it follows that the forgetful functor $ V \colon\ \FdOScb
\to \FdVect$ gives an equivalence of categories $\FdOScb \simeq \FdVect.$
Because of this, $\FdOScb$ is a less interesting, and more degenerate (from a
MALL perspective) category compared to $\FdOS$, e.g. the three tensors
$\ptimes$, $\hagertimes$, $\itimes$ are naturally isomorphic to each other, via the identity map, in
$\FdOScb$ and each of them gives a compact closed structure on $\FdOScb.$ Of
course, this is not true in $\FdOS,$ which is the primary category of interest.

%%%%%%%%%%%%%%%%%%%%%%%%%%%%%%%%%%%%%%%%%%%%%%%%%%%%%%%%%%%%%%%%%%%%%%%%%%%%%%
\section{Novel Categorical Approach to von Neumann Algebras}
\label{sec:vNAlg}
%%%%%%%%%%%%%%%%%%%%%%%%%%%%%%%%%%%%%%%%%%%%%%%%%%%%%%%%%%%%%%%%%%%%%%%%%%%%%%

Operator spaces alone do not have sufficient structure for quantum computation.
In the Heisenberg picture, we can use \emph{von Neumann algebras (vN-algebras)}
to define the relevant quantum operations. Every such algebra can be equipped
with a canonical o.s.s., so we may think of them as operator spaces with
additional structure. In this section, we consider f.d. vN-algebras and show
that they may be equivalently defined as certain kinds of involutive monoid
objects in $(\FdOS, \mathbb C, \hagertimes).$ This serves two purposes: (1) it
can help readers who are familiar with category theory; (2) it makes it
easy to dualise the definition and formulate \emph{von Neumann coalgebras} (\S
\ref{sec:vNCoalg}) which we use for the \Schrod{} picture.

We are only concerned with \emph{finite-dimensional} vN-algebras. They coincide
with the finite-dimensional C*-algebras and our next definition is standard \cite{takesaki1,blackadar}.

\begin{definition}
  \label{def:vN-alg-standard}
  A f.d. \emph{vN-algebra} is a f.d. complex vector space $A$ together with:
  \begin{enumerate}
    \item a bilinear binary operation, called \emph{multiplication}, and written
  via juxtaposition, which is associative, i.e. $a(bc) = (ab)c$;
\item an element $1 \in A$, called \emph{unit}, such that $1a = a = a1;$
\item an anti-linear unary operation $(\cdot)^*$, called \emph{involution}, such that
  $(a^*)^* = a$ and $(ab)^* = b^*a^*;$
\item a (necessarily complete) norm which is \emph{submultiplicative}, i.e. $\norm{ab} \leq \norm{a} \norm{b}$,
  and which satisfies the \emph{C*-identity}, i.e. $\norm{a^*a} = \norm{a}^2.$
  \end{enumerate}
\end{definition}

If $A$ is a f.d. vN-algebra, then so is $\MM_n(A)$ -- the multiplication, unit and
involution are defined in the obvious way and then there is a \emph{unique} way
to define a norm on $\MM_n(A)$ such that it becomes a vN-algebra \cite[1.2.3]{blecher-merdy}. These norms
give the canonical o.s.s of $A$ \cite[1.2.3]{blecher-merdy} and henceforth we view vN-algebras as operator
spaces with additional structure.

\begin{example}
  \label{ex:fd-vN-algebra}
  The operator space $M_n$ is a vN-algebra with unit given by the identity
  matrix, multiplication by matrix multiplication, and involution given by
  $(\cdot)^*.$ The operator space $\linftyplus_{1 \leq i \leq n} M_{k_i}$
  is also a vN-algebra with unit, multiplication, and involution defined
  pointwise. In fact, modulo a vN-algebra isomorphism (defined below),
  every f.d. vN-algebra is of this form \cite[p. 50]{takesaki1}.
\end{example}

For our categorical definition of a f.d. vN-algebra, we take inspiration
from \emph{reversing involutive monoids} \cite{Jacobs_2011} and \emph{star
algebras} in a bar category \cite{bm2007}. More importantly, the
(f.d.) unital operator algebras have been characterised as the monoid objects
w.r.t. the Haagerup tensor in $\mathbf{(Fd)OS}$
\cite[Theorem 6.1]{Pisier_2003}, \cite{operator-algebras-haagerup}. Since f.d. vN-algebras are special kinds of unital operator
algebras, it is natural to consider them as special kinds of monoid objects in
$(\FdOS, \mathbb C, \hagertimes).$ We do so in our subsequent definition.

\begin{definition}\label{def:vN-alg-cat}
  A f.d. \emph{vN-algebra} is a f.d. operator space $A$ together with
  \begin{enumerate}
    \item a complete contraction $\mu \colon A \hagertimes A \to A$, called \emph{multiplication};
    \item a complete contraction $\eta \colon \mathbb C \to A$, called \emph{unit};
    \item a complete contraction $i \colon A_c \to A_o,$ called \emph{involution},
  \end{enumerate}
  such that $(A, \eta, \mu)$ is a monoid object in $(\FdOS, \mathbb C, \hagertimes)$, the following diagrams
    \begin{equation}
      \adjustbox{scale=0.8}{
        \begin{tikzcd}[ampersand replacement= \&]
          A \arrow[equal]{d} \arrow[equal]{rr}           \&           \& A      \& A_c \hagertimes A_c \arrow[d, "i\otimes i"'] \arrow[equal]{r} \& (A\hagertimes A)_c \arrow[r, "\mu"] \& A_c \arrow[d, "i"] \\
          A_{cc} \arrow[r, "i_c"'] \& A_{oc} = A_{co} \arrow[r, "i_o"'] \& A_{oo} \arrow[equal]{u} \& A_o \hagertimes A_o \arrow[r, "{\gamma}"']              \& (A\hagertimes A)_o \arrow[r, "\mu"'] \& A_o
        \end{tikzcd}}
      \end{equation}
    commute, and for every complete isometry (i.e. strong mono) $a: \mathbb{C} \to A$, the following composite in $\FdOS$ is also a
    complete isometry (i.e. strong mono)
    \begin{equation}
      \mathbb{C} \cong \mathbb{C}_o \hagertimes \mathbb{C}_c
      \xrightarrow{a \otimes a} A_o \hagertimes A_c
      \xrightarrow{A_o \otimes i} A_o \hagertimes A_o
      \xrightarrow{\gamma} (A \hagertimes A)_o
      \xrightarrow{\mu} A_o.
    \end{equation}
\end{definition}
\begin{restatable}{theorem}{vNAlgEquiv}
  Definitions \ref{def:vN-alg-standard} and \ref{def:vN-alg-cat} are equivalent, i.e. the data in one definition uniquely determines the data in the other.
\end{restatable}

Our theorem is proven in Appendix \ref{app:vN-alg} and its proof shows that the
two (equivalent) notions may be used interchangeably. We do so in the sequel.

\begin{definition}
  \label{def:vNalg}
  Let $\vNalg$ denote the subcategory of $\fdOS$ whose objects
  are f.d. vN-algebras and whose morphisms $\varphi \colon A \to B$
  are \emph{multiplicative, unital,} and \emph{involutive}, i.e. the following diagrams in $\FdOS$ commute\footnote{It suffices to assume $\varphi$ is linear and the three diagrams imply it is also c.c.}.
  \begin{equation*}
    \adjustbox{scale=0.8}{
      \begin{tikzcd}
        A \hagertimes A \arrow[r, "\varphi \otimes \varphi"] \arrow[d, "\mu_A"'] & B \hagertimes B \arrow[d, "\mu_B"]
        & & \mathbb{C} \arrow[ld, "\eta_A"'] \arrow[rd, "\eta_B"] &
        & A_c \arrow[r, "\varphi_c "] \arrow[d, "i_A"'] & B_c \arrow[d, "i_B"] \\
        A \arrow[r, "\varphi"'] & B
        & A \arrow[rr, "\varphi"'] & & B
        & A_o \arrow[r, "\varphi_o"'] & B_o \\
    \end{tikzcd}}
  \end{equation*}
\end{definition}

In the literature, these morphisms are known as ``\emph{unital $*$-homomorphisms}'',
but we simply call them vN-algebra morphisms. It is
already known that (possibly infinite-dimensional) vN-algebras form a symmetric
monoidal category with small products \cite{kornell-vn}, so the next
propositions follow easily. However, we present them in our categorical style, as
this helps for the development in \S \ref{sec:vNCoalg}.

\begin{restatable}{proposition}{vNtens}\label{prop:tensor-alg}
  If $A, B \in \vNAlg$, then
  the operator space $A \itimes B$ with
  \begin{itemize}
    \item unit $\eta \eqdef \left( \mathbb C \cong \mathbb C \itimes \mathbb C \xrightarrow{\eta_A \otimes \eta_B} A \itimes B \right)$
    \item involution  $i \eqdef \left( (A \itimes B)_c = A_c \itimes B_c \xrightarrow{i_A \otimes i_B} A_o \itimes B_o = (A \itimes B)_o \ \right)$
    \item mult. $\mu \eqdef \left((A\injtens B) \hagertimes (A\injtens B) \xrightarrow{v} (A \hagertimes A) \injtens (B \hagertimes B) \xrightarrow{\mu_A \otimes \mu_B} A \injtens B \right) $
  \end{itemize}
  is also a vN-algebra. Also, the operator space $A \linftyplus B$ is a vN-algebra with
  \begin{itemize}
    \item unit $\eta \eqdef \left(\mathbb C \xrightarrow{\langle \text{id}, \text{id} \rangle} \mathbb{C} \linftyplus \mathbb{C} \xrightarrow{\eta_A \oplus \eta_B} A \linftyplus B\right)$
    \item involution  $i \eqdef \left( (A \linftyplus B)_c = A_c \linftyplus B_c \xrightarrow{i_A \oplus i_B} A_o \linftyplus B_o = (A \linftyplus B)_o\right)$
    \item mult. $\mu \eqdef \left( (A \linftyplus B)\hagertimes (A \linftyplus B)
    \xrightarrow{v'} (A \hagertimes A) \linftyplus (B \hagertimes B) \xrightarrow{\mu_A \oplus \mu_B} A \linftyplus B \right)$
  \end{itemize}
  where $v' \eqdef \left( (A \linftyplus B)\hagertimes (A \linftyplus B)
  \xrightarrow{\langle (\pi_A \otimes \pi_A), (\pi_B \otimes \pi_B) \rangle} (A \hagertimes A) \linftyplus (B \hagertimes B)\right) $
\end{restatable}

\begin{restatable}{proposition}{vNsmc}
  \label{prop:vN-monoidal}
  The category $\vNalg$ is symmetric monoidal and has finite products with
  the action of $\itimes$ and $\linftyplus$ extended as in Proposition \ref{prop:tensor-alg}.
\end{restatable}

\paragraph{Quantum Operations.} The quantum operations in the Heisenberg picture of f.d. quantum theory are given by the completely positive unital (CPU) maps.
We already defined unital maps, so now we recall completely positive maps.
Given $A \in \vNAlg$, we say that $p \in A$ is a \emph{positive
element} if there exists $a \in A$, such that $p = a^*a$. An equivalent
way of saying this is that the diagram in $\FdOScb$
\begin{equation}
  \label{eq:positive-element}
\begin{tikzcd}
  \mathbb{C}  \arrow[equal]{d} \arrow[r, "\cong"]   & \mathbb{C}_o \hagertimes \mathbb{C}_c \arrow[r, "a \otimes a"] & A_o \hagertimes A_c \arrow[r, "A_o \otimes i"] & A_o \hagertimes A_o  \arrow[r, "\gamma"] & (A \hagertimes A)_o \arrow[d, "\mu"] \\
  \mathbb{C}_o \arrow[rrrr, "p"'] &                                                                    &                                               &                                          & A_o
\end{tikzcd}
\end{equation}
commutes, where we have made the obvious identification of elements of $A$ as morphisms in $\FdOScb.$
A linear map
$\varphi \colon A \to B$ between two f.d. vN-algebras is called \emph{positive} if it
preserves positive elements and \emph{completely positive (c.p.)} if
$M_n \itimes \varphi \colon M_n \itimes A \to M_n \itimes B$ is positive for
every $n \in \mathbb N.$

Finally, we write $\HH$ for the category whose objects are f.d. vN-algebras and
whose morphisms are the CPU maps. The next proposition, which is a special
case of \cite[Corollary 5.1.2]{er2000operator}, shows that $\HH$ is a subcategory of
$\FdOS.$ The proposition after it summarises the categorical properties of
$\HH$ relevant to MALL.

\begin{restatable}{proposition}{eqCPCC}
  \label{prop:alg-cc-cp}
  Let $\varphi \colon A_1 \to A_2$ be a linear unital map between two von Neumann algebras. Then, $\varphi$ is completely positive iff $\varphi$ is completely contractive.
\end{restatable}

\begin{restatable}{proposition}{Hsmcc}
  \label{prop:Hsmcc}
  The category $\HH$ is symmetric monoidal (w.r.t. $\itimes$), has finite products (given by $\linftyplus$), and the inclusions
  $\vNAlg \hookrightarrow \HH \hookrightarrow \FdOS$ preserve this data.
\end{restatable}

%%%%%%%%%%%%%%%%%%%%%%%%%%%%%%%%%%%%%%%%%%%%%%%%%%%%%%%%%%%%%%%%%%%%%%%%%%%%%%
\section{von Neumann Coalgebras}
\label{sec:vNCoalg}
%%%%%%%%%%%%%%%%%%%%%%%%%%%%%%%%%%%%%%%%%%%%%%%%%%%%%%%%%%%%%%%%%%%%%%%%%%%%%%

In \cite{er-shuffle}, a vN-coalgebra structure is \emph{defined} as the
predual of the structure of a vN-algebra. Instead of doing this, we introduce f.d. vN-coalgebras
\emph{independently} of vN-algebras and we \emph{derive} the duality later.
Thanks to the preparatory categorical work we did in Sections \ref{sec:fdos}--\ref{sec:vNAlg}, and
keeping in mind that the Haagerup tensor is self-dual
(\cref{prop:hager-self-dual}), we can now easily dualise the
definition of a f.d. vN-algebra to obtain that of a f.d. vN-coalgebra.

\begin{definition}\label{cat.von-neum-co}
  A f.d. \emph{vN-coalgebra} is a f.d. operator space $C$ together with
  \begin{enumerate}
    \item a complete contraction $\delta \colon C \to C \hagertimes C$, called \emph{comultiplication};
    \item a complete contraction $\varepsilon \colon C \to \mathbb C$, called \emph{counit};
    \item a complete contraction $j \colon C_o \to C_c,$ called \emph{involution},
  \end{enumerate}
  such that $(C, \varepsilon, \delta)$ is a comonoid object in $(\FdOS, \mathbb C, \hagertimes)$,
  the following diagrams
      \begin{equation}
        \adjustbox{scale=0.8}{
          \begin{tikzcd}
            C \arrow[equal]{rr}                                           &                                                                                   & C         &    C_o \arrow[d, "\delta"'] \arrow[r, "j"]                 & C_c \arrow[r, "\delta"]                                     & (C\hagertimes C)_c \arrow[equal]{d} \\
            C_{oo} \arrow[equal]{u} \arrow[r, "j_o"'] & C_{co} = C_{oc} \arrow[r, "j_c"'] & C_{cc} \arrow[equal]{u} & (C\hagertimes C)_o \arrow[r, "{\gamma}"'] & C_o\hagertimes C_o \arrow[r, "j\otimes j"'] & C_c\hagertimes C_c
          \end{tikzcd}
        }
      \end{equation}
    commute, and for every complete quotient map (i.e. strong epi) $e: C \to \mathbb{C}$, the following composite in $\FdOS$ is
    a complete quotient map (i.e. strong epi)
    \begin{equation}
      C_o \xrightarrow{\delta}
      (C\hagertimes C)_o \xrightarrow{\gamma}
      C_o \hagertimes C_o \xrightarrow{C_o \otimes j}
      C_o\hagertimes C_c \xrightarrow{e \otimes e}
      \mathbb{C}_o \hagertimes \mathbb{C}_c \cong \mathbb{C} .
    \end{equation}
\end{definition}

We dualise \cref{def:vNalg} in order to obtain the next one.

\begin{definition}
  Let $\vNCoalg$ denote the subcategory of $\fdOS$ whose objects
  are f.d. vN-coalgebras and whose morphisms $\varphi \colon C \to D$
  are \emph{comultiplicative, counital,} and \emph{involutive}, i.e. the following diagrams in $\FdOS$ commute.
  \begin{equation*}
    \adjustbox{scale=0.8}{
      \begin{tikzcd}
        C \arrow[r, "\varphi"] \arrow[d, "\delta_C"'] & D \arrow[d, "\delta_D"]
        & C \arrow[rd, "\varepsilon_C"'] \arrow[rr, "\varphi"] & & D \arrow[ld, "\varepsilon_D"]
        & C_o \arrow[r, "\varphi_o"] \arrow[d, "j_C"'] & D_o \arrow[d, "j_D"] \\
        C \hagertimes C \arrow[r, "\varphi \otimes \varphi"'] & D \hagertimes D
        & & \mathbb{C} &
        & C_c \arrow[r, "\varphi_c"'] & D_c \\
    \end{tikzcd}
  }
  \end{equation*}
\end{definition}

We also dualise the construction of tensors and direct sums from Proposition
\ref{prop:tensor-alg} and we now use the remaining shuffle map from Proposition
\ref{prop:shuffle}.
\begin{restatable}{proposition}{vNcoalgtens}
  \label{prop:tensor-coalg}
  If $C,D \in \vNcoalg$, then the operator space $C \ptimes D$ with
  \begin{itemize}
    \item counit $ \varepsilon \eqdef \left( C \ptimes D \xrightarrow{\varepsilon_C \otimes \varepsilon_D} \mathbb C \ptimes \mathbb C \cong \mathbb C \right)$
    \item involution $j \eqdef \left( (C \ptimes D)_o = C_o \ptimes D_o \xrightarrow{j_C \otimes j_D} C_c \ptimes D_c = (C \ptimes D)_c \right)$
    \item comult. $ \delta \eqdef \left(C \projtens D \xrightarrow{\delta_C \otimes \delta_D} (C \hagertimes C) \projtens (D \hagertimes D) \xrightarrow{w} (C\projtens D) \hagertimes (C\projtens D)\right) $
  \end{itemize}
  is also a vN-coalgebra. Also, the operator space $C \loneplus D$ is a vN-coalgebra with
  \begin{itemize}
    \item counit $\varepsilon \eqdef \left( C \loneplus D \xrightarrow{\varepsilon_C \oplus \varepsilon_D}\mathbb{C}\loneplus \mathbb{C} \xrightarrow{[ \text{id}, \text{id}]} \mathbb{C}\right) $
    \item involution  $j \eqdef \left( (C \loneplus D)_o = C_o \loneplus D_o \xrightarrow{j_C \oplus j_D} C_c \loneplus D_c = (C \loneplus D)_c \right)$
    \item comult. $\delta \eqdef \left(C \loneplus D \xrightarrow{\delta_C \oplus \delta_D} (C \hagertimes C) \loneplus (D \hagertimes D)
    \xrightarrow{w'} (C \loneplus D) \hagertimes (C \loneplus D)\right)$
  \end{itemize}
  where $w' \eqdef \left( (C \hagertimes C) \loneplus (D \hagertimes D) \xrightarrow{[(\iota_C \otimes \iota_C) , (\iota_D \otimes \iota_D)]} (C \loneplus D) \hagertimes (C \loneplus D)\right)$
\end{restatable}

\begin{restatable}{proposition}{vNcoalgsmcc}
  \label{prop:vNcoalgsmcc}
  The category $\vNcoalg$ is symmetric monoidal and has finite coproducts with
  the action of $\ptimes$ and $\loneplus$ extended as in Proposition \ref{prop:tensor-coalg}.
\end{restatable}

\begin{restatable}{example}{TnCoalgebra}
  \label{ex:coalg}
  The operator space $T_n$ is a vN-coalgebra with:
    (1) counit given by the trace functional $\varepsilon \eqdef \trace{} \colon T_n \to \mathbb C :: t \mapsto \trace(t);$
    (2) involution given by the conjugate transpose $j \eqdef (\cdot)^* \colon T_n^o \to T_n^c ; $
    (3) comultiplication given by $\delta \colon T_n \to T_n \hagertimes T_n :: e_{ij} \mapsto \sum_k  (e_{kj} \otimes e_{ik}),$
  where $e_{ij} \in T_n$ is the matrix that has $1$ in position $(i,j)$ and $0$ elsewhere. Using Proposition \ref{prop:tensor-coalg},
  we see that $\loneplus_{1 \leq i \leq n} T_{k_i}$ is also a vN-coalgebra and using Theorem \ref{thm:alg-coalg-dual} and Example \ref{ex:fd-vN-algebra} it follows
  that each f.d. vN-coalgebra is of this form, modulo vN-coalgebra isomorphism.
\end{restatable}

\paragraph{Quantum Operations.} We use f.d. vN-coalgebras to define the quantum operations in the
\Schrod{} picture, i.e. the CPTP maps. This can be done by dualising the CPU maps from \S \ref{sec:vNAlg}.
The dual of a unital map is obviously a counital one. Unfortunately, in quantum information theory,
the term ``trace-preserving'' is standard, so we simply say that a linear map $\varphi \colon
C \to D$ between two vN-coalgebras is \emph{trace-preserving} if it is
counital. We explain why this is justified. Note that if $C = T_n$ and $D =
T_m$, then from Example \ref{ex:coalg}, we see that a counital map is precisely
a trace-preserving one. More generally, the direct sum $\loneplus_{i} T_{k_i}$
can be linearly embedded block-diagonally in a sufficiently large matrix space
$\MM_{m}$ and then counitality is equivalent to trace-preservation with
respect to this view. Using Example \ref{ex:coalg}, this
covers all f.d. vN-coalgebras, modulo isomorphism.

We proceed with the (dual) definition of a completely positive map between two f.d. vN-coalgebras.
We work again in the category $\FdOScb.$ Given a vN-coalgebra $C$, we say that a linear functional $p \colon C \to \mathbb C$
is \emph{positive} if there exists a linear functional $a \colon C \to \mathbb C$, such that the following diagram:
\begin{equation}
  \label{eq:positive-functional}
\begin{tikzcd}
  (C\hagertimes C)_o \arrow[r, "\gamma"]    & C_o\hagertimes C_o \arrow[r, "C_o\otimes j"] & C_o \hagertimes C_c \arrow[r, "a \otimes a"] & \mathbb{C}_o \hagertimes \mathbb{C}_c \arrow[r, "\cong"] & \mathbb{C} \arrow[equal]{d} \\
  C_o \arrow[u, "\delta"] \arrow[rrrr, "p"'] &                                              &                                                  &                                                         & \mathbb{C}_o
\end{tikzcd}
\end{equation}
commutes. Note that \eqref{eq:positive-functional} is dual to \eqref{eq:positive-element}. Continuing in this fashion, we say that
a linear map $\varphi \colon C \to D$ between vN-coalgebras is \emph{positive} if it reflects positive linear functionals, i.e.
$p \circ \varphi \colon C \to \mathbb C$ is a positive linear functional on $C$ for every positive linear functional $p \colon D \to \mathbb C.$
Finally, we say $\varphi$ is \emph{completely positive} if $T_n \ptimes \varphi \colon T_n \ptimes C \to T_n \ptimes D$ is positive for every
$n \in \mathbb N.$
This notion of (complete) positivity is equivalent to the usual/concrete one (see Appendix \ref{app:vn-coalg-cp}).

We write $\SSS$ for the category whose objects are f.d. vN-coalgebras
and whose morphisms are the CPTP maps between them. Finally, we show that
$\SSS$ is a subcategory of $\FdOS$ and summarise its structure relevant to
MALL.

\begin{proposition}
  \label{prop:coalg-cc-cp}
  Let $\varphi \colon C_1 \to C_2$ be a linear trace-preserving map between f.d. vN-coalgebras. Then, $\varphi$ is completely positive iff $\varphi$ is completely contractive.
\end{proposition}

\begin{restatable}{proposition}{Ssmcc}
  The category $\SSS$ is symmetric monoidal (w.r.t. $\ptimes$), has finite coproducts (given by $\loneplus$), and the inclusions
  $\vNCoalg \hookrightarrow \SSS \hookrightarrow \FdOS$ preserve this data.
\end{restatable}

%%%%%%%%%%%%%%%%%%%%%%%%%%%%%%%%%%%%%%%%%%%%%%%%%%%%%%%%%%%%%%%%%%%%%%%%%%%%%%
\section{Heisenberg-\Schrod{} Duality and vN-(co)algebras}
\label{sec:hs}
%%%%%%%%%%%%%%%%%%%%%%%%%%%%%%%%%%%%%%%%%%%%%%%%%%%%%%%%%%%%%%%%%%%%%%%%%%%%%%

In f.d. quantum theory, the quantum operations in the Heisenberg picture
are precisely the CPU maps, i.e. the morphisms of $\HH$, whereas the quantum
operations in the \Schrod{} picture are precisely the CPTP maps, i.e. the
morphisms of $\SSS.$ Both $\HH$ and $\SSS$ are subcategories of $\FdOS$ and the
dual functor $(\cdot)^* \colon \FdOS \to \FdOS$ can be defined also on
vN-(co)algebras, as we show next.

\begin{restatable}{proposition}{vNdual}
  \label{prop:duals}
  If $A \in \vNAlg$, then the o.s. dual $A^*$ is a vN-coalgebra with
  counit $\varepsilon \eqdef \left( A^* \xrightarrow{\eta^*} \mathbb C^* \cong \mathbb C\right) , $
  involution $j \eqdef \left( (A^*)_o = (A_o)^* \xrightarrow{i^*} (A_c)^* \cong (A^*)_c  \right) $, and
  comultiplication $\delta \eqdef \left(A^* \xrightarrow{\mu^*} (A \hagertimes A)^* \cong A^* \hagertimes A^*\right)$.
  Conversely, if $C \in \vNCoalg$, then the o.s. dual $C^*$ is a vN-algebra with unit
  $\eta \eqdef \left(\mathbb C \cong \mathbb C^* \xrightarrow{\varepsilon^*} C^* \right) , $
  involution $i \eqdef \left( (C^*)_c \cong (C_c)^* \xrightarrow{j^*} (C_o)^* = (C^*)_o \right) $,
  and multiplication given by $\mu \eqdef \left( C^* \hagertimes C^* \cong (C \hagertimes C)^* \xrightarrow{\delta^*} C^* \right).$
\end{restatable}

\begin{restatable}{theorem}{algcoalgdual}
  \label{thm:alg-coalg-dual}
  We have equivalences of categories $(\cdot)^* \colon \SSS \simeq \HH^\mathrm{op} \colon (\cdot)^*$
  and also $(\cdot)^* \colon \vNcoalg \simeq \vNalg^\mathrm{op} \colon (\cdot)^*$,
  where the action on objects of $(\cdot)^*$ is defined as in \cref{prop:duals} and on morphisms
  in the usual way, i.e. as in $\FdOS.$ Moreover, these equivalences are strong monoidal (and (co)product preserving).
\end{restatable}

\begin{example}
  The c.i.i. from \cref{ex:Tn} is a vN-coalgebra isomorpism $T_n \cong M_n^*$, which
  determines a vN-algebra isomorphism $T_n^* \cong M_n$ by duality.
\end{example}

Theorem \ref{thm:alg-coalg-dual} gives a categorical formulation of the Heisenberg-\Schrod{}
duality and shows precisely how the MALL related structure of $\vNCoalg$ and
$\SSS$ from \S \ref{sec:vNCoalg} is dual to that of $\vNAlg$ and $\HH$ from \S
\ref{sec:vNAlg}, and vice-versa. The duality $\vNcoalg \simeq \vNalg^{\mathrm{op}}$ shows that
vN-(co)algebras are dual to each other in a strong mathematical sense, whereas
the duality $\SSS \simeq \HH^{\mathrm{op}}$ is more relevant
for quantum computation and we use it in the next section.

%%%%%%%%%%%%%%%%%%%%%%%%%%%%%%%%%%%%%%%%%%%%%%%%%%%%%%%%%%%%%%%%%%%%%%%%%%%%%%
\section{Revisiting Quantum Coherence Spaces}
\label{sec:qcs}
%%%%%%%%%%%%%%%%%%%%%%%%%%%%%%%%%%%%%%%%%%%%%%%%%%%%%%%%%%%%%%%%%%%%%%%%%%%%%%

The preparatory categorical work in Sections \ref{sec:fdos} -- \ref{sec:hs} shows
that the (sub)categories $\SSS \hookrightarrow \FdOS \hookleftarrow \HH$ have
the right categorical structure for a model of MALL whose duality is induced
by the Heisenberg-\Schrod{} duality. However, there is one problem: neither
$\SSS$, nor $\HH$, is a \emph{full} subcategory of $\FdOS.$
This means that we have morphisms between vN-algebras that are not
CPU maps and morphisms between vN-coalgebras that are not CPTP maps in $\fdOS$.
For example, the map $a \mapsto -a:M_2 \to M_2$ is not positive but it is in
$\fdOS(M_2, M_2)$.
Luckily, there
is a simple solution: a semantic technique of Hyland and Schalk \cite{HS},
based on \emph{gluing and orthogonality}, allows us to carve out a category
$\QQ$ from $\FdOS,$ such that we get fully faithful inclusions
$\SSS \hookrightarrow \QQ \hookleftarrow \HH$, ensuring that $\QQ$ contains
precisely the quantum operations in the relevant homsets, while still preserving all the
MALL structure. The main idea, in a nutshell, is to ensure trace-preservation
(unitality) for morphisms between vN-coalgebras (vN-algebras) and then fullness
is guaranteed by Propositions \ref{prop:alg-cc-cp} and \ref{prop:coalg-cc-cp}.
This is precisely what we do by defining a suitable notion of orthogonality/polarity
based on ideas from \cite{HS}. The resulting model has an obvious
resemblance to (probabilistic) coherence spaces.

More specifically, we are interested in pairs $(X, S)$, where $X$ is
an operator space and $S \subseteq \Ball(X)$ is a subset of the unit ball of
$X$, which satisfies an additional condition that we proceed to explain.
Our next definition has similarities with \emph{polar
sets} from functional and convex analysis, which motivates its name.

\begin{definition}[Polar]
  Given an operator space $X$ and a subset $S\subseteq \Ball(X),$
  we define the \emph{polar} of $S$ to be
  $ S^\circ \eqdef \{f \in \Ball(X^*) \ | \ \forall s \in S. \ f(s) = 1 \} . $
  We say that $S$ is \emph{bipolar} if $S=d^{-1}[S^{\circ\circ}]$,
  where $d \colon X \cong X^{**}$ is the canonical c.i.i.
\end{definition}

To avoid notational overhead, when dealing with multiple applications of the polar construction,
we implicitly consider the image of such polar sets to be taken under $d^{-1}$. With this
convention, we simply say that $S$ is bipolar if $S = S^{\circ\circ}.$
As expected, the polar introduces a Galois connection \cite[\S 5.1]{HS} and we have:
for $R \subseteq S\subseteq \Ball(X)$, it follows that
$S^\circ \subseteq R^\circ$, $S\subseteq S^{\circ\circ},$ and $S^{\circ\circ\circ} = S^{\circ}$.
We can now introduce the main category that is computationally relevant.
\begin{definition}
  Let $\QQ$ be the category whose objects are pairs $(X,S)$ with $X \in \FdOS$
  and $S \subseteq \Ball(X)$ a bipolar set, i.e. $S = S^{\bipolar},$
  and whose morphisms $f: (X, S) \to (Y, R)$ are complete contractions
  $f: X \to Y$ such that $f[S] \subseteq R$.
\end{definition}

Our notion of polarity is instead called ``orthogonality'' in \cite{HS} and
it enjoys useful properties, e.g. it is \emph{tight}, \emph{stable}, \emph{focused},
\emph{precise}, in the sense of \cite{HS}. By using results from
\cite{HS}, it follows that the MALL structure of $\FdOS$ is preserved.

\begin{restatable}{theorem}{qMall}
\label{thm:q-mall}
  The category $\QQ$ is $*$-autonomous and has finite (co)products:
  duals are given by $(X,S)^* \eqdef (X^*, S^\circ);$
  the initial object by $(0, \varnothing)$;
  the terminal object by $(0, \{0\});$
  binary products by $(X,S) \linftyplus (Y, R) \eqdef (X \linftyplus Y, S \times R)$;
  binary coproducts by $(X,S) \loneplus (Y,R) \eqdef (X \loneplus Y, (S^\circ + R^\circ)^\circ),$
  where
  \[ S^\circ + R^\circ \eqdef \{ [ s', r' ] \colon X \loneplus Y \to \mathbb{C} \ | \ s' \in S^\circ, r' \in R^\circ \} \subseteq \Ball((X \loneplus Y)^*) ; \]
  multiplicative conjunction is given by $(X,S) \ptimes (Y, R) \eqdef (X \ptimes Y, (S \otimes R)^{\circ\circ})$, where
  $S\otimes R \eqdef \{s\otimes r \ | \ s\in S, \ r \in R\}$;
  multiplicative disjunction by $(X,S) \itimes (Y, R) \eqdef (X \itimes Y, (S^\circ \otimes R^\circ)^\circ);$
  tensor unit for both tensors is given by $(\mathbb{C}, \{1\})$.
\end{restatable}

The polarity that we introduced behaves very well w.r.t. vN-(co)algebras.
In order to see this, we first introduce some notation.
For a vN-coalgebra $C$, we say that an element $p \in C$ is \emph{positive}
if the map $1 \mapsto p \colon \mathbb C \to C$ is positive and we write $p \geq 0$ to indicate this.
The \emph{density operators} on $C$ are defined to be the set $P_C \eqdef \{ p \in C \ |\ \varepsilon(p) = 1 \text{ and } p \geq 0 \}$.
Note that if $C = T_n$, then $P_C$ is precisely the set $D_n$ of $n \times n$ density matrices, i.e. positive matrices whose trace is $1.$

\begin{restatable}{theorem}{QvNAlg}
  \label{thm:vn-alg-include}
  If $A \in \vNAlg$, then $(A, \{1_A\}) \in \Ob(\QQ).$
  Moreover, the functor $H \colon \HH \to \QQ$ defined by $H(A) \eqdef (A, \{1_A\})$
  and $H(f) \eqdef f,$ is fully faithful, strict monoidal w.r.t $\itimes$, and it strictly preserves
  finite products.
\end{restatable}

\begin{restatable}{theorem}{QvNCoalg}
  If $C \in \vNCoalg$, then $(C, P_C) \in \Ob(\QQ).$
  Moreover, the functor $S \colon \SSS \to \QQ$ defined by $S(C) \eqdef (C, P_C)$
  and $S(f) \eqdef f,$ is fully faithful, strict monoidal w.r.t $\ptimes$, and it strictly preserves
  finite coproducts.
\end{restatable}

\begin{example}
  For a unitary matrix $u \in \MM_{2^n},$ the c.i.i. $a \mapsto u a u^* : T_{2^n} \to T_{2^n}$
  describes a unitary evolution of the system with respect to $u$ in the
  Schrödinger picture. It preserves density operators and is thus a map in
  $\QQ((T_{2^n}, D_{2^n}), (T_{2^n}, D_{2^n}))$.
  The corresponding c.i.i. in the Heisenberg picture, $a \mapsto u^* a u : M_{2^n} \to M_{2^n}$,
  preserves the unit and is thus a map in
  $\QQ((M_{2^n}, \{1_{2^n}\}), (M_{2^n}, \{1_{2^n}\}))$.
  Since $(T_2, D_2) \ptimes \cdots \ptimes (T_2,D_2) \cong (T_{2^n}, D_{2^n})$ and likewise $(M_{2}, \{1_2\}) \itimes \cdots \itimes (M_2, \{1_2\}) \cong
  (M_{2^n}, \{1_{2^n}\})$, these objects of $\QQ$ can be
  understood as representing an array of $n$ qubits in the \Schrod{}/Heisenberg picture.
\end{example}

%%%%%%%%%%%%%%%%%%%%%%%%%%%%%%%%%%%%%%%%%%%%%%%%%%%%%%%%%%%%%%%%%%%%%%%%%%%%%%
% THE TWO TABLES
%%%%%%%%%%%%%%%%%%%%%%%%%%%%%%%%%%%%%%%%%%%%%%%%%%%%%%%%%%%%%%%%%%%%%%%%%%%%%%

\begin{figure}[t]
  \begin{tabular}{ |c|c|c| }
    \hline
    {\scriptsize \Schrod{} Picture} & $\SSS \fullsubright \QQ$ & $\text{LL}_{+}$ \\
    \hline
    {\scriptsize System description} & $\mathcal C, \mathcal D$ & $P, R$\\
    \hline
    {\scriptsize Quantum composition} & $\mathcal C \ptimes \mathcal D $ &  $P \otimes R$ \rule{0pt}{2.5ex} \\
    \hline
    {\scriptsize Classical composition} & $\mathcal C \loneplus \mathcal D$ & $P \oplus R$ \rule{0pt}{2.5ex} \\
    \hline
    {\scriptsize Quantum operation} & $\mathcal C \xrightarrow{\mathrm{CPTP}} \mathcal D$ & $P \vdash R$ \rule{0pt}{2.5ex} \\
    \hline
  \end{tabular}
  \quad
  \begin{tabular}{ |c|c|c| }
    \hline
    {\scriptsize Heisenberg Picture} & $\HH \fullsubright \QQ$ & $\text{LL}_{-}$ \\
    \hline
    {\scriptsize System description} & $\mathcal A, \mathcal B$ & $N, M$ \\
    \hline
    {\scriptsize Quantum composition} & $\mathcal A \itimes \mathcal B $ & $N \parr M$ \rule{0pt}{2.5ex} \\
    \hline
    {\scriptsize Classical composition} & $\mathcal A \linftyplus \mathcal B$ & $N \aconj M$ \rule{0pt}{2.5ex} \\
    \hline
    {\scriptsize Quantum operation} & $\mathcal B \xrightarrow{\mathrm{CPU}} \mathcal A$ & $M \vdash N$ \rule{0pt}{2.5ex} \\
   \hline
  \end{tabular}
  \caption{\Schrod{}/Heisenberg picture and positive/negative logical polarity.}
  \label{fig:the-two-tables}
\end{figure}

We are now justified in presenting the summary provided by Figure
\ref{fig:the-two-tables}, on which we elaborate. In our model $\QQ$, formulas
in MALL admit interpretations as objects of $\QQ$ and proofs are interpreted as
morphisms of $\QQ.$ Formulas with positive (negative) logical polarities admit
natural interpretations as vN-coalgebras $\mathcal C, \mathcal D$ (vN-algebras
$\mathcal A, \mathcal B$) and proofs between such formulas correspond precisely
to the CPTP (CPU) maps, i.e. the quantum operations in the \Schrod{}
(Heisenberg) picture. Moreover, this correspondence is preserved by classical
composition (spacewise composition of systems where only a limited amount of
classical interactions are possible) and by quantum composition (spacewise
composition with the full range of quantum interactions possible, e.g.
entanglement), in both pictures, whose interpretations are provided by the
respective tables.

\paragraph{Pure Quantum Computation.} We showed that $\QQ$ captures \emph{mixed state} quantum computation in both pictures.
In fact, $\QQ$ also has interesting properties that are relevant to \emph{pure
state} quantum computation, where unitarity is very important. If $M$ is a f.d.
vN-algebra and $u \in M$ a unitary element, i.e. $uu^* = 1_M = u^*u,$ then
$(M, \{u\}) \in \Ob(\QQ),$ i.e. the set $\{ u\}$ is bipolar in $M$.
We can reason in $\QQ$ about interesting higher-order (pure state) maps,
such as the pure state \emph{quantum switch}. This is the linear map
$ \qsw \colon M_n \ptimes M_n \to M_{2n} $ defined by
$\qsw(a \otimes b) \eqdef (\ket 0 \bra 0 \otimes (ab) ) + (\ket 1 \bra 1 \otimes (ba) ) $
and it is a complete contraction (special case of \cite[\S IV]{os-lics}). For any
$n \times n$ unitary matrices $u$ and $v$,
we have object equalities $(M_n, \{u\}) \ptimes (M_n, \{v\}) = (M_n \ptimes M_n, \{u \otimes v\}^\bipolar) = (M_n \ptimes M_n, \{u \otimes v\})$
and $\qsw$ becomes a morphism of our model $\QQ$ (see Appendix \ref{app:qsw} for a proof) with type
$ \qsw \colon (M_n, \{u\}) \ptimes (M_n, \{v\}) \to \left(M_{2n}, \{ \ket 0 \bra 0 \otimes (uv)  + \ket 1 \bra 1 \otimes (vu) \} \right) , $
where all three bipolar subsets consist of unitary matrices.
Therefore $\qsw$ may be recognised as a valid (unitary-preserving) higher-order map in $\QQ.$
The ``higher-order'' description of $\qsw$ is justified by the fact that in pure state computation, we think of the unitary elements $u,v \in M_n$
as first-order (reversible) functions.

One of the reasons $\qsw$ is interesting, is that it does not admit any reasonable multi-linear decomposition, i.e.
there are no complete contractions $\varphi_1 \colon M_n \ptimes M_n \to M_{2n,m}$ and $\varphi_2 \colon M_n \ptimes M_n \to M_{m,2n}$
such that $\qsw(a \otimes b) = \varphi_1(a) \varphi_2(b).$ In fact, this can be determined by replacing $\ptimes$ with $\hagertimes$ and checking if $\qsw$ is still a complete contraction
(see \cite[\S IV]{os-lics} for more details).
The main reason why this is the case is that the Haagerup tensor actually \emph{characterises} maps that do admit similar multi-linear decompositions \cite[\S 9.4]{er2000operator}.
Such reasoning can also be done in $\QQ$, because we have a valid object
$(M_n \hagertimes M_n, \{u \otimes v\}) \in \Ob(\QQ),$ but $\qsw \colon (M_n \hagertimes M_n, \{u \otimes v\}) \to \left(M_{2n}, \{ \ket 0 \bra 0 \otimes (uv)  + \ket 1 \bra 1 \otimes (vu) \} \right)$
is \emph{not} a complete contraction, so we can conclude following the same arguments as in \cite[\S IV]{os-lics}.
This shows our model $\QQ$ can be used to reason about complex higher-order behaviour via the Haagerup tensor, which is outside the scope of Linear Logic,
and suggests a new direction for future work.

\section{Discussion and Future Work}
We described a new approach to quantum coherence spaces, in which
the proofs of formulas with positive logical polarity correspond precisely to CPTP
maps and proofs of formulas with negative logical polarity correspond precisely
to CPU maps, remedying the issues identified by Selinger \cite{qcs-selinger} in
the approach proposed by Girard \cite{qcs1,qcs2,qcs3}.
Other relevant works related to quantum theory include
\cite{mall-simons-kissinger} which describes a model of MALL and \cite{qfpc}
which describes a model of LL. These papers use techniques based on gluing and
orthogonality or are inspired by this. Our model exhibits two main differences:
first, all of our polarised formulas have physical interpretations and a
mathematical formulation based on mathematical physics; second, complete
positivity is not assumed a priori, but it is \emph{derived}
only where necessary, i.e. for \emph{mixed state} quantum
computation. We showed in \secref{sec:qcs} that our model has the
advantage that it can also be used to reason about \emph{pure state} quantum
computation, where complete positivity makes no sense.
Our model can be used to interpret a quantum lambda calculus using standard techniques, 
and since it can talk about both pure state and mixed state quantum computation, one direction
for future work is to use it to extract suitable design for type systems
and/or logics that can combine classical and quantum control.

Another direction for future work is related to general recursion -- our morphisms
should be seen as ``total'' rather than ``partial'' from the point of view of
domain theory \cite{abramsky-jung}. To get more interesting domain-theoretic
properties, one should replace trace-preserving maps with
\emph{trace-non-increasing} (TNI) maps in the
\Schrod{} picture, replace unital
maps with \emph{subunital} (SU) ones in the Heisenberg picture.
Indeed, in finite
dimensions, the Loewner order gives both the CPTNI maps and the CPSU maps the
structure of a continuous domain \cite{qcs-selinger,vqpl}, which is a very nice
kind of dcpo from a domain-theoretic point of view.
To do this, we would need to modify the
polarity/orthogonality in \secref{sec:qcs} accordingly. We believe that this
might be achieved through the noncommutative generalisation
of convexity, called \emph{matrix convexity} \cite[\S 5.5]{er2000operator}\cite{EW1997117}. This would hopefully strengthen the links to probabilistic coherence spaces, which can indeed be used to model general recursion.

We note that a categorical approach to (infinite dimensional)
\emph{pre-C$^*$-algebras} as monoid objects in a monoidal category with an
involution is presented in \cite{fl2025},
however, their definition requires them to work with unbounded linear maps.
By comparison, our approach based on using the Haagerup tensor
does not have the same issue, even if we were to use infinite-dimensional spaces, because
this tensor may be used to \emph{characterise} unital operator algebras \cite{Pisier_2003,operator-algebras-haagerup}.
Note that we still get some structure relevant to categorical probability theory
as subcategories of (co)commutative vN-(co)algebras in
$\vNalg^\text{op} (\simeq \vNcoalg)$ and $\HH^\text{op} (\simeq \SSS)$ form Markov
categories w.r.t. the completely injective (projective) tensors.

Another direction for future work would be to try to generalise the categorical definition of
vN-algebras and vN-coalgebras to infinite dimensions, where the dual and
predual do not coincide.
Note that this would also require us to adapt the choice of some of the monoidal
structures used, as in \cite{er-shuffle}, because the Haagerup tensor is not
self-dual in a De Morgan sense for infinite-dimensional operator spaces (in general). Developing the theory of
infinite-dimensional vN-coalgebras can hopefully be used to further strengthen the
observations in \cite{os-lics} regarding the connections between polarised LL
and the Heisenberg-Schrödinger duality.
Finally, one could investigate if and how
our linear logic approach to the Heisenberg-Schrödinger duality relates to
the categorical logic approach of state-and-effect
triangles \cite{Jacobs_2017}, given that $\HH^\text{op} \simeq \SSS$
is indeed an effectus.

\vspace{1mm}
\noindent\textbf{Acknowledgements.} We thank the anonymous reviewers for
their feedback which led to multiple improvements of the paper.
We also thank Benoît Valiron, Bert Lindenhovius, Titouan Carette, and James Hefford for
discussions and/or useful feedback.
This work has been partially funded by the French
National Research Agency (ANR) within the framework of ``Plan France 2030'',
under the research projects EPIQ ANR-22-PETQ-0007, HQI-Acquisition
ANR-22-PNCQ-0001, HQI-R\&D ANR-22-PNCQ-0002, and also by CIFRE 2022/0081.

%%%%%%%%%%%%%%%%%%%%%%%%%%%%%%%%%%%%%%%%%%%%%%%%%%%%%%%%%%%%%%%%%%%%%%%%%%%%%%
% END of main body
%%%%%%%%%%%%%%%%%%%%%%%%%%%%%%%%%%%%%%%%%%%%%%%%%%%%%%%%%%%%%%%%%%%%%%%%%%%%%%

\newpage
\bibliographystyle{splncs04}
\bibliography{bib}
\newpage

%%%%%%%%%%%%%%%%%%%%%%%%%%%%%%%%%%%%%%%%%%%%%%%%%%%%%%%%%%%%%%%%%%%%%%%%%%%%%%
% Appendix
%%%%%%%%%%%%%%%%%%%%%%%%%%%%%%%%%%%%%%%%%%%%%%%%%%%%%%%%%%%%%%%%%%%%%%%%%%%%%%

\appendix

\section{Note about these appendices}

We intend to polish the material presented in these appendices and submit a
journal version of this paper in the future. In the meantime, we hope that
readers who are interested in the technical details find the information
contained herein useful.

\newpage
\section{The category $\FdOS$}
\label{app:fdos}

In this appendix we provide some omitted proofs about the categorical properties
of $\FdOS.$

%%%%%%%%%%%%%%%%%%%%%%%%%%%%%%%%%%%%%%
\subsection{Model of MALL}
%%%%%%%%%%%%%%%%%%%%%%%%%%%%%%%%%%%%%%

The definition of the completely projective tensor that we provided in
\secref{sec:fdos} is not the standard one that may be found in the literature.
For \emph{finite-dimensional} operator spaces, it is equivalent to the standard
one and we prove this here. The reason for presenting the simpler definition in
\secref{sec:fdos} is that it has a more logical/categorical character and we
believe it makes the paper more accessible to experts in these topics.

We begin the development by recalling the standard definition for the
completely projective tensor product $X \ptimes Y$ and we prove at the end of
this subsection that this determines the same o.s.s on $X \otimes Y$ as the one
in \cref{def:ptimes-simple}.

\begin{definition}[{Completely projective tensor \cite[\S 7.1]{er2000operator}}]
  \label{def:ptimes-standard}
  Let $X, Y \in \FdOS$. The \emph{completely projective tensor product},
  written $X \ptimes Y$, is the vector space
  $X \otimes Y$ together with the o.s.s defined for
  $v \in \MM_n(X \otimes Y)$ by:
  \begin{align*}
    \norm{v}_{\land} \eqdef \inf \{ & \norm{\alpha} \norm{x} \norm{y} \norm{\beta} \  :\
    p \in \mathbb N, q \in \mathbb N, x \in M_p(X) , y \in M_q(Y), \\
    & \alpha \in M_{n, pq}, \beta \in M_{pq, n} , \text{ and } v = \alpha (x \otimes y) \beta\} .
  \end{align*}
  Here $x\otimes y\in M_{pq}(X\otimes Y)$ stands for the ``tensor product of matrices'' given by
  $x\otimes y\eqdef [x_{ij}\otimes y_{kl}]_{(i,k),(j,l)},$ which generalises the usual Kronecker product.
\end{definition}

\begin{proposition}[{\cite[Sec. 7-9]{er2000operator}}]
  \label{prop:tensorprops}
  The tensor product $\alphatens$
  (for $\alpha = \widehat{\ }, \widecheck{\ }, h$) is
  \begin{itemize}
    \item Functorial wrt c.b maps and c.c maps
    \item Associative: $(A \alphatens B) \alphatens C \cong A \alphatens (B \alphatens C)
    : (a \otimes b) \otimes c \mapsto a \otimes (b \otimes c)$
    \item Has unitors: $\mathbb{C} \alphatens A \cong A \cong A \alphatens \mathbb{C} :
    \lambda \otimes a \mapsto \lambda a \mapsfrom a \otimes \lambda$
  \end{itemize}
  and for $\alpha = \widehat{\ }, \widecheck{\ }$
  \begin{itemize}
    \item Symmetric: $A \alphatens B \cong B \alphatens A: a \otimes b \mapsto b \otimes a$
  \end{itemize}
\end{proposition}

\fdosMonoidal*

\begin{proof}
  All three tensors can be extended to bifunctors and the associators, unitors
  and symmetries are the ones given in \cref{prop:tensorprops}. The coherence
  diagrams are straightforward to verify.
  To see that the functor $V: \fdOS \to \FdVect$ is a strict monoidal, note that
  the underlying vector space of $A \alphatens B$
  (for $\alpha = \widehat{\ }, \widecheck{\ }, h$) is $A\otimes B$ for any f.d.
  operator spaces $A$ and $B$ and
  that (symmetric) monoidal natural isomorphisms in \cref{prop:tensorprops}
  coincide with the ones for $\FdVect$.
  \qed
\end{proof}

\begin{lemma}
  \label{lemma:fdos-smcc}
  The symmetric monoidal category $(\fdOS, \mathbb{C}, \projtens)$
  is closed w.r.t. $\CB(-,-): \fdOS^\mathrm{op} \times \fdOS \to \fdOS$.
\end{lemma}

\begin{proof}
  Follows from restricting \cite[Prop. III.33]{os-lics} to $\fdOS$.
  \qed
\end{proof}

\begin{proposition}
  \label{prop:fdos-star-aut}
  The category \fdOS \ is $*$-autonomous with $\mathbb{C}$ as dualizing object.
\end{proposition}

\begin{proof}
  Because of \cref{lemma:fdos-smcc}, it suffices to prove that $\mathbb{C}$
  is a global dualizing object.
  Let $A \in \fdOS$, currying the evaluation map
  $\text{ev}_{A,\mathbb{C}}: \CB(A,\mathbb{C}) \otimes A \to \mathbb{C}$, gives
  us the canonical embedding of $A$ into its double dual
    \begin{equation*}
      \begin{array}{r c l}
        d_A : A & \to & \CB(\CB(A,\mathbb{C}),\mathbb{C})\\
        a & \mapsto & (f \mapsto f(a))
      \end{array}
    \end{equation*}
  which by \cite[Prop. 3.1.2]{er2000operator} is a c.i. inclusion.
  As $A$ is f.d, we have by the rank-nullity
  argument that this is a isomorphism.
  \qed
\end{proof}

\fdosMall*

\begin{proof}
  Follows by combining \cref{prop:fdos-star-aut} with results about (co)products
  from \cref{sub:direct-sums}.
  \qed
\end{proof}

We can now easily prove that the two definitions of the completely projective
tensor are equivalent on finite-dimensional operator spaces, i.e. they define
the same o.s.s. on $X \otimes Y.$

\begin{proposition}
  Let $X$ and $Y$ be two f.d. operator spaces. Then the o.s.s
  on $X \otimes Y$ given by
  Definition \ref{def:ptimes-simple} coincides with the o.s.s of Definition \ref{def:ptimes-standard}.
\end{proposition}
\begin{proof}
  From Theorem \ref{thm:fdos-mall}, we know that $\FdOS$ is $*$-autonomous and
  $X \ptimes Y \cong (X^* \itimes Y^*)^*$ when $\ptimes$ is defined in the
  standard way, i.e. as in Definiton \ref{def:ptimes-standard}. But this c.i.i
  is exactly how the o.s.s. of $X \ptimes Y$ is defined in
  \cref{def:ptimes-simple}, so the two operator space structures are exactly
  equal. \qed
\end{proof}

%%%%%%%%%%%%%%%%%%%%%%%%%%%%%%%%%%%%%%
\subsection{BV-categorical structure}
%%%%%%%%%%%%%%%%%%%%%%%%%%%%%%%%%%%%%%

%%%%%%%%%%%%%%%% BV STUFF %%%%%%%%%%%%%%
\begin{definition}
  Let $\mathbf{C}$ be a category with two monoidal structure
  $(\mathbf{C}, I, \otimes)$ and $(\mathbf{C}, J, \vartriangleleft)$, then
  $(J, \vartriangleleft)$ is \emph{normal duoidal} to $(I, \otimes)$ if
  we have a natural transformation:
  \begin{equation*}
    w: (A\vartriangleleft B)\otimes (C\vartriangleleft D) \to
    (A \otimes C) \vartriangleleft (B\otimes D)
  \end{equation*}
  (the morphism $w$ is often called \emph{weak interchange})
  and and morphisms
  \begin{equation*}
    i:I \cong J \qquad \delta_I : I \to I \vartriangleleft I \qquad \mu_J : J \otimes J \to J
  \end{equation*}

  Along with coherence conditions for associativity,
  \begin{equation}
    \label{eq:assoc-w-seq}
    \begin{tikzcd}
      ((A\vartriangleleft B) \otimes (C \vartriangleleft D)) \otimes (E \vartriangleleft F) \arrow[r, "\alpha"] \arrow[d, "w \otimes \mathrm{id}"']             & (A\vartriangleleft B) \otimes ((C \vartriangleleft D) \otimes (E \vartriangleleft F)) \arrow[d, "\mathrm{id} \otimes w"]      \\
      ((A \otimes C) \vartriangleleft (B\otimes D)) \otimes (E \vartriangleleft F) \arrow[d, "w"']                                 & (A\vartriangleleft B) \otimes ((C \otimes E) \vartriangleleft (D \otimes F)) \arrow[d, "w"] \\
      ((A \otimes C) \otimes E) \vartriangleleft ((B \otimes D) \otimes F)  \arrow[r, "\alpha \vartriangleleft \alpha"'] & (A \otimes (C \otimes E)) \vartriangleleft (B \otimes (D \otimes F))
      \end{tikzcd}
  \end{equation}
  \bigskip

  \begin{equation}
    \label{eq:assoc-w-tens}
    \begin{tikzcd}
      ((A \vartriangleleft C) \vartriangleleft E) \otimes ((B \vartriangleleft D) \vartriangleleft F)  \arrow[r, "\alpha \otimes \alpha"] \arrow[d, "w"']   & (A \vartriangleleft (C \vartriangleleft E)) \otimes (B \vartriangleleft (D \vartriangleleft F)) \arrow[d, "w"]  \\
      ((A \vartriangleleft C) \otimes (B \vartriangleleft D)) \vartriangleleft (E \otimes F) \arrow[d, "w \vartriangleleft \mathrm{id}"']   & (A\otimes B) \vartriangleleft ((C \vartriangleleft E) \otimes (D \vartriangleleft F)) \arrow[d, "\mathrm{id} \vartriangleleft w"] \\
      ((A\otimes B) \vartriangleleft (C \otimes D)) \vartriangleleft (E \otimes F) \arrow[r, "\alpha"']           & (A\otimes B) \vartriangleleft ((C \otimes D) \vartriangleleft (E \otimes F))     \\
    \end{tikzcd}
  \end{equation}

  and unitality,

  \begin{equation}
    \label{eq:unital-tens}
    \begin{tikzcd}
      I \otimes (A \vartriangleleft B) \arrow[r, "\delta_I \otimes \text{id} "]  & (I \vartriangleleft I) \otimes (A \vartriangleleft B) \arrow[d, "w"] \\
      A \vartriangleleft B \arrow[u, "\lambda_{A\vartriangleleft B}"] \arrow[r, "\lambda_A \vartriangleleft \lambda_B"'] & (I \otimes A) \vartriangleleft (I \otimes B)
    \end{tikzcd}
    \qquad
    \begin{tikzcd}
     (A \vartriangleleft B) \otimes I \arrow[r, "\text{id} \otimes \delta_I "]  &  (A \vartriangleleft B) \otimes (I \vartriangleleft I) \arrow[d, "w"] \\
      A \vartriangleleft B \arrow[u, "\rho_{A\vartriangleleft B}"] \arrow[r, "\rho_A \vartriangleleft \rho_B"'] & (A \otimes I) \vartriangleleft (B \otimes I)
      \end{tikzcd}
  \end{equation}
  \bigskip

  \begin{equation}
    \label{eq:unital-seq}
    \begin{tikzcd}
      J \vartriangleleft (A \otimes B) & \arrow[l, "\mu_J \vartriangleleft \text{id}"']  (J \otimes J) \vartriangleleft (A \otimes B) \\
      A \otimes B \arrow[u, "\lambda_{A\otimes B}"] \arrow[r, "\lambda_A \otimes \lambda_B"'] & (J \vartriangleleft A) \otimes (J \vartriangleleft B) \arrow[u, "w"']
    \end{tikzcd}
    \qquad
    \begin{tikzcd}
      (A \otimes B) \vartriangleleft J & \arrow[l, "\text{id} \vartriangleleft \mu_J"']  (A \otimes B) \vartriangleleft (J \otimes J) \\
      A \otimes B \arrow[u, "\rho_{A\otimes B}"] \arrow[r, "\rho_A \otimes \rho_B"'] & (A \vartriangleleft J) \otimes (B \vartriangleleft J) \arrow[u, "w"']
    \end{tikzcd}
  \end{equation}
  as well as diagrams saying that $I$ is a comonoid in
  $(\mathcal{C}, J, \vartriangleleft)$ and that $J$ is a monoid in
  $(\mathcal{C}, I, \otimes)$.
\end{definition}

\begin{definition}\cite[Sec. 3]{bv-category}
  A \emph{BV-category with negation} is a $*$-autonomous category
  ${(\mathbf{C}, I, \otimes, \multimap)}$ with an additional monoidal
  structure
  $(\mathbf{C}, J, \vartriangleleft)$ that is \emph{normal duoidal} to
  $(\mathbf{C}, I, \otimes)$, such that the interchange commutes with the
  symmetry of $\otimes$ in the following way:
  \begin{equation}
    \label{def:w-comp-sym}
    \begin{tikzcd}
      (A\vartriangleleft B)\otimes (C\vartriangleleft D) \arrow[d, "w"'] \arrow[r, "\cong"] & (C\vartriangleleft D) \otimes (A\vartriangleleft B) \arrow[d, "w"]\\
      (A \otimes C) \vartriangleleft (B\otimes D) \arrow[r, "\cong"'] & (C \otimes A) \vartriangleleft (D \otimes B)
    \end{tikzcd}
  \end{equation}
\end{definition}

\begin{remark}
  \label{rem:parr-duoidal}
  Note that in a \emph{BV-category with negation} we also have a duoidal structure
  between $(\mathbf{C}, \bot, \parr)$ and $(\mathbf{C}, J, \vartriangleleft)$,
  where $(-) \parr (-) \eqdef (-)^* \multimap (-)$.
\end{remark}

\begin{proposition}\label{BV-cat-fdVect}
  The category $\FdVect$  is a \emph{BV-category with negation} w.r.t.
  $(\FdVect, \mathbb{C}, \otimes, L(-,-))$ and
  $(\FdVect, \mathbb{C}, \otimes)$.
\end{proposition}

\begin{proof}
  $\FdVect$ is compact closed, hence $*$-autonomous, and any symmetric monoidal
  structure is normal duoidal to itself (special case of \cite[Prop. 6.10.]{AM2010}).
  \qed
\end{proof}

\begin{proposition}\label{BV-cat}
  The category $\fdOS$ is a \emph{BV-category with negation} w.r.t.
  $(\fdOS, \mathbb{C}, \projtens, \CB(-,-))$ and
  $(\fdOS, \mathbb{C}, \hagertimes)$.
\end{proposition}

\begin{proof}
  We know by \cref{prop:fdos-star-aut} that $\fdOS$ is $*$-autonomous.
  By \cref{prop:shuffle} we have natural weak interchange maps in
  $\FdOS$ and as we have a common unit for all tensors, the result now follows from
  from \cref{prop:fdos-monoidal} as the appropriate morphisms are mapped to the
  BV-stucture on $\FdVect$.
  \qed
\end{proof}

%%%%%%%%%%%%%%%%%%%%%%%%%%%%%%%%%%%%%%
\subsection{The conjugate/opposite functorial involutions}
%%%%%%%%%%%%%%%%%%%%%%%%%%%%%%%%%%%%%%

We are going to use the following proposition multiple times, so we recall it
here for convenience.

\opConjugateCB*

In particular, we use Proposition \ref{prop:op-conjugate-cb} to easily prove
many of the following results.

\functorialInvolutions*
\begin{proof}
  The last two facts are already known, so we provide references.

  \textbf{Case} $(X \hagertimes Y)_c = X_c \hagertimes Y_c$. It is claimed
  in \cite[p. 97]{Pisier_2003} that this is a completely isometric isomorphism. In our
  presentation, $X_c$ has the same underlying set as $X$, so this is actually an equality.

  \textbf{Case} $(X \hagertimes Y)_o \cong Y_o \hagertimes X_o$. This is claimed
  in \cite[pp. 96--97]{Pisier_2003}.

  \textbf{Case} $(X \itimes Y)_c = X_c \itimes Y_c.$ This is easy to prove using Proposition \ref{prop:op-conjugate-cb}.
  Let $\iota_{X,Y} \colon X \itimes Y \cong \CB(Y^*, X)$ be the canonical completely isometric isomorphism from Definition \ref{def:itimes}.
  More specifically, $\iota_{X,Y}(x \otimes y) \eqdef (\varphi \mapsto \varphi(y)x).$
  The following composite
  \[ (X \itimes Y)_c \xrightarrow{\iota_c} \CB(Y^*, X)_c = \CB((Y^*)_c, X_c) \cong \CB((Y_c)^*, X_c) \xrightarrow{\iota^{-1}} X_c \itimes Y_c \]
  consists of completely isometric isomorphisms. Evaluating on elementary tensors
  \begin{align*}
    x \otimes y &\xmapsto{\iota_c} (\varphi \mapsto \varphi(y)x) \\
                &= \left(\varphi \mapsto \overline{\varphi(y)} \star x \right) \\
                &\xmapsto{\cong} \left(\varphi \mapsto \varphi(y) \star x \right) \\
                &\xmapsto{\iota^{-1}} x \otimes y
  \end{align*}
  shows the above composite is the identity map. Here we have written $\star$ for the scalar multiplication in $X_c$ and the vertically aligned symbols
  indicate the action of the isomorphisms.
  Since this c.i.i. is the identity map, the equality follows.

  \textbf{Case} $(X \itimes Y)_o = X_o \itimes Y_o.$ The proof is analogous to the previous case.
  The following composite
  \[ (X \itimes Y)_o \xrightarrow{\iota_o} \CB(Y^*, X)_o = \CB((Y^*)_o, X_o) = \CB((Y_o)^*, X_o) \xrightarrow{\iota^{-1}} X_o \itimes Y_o \]
  consists of completely isometric isomorphisms. Evaluating on elementary tensors, we see that
  this c.i.i. is the identity map and the equality follows.

  \textbf{Case} $(X \ptimes Y)_c = X_c \ptimes Y_c.$ This can be proven by using the fact that $(X \itimes Y)_c = X_c \itimes Y_c$ and adapting the proof of the latter fact.
  Let $\iota_{X,Y} \colon X \ptimes Y \cong (X^* \itimes Y^*)^*$ be the canonical c.i.i. defined by
  $\iota_{X,Y}(x \otimes y) \eqdef (\varphi \otimes \psi) \mapsto \varphi(x)\psi(y).$
  The following composite
  \[ (X \ptimes Y)_c \xrightarrow{\iota_c} ((X^* \itimes Y^*)^*)_c \cong ((X^* \itimes Y^*)_c)^* = ((X^*)_c \itimes (Y^*)_c)^* \cong
  ((X_c)^* \itimes (Y_c)^*)^* \xrightarrow{\iota^{-1}} X_c \ptimes Y_c \]
  consists of completely isometric isomorphisms. Evaluating on elementary tensors
  \begin{align*}
    (x \otimes y) &\xmapsto{\iota_c} \left( (\varphi \otimes \psi) \mapsto \varphi(x) \cdot \psi(y) \right) \\
                  &\xmapsto{\cong} \left( (\varphi \otimes \psi) \mapsto \overline{\varphi(x) \cdot \psi(y)} \right) \\
                  &\xmapsto{\cong} \left( (\varphi \otimes \psi) \mapsto \overline{\overline{\varphi(x)} \cdot \overline{\psi(y)}} \right) \\
                  &= \left( (\varphi \otimes \psi) \mapsto \varphi(x) \cdot \psi(y) \right) \\
                  &\xmapsto{\iota^{-1}} (x \otimes y)
  \end{align*}
  shows that this composite is the identity map and the equality follows.

  \textbf{Case} $(X \ptimes Y)_o = X_o \ptimes Y_o.$ This can be proven analogously to the previous case.
  The following composite
  \[ (X \ptimes Y)_o \xrightarrow{\iota_o} ((X^* \itimes Y^*)^*)_o =
  ((X_o)^* \itimes (Y_o)^*)^* \xrightarrow{\iota^{-1}} X_o \ptimes Y_o \]
  consists of completely isometric isomorphisms. After evaluating on elementary tensors
  we see
  that the composite is the identity map and the equality follows.

  \textbf{Case} $(X \linftyplus Y)_o = X_o \linftyplus Y_o.$ Clearly, the identity function
  gives a linear isomorphism $\mathrm{id} \colon (X \linftyplus Y)_o \cong X_o \linftyplus Y_o.$
  To show that it is completely isometric, let $[(x_{ij}, y_{ij})] \in M_n((X \linftyplus Y)_o).$ We verify:
  \begin{align*}
    \norm{[(x_{ij}, y_{ij})]}_{M_n((X \linftyplus Y)_o)} &= \norm{[(x_{ji}, y_{ji})]}_{M_n((X \linftyplus Y))} \\
                                                         &= \norm{([x_{ji}], [y_{ji}])}_{M_n(X) \linftyplus M_n(Y)} \\
                                                         &= \mathrm{max}(\norm{[x_{ji}]}_{M_n(X)}, \norm{[y_{ji}]}_{M_n(Y)}) \\
                                                         &= \mathrm{max}(\norm{[x_{ij}]}_{M_n(X_o)}, \norm{[y_{ij}]}_{M_n(Y_o)}) \\
                                                         &= \norm{([x_{ij}], [y_{ij}])}_{M_n(X_o) \linftyplus M_n(Y_o)} \\
                                                         &= \norm{[(x_{ij}, y_{ij})]}_{M_n(X_o \linftyplus Y_o)} ,
  \end{align*}
  where each equality holds by definition of the relevant (matrix) norm.

  \textbf{Case} $(X \linftyplus Y)_c = X_c \linftyplus Y_c.$ The proof is fully analogous to the previous case.
  The identity function a linear isomorphism $\mathrm{id} \colon (X \linftyplus Y)_c \cong X_c \linftyplus Y_c.$
  To show that it is completely isometric, one can simply follow the same steps as in the previous proof while using the o.s.s of the conjugated spaces.

  \textbf{Case} $(X \loneplus Y)_c = X_c \loneplus Y_c.$ Let $\iota_{X,Y} \colon X \loneplus Y \cong (X^* \linftyplus Y^*)^*$ be the canonical completely
  isometric isomorphism that defines the o.s.s of $X \loneplus Y$. It is given by $\iota_{X,Y}(x,y) = \left( (\varphi, \psi) \mapsto \varphi(x) + \psi(y) \right).$
  The following composite
  \[ (X \loneplus Y)_c \xrightarrow{\iota_c} ((X^* \linftyplus Y^*)^*)_c \cong ((X^* \linftyplus Y^*)_c)^* \cong ( (X_c)^* \linftyplus (Y_c)^*)^* \xrightarrow{\iota^{-1}} X_c \loneplus Y_c \]
  consists of completely isometric isomorphisms. Evaluating on elements
  \begin{align*}
    (x , y)       &\xmapsto{\iota_c} \left( (\varphi , \psi) \mapsto \varphi(x) + \psi(y) \right) \\
                  &\xmapsto{\cong} \left( (\varphi , \psi) \mapsto \overline{\varphi(x) + \psi(y)} \right) \\
                  &\xmapsto{\cong} \left( (\varphi , \psi) \mapsto \overline{\overline{\varphi(x)} + \overline{\psi(y)}} \right) \\
                  &= \left( (\varphi , \psi) \mapsto \varphi(x) + \psi(y) \right) \\
                  &\xmapsto{\iota^{-1}} (x , y)
  \end{align*}
  shows that this composite is the identity map and the equality follows.

  \textbf{Case} $(X \loneplus Y)_o = X_o \loneplus Y_o.$ This is analogous to the previous case.
  The following composite
  \[ (X \loneplus Y)_o \xrightarrow{\iota_o} ((X^* \linftyplus Y^*)^*)_o = ((X^* \linftyplus Y^*)_o)^* = ( (X_o)^* \linftyplus (Y_o)^*)^* \xrightarrow{\iota^{-1}} X_o \loneplus Y_o \]
  consists of completely isometric isomorphisms. Evaluating on elements like in the previous case
  shows that this is the identity map and the equality follows.
  \qed
\end{proof}

\begin{lemma}
  There is a natural isomorphism $((-)_c)^* \cong ((-)^*)_c: \fdOS \to \fdOS^\mathrm{op}$
\end{lemma}

\begin{proof}
  The natural isomorphism is given by postcomposing with the isomorphism
  $\mathbb{C} \cong \mathbb{C}_c$. To see this, let $X \in \fdOS$, we unfold the
  definition of the dual to get the following chain of isomorphisms
  \begin{equation*}
    \begin{array}{rclr}
      ((X)_c)^* & \eqdef & \CB (X_c , \mathbb{C}) & \\
                & \cong  & \CB (X_c, \mathbb{C}_c) & \qquad \text{(postcomposing with $\mathbb{C} \cong \mathbb{C}_c$)} \\
                & =      & (\CB (X, \mathbb{C}))_c &  \\
                & \eqdef & (X^*)_c &
    \end{array}
  \end{equation*}
  All of these are natural, hence we are done.
  \qed
\end{proof}

\begin{lemma}
  Given operator spaces $X, Y \in \FdOS$ we have that the following commutation
  property between the \emph{conjugate/opposite} involution and the self-duality
  of the Haagerup tensor:
  \begin{equation*}
    \begin{tikzcd}
      (X^* \hagertimes Y^*)_c \arrow[equal]{d} \arrow[r, "\cong"] & ((X \hagertimes Y)^*)_c \arrow[d, "\cong"] &  & (X^* \hagertimes Y^*)_o \arrow[d, "\gamma"'] \arrow[r, "\cong"] & ((X \hagertimes Y)^*)_o \arrow[equal]{d}  \\
      (X^*)_c \hagertimes (Y^*)_c \arrow[d, "\cong"']      & ((X \hagertimes Y)_c)^* \arrow[equal]{d}          &  & (Y^*)_o \hagertimes (X^*)_o \arrow[equal]{d}                          & ((X \hagertimes Y)_o)^* \arrow[d, "\gamma^*"]         \\
      (X_c)^* \hagertimes (Y_c)^* \arrow[r, "\cong"']      & (X_c \hagertimes Y_c)^*                           &  & (Y_o)^* \hagertimes (X_o)^* \arrow[r, "\cong"']                          & (Y_o \hagertimes X_o)^*
      \end{tikzcd}
  \end{equation*}

\end{lemma}

\begin{proof}
  Straight forward diagram-chase after reducing to elementary tensors.
  \qed
  \end{proof}

\newpage

%%%%%%%%%%%%%%%%%%%%%%%%%%%%%%%%%%%%%%%%%%%%%%%%%%%%%%%%%%%%%%%%%%%%%%%%%%%%%%
\section{Proofs related to von Neumann (co)algebras}
\label{app:vN-alg}
%%%%%%%%%%%%%%%%%%%%%%%%%%%%%%%%%%%%%%%%%%%%%%%%%%%%%%%%%%%%%%%%%%%%%%%%%%%%%%

In this appendix we provide some of the omitted proofs that are relevant for vN-(co)algebras.

\subsection{von Neuman algebras}
Let us recall the definition of a vN-algebra.

\begin{definition}
  A f.d. \emph{vN-algebra} is an operator space $A$ together with
  \begin{enumerate}
    \item a complete contraction $\mu \colon A \hagertimes A \to A$, called \emph{multiplication};
    \item a complete contraction $\eta \colon \mathbb C \to A$, called \emph{unit};
    \item a complete contraction $i \colon A_c \to A_o,$ called \emph{involution},
  \end{enumerate}
  such that
  \begin{itemize}
    \item $(A, \eta, \mu)$ is a monoid object in $(\FdOS, \mathbb C, \hagertimes)$;
    \item $i$ is an \emph{involution}, meaning that the following diagram commutes:
    \begin{equation}
      \tag{involution}
      \label{eq:vNinv}
        \begin{tikzcd}[ampersand replacement= \&]
          A \arrow[equal]{d} \arrow[equal]{rr}           \&           \& A       \\
          A_{cc} \arrow[r, "i_c"'] \& A_{oc} = A_{co} \arrow[r, "i_o"'] \& A_{oo} \arrow[equal]{u}
        \end{tikzcd}
      \end{equation}
      \item $i$ is \emph{reverse multiplicative}, meaning that the following
      diagram commutes:
      \begin{equation}
        \tag{reverse multiplicative}
        \label{eq:vNrevmult}
        \begin{tikzcd}[ampersand replacement= \&]
          A_c \hagertimes A_c \arrow[d, "i\otimes i"'] \arrow[equal]{r} \& (A\hagertimes A)_c \arrow[r, "\mu"] \& A_c \arrow[d, "i"] \\
          A_o \hagertimes A_o \arrow[r, "{\gamma}"']              \& (A\hagertimes A)_o \arrow[r, "\mu"'] \& A_o
        \end{tikzcd}
      \end{equation}
      \item for every complete isometry (i.e. strong mono) $a: \mathbb{C} \to A$, the following composite in $\FdOS$ is also a
      complete isometry (i.e. strong mono)
      \begin{equation}
        \tag{C$^*$-id}
        \label{eq:c-star-id}
        \mathbb{C} \cong \mathbb{C}_o \hagertimes \mathbb{C}_c
        \xrightarrow{a \otimes a} A_o \hagertimes A_c
        \xrightarrow{A_o \otimes i} A_o \hagertimes A_o
        \xrightarrow{\gamma} (A \hagertimes A)_o
        \xrightarrow{\mu} A_o.
      \end{equation}
  \end{itemize}
\end{definition}

\begin{remark}
  The definition above is the same as the one in \cref{sec:vNAlg} only
  slightly expanded to give names to the coherence conditions.
\end{remark}

\vNAlgEquiv*
\begin{proof}
  $(\Rightarrow).$ Let $A$ be a f.d. vN-algebra in the sense of Definition \ref{def:vN-alg-standard}.
  We already explained what is the canonical o.s.s. of $A.$
  We define $\eta \colon \mathbb C \to A :: 1 \mapsto 1_A$ and
  this is clearly a contraction and therefore also a complete one.
  Since the multiplication of $A$ is bilinear, it admits a unique
  extension to a linear map $\mu \colon A \otimes A \to A.$ This
  linear map is actually a complete contraction with respect to the
  Haagerup tensor, i.e. $\mu \colon A \hagertimes A \to A$ is a
  complete contraction \cite[Proof of Theorem 17.1.2]{er2000operator}.
  The involution $i \colon A_c \to A_o$ is defined in the obvious way $i(a) \eqdef a^*$
  and this is a completely isometric isomorphism \cite[p. 65]{Pisier_2003}. The commutativity
  of diagrams \eqref{eq:vNinv} and \eqref{eq:vNrevmult} follows easily from Definition \ref{def:vN-alg-standard} (3.).
  From Definition \ref{def:vN-alg-standard} (1.) and (2.) it follows that $(A, \eta, \mu)$ is a monoid object.
  Finally, any (complete) isometry $a \colon \mathbb C \to A$ may be identified with an element $a \in A$ with $\norm{a} = 1$
  and then the composite in \eqref{eq:c-star-id} is identified with the element $a^*a \in A$. The requirement that this composite
  is a complete isometry means that $\norm{a^*a} = 1$ which follows from the C*-identity of $A,$ i.e. Definition \ref{def:vN-alg-standard} (4.).

  $(\Leftarrow).$ This direction is easier. It is obvious how to define the unit, multiplication and involution from $\eta$, $\mu$ and $i$.
  Conditions (1.) -- (3.) of Definition \ref{def:vN-alg-standard} now follow easily. The C*-norm of $A$ is simply the ground norm of $A$.
  That this norm is submultiplicative follows from \cite[Proof of Theorem 17.1.2]{er2000operator}. It only remains to show that the C*-identity holds.
  Let $a \in A$. If $a = 0$, then obviously $\norm{a^*a} = \norm{\mu(0 \otimes 0)} = 0 = \norm{a}^2$, so let us assume that $a \neq 0.$
  We can now use a simple renormalisation argument to complete the proof.
  Let $b = \frac{a}{\norm{a}}$ so that $\norm{b} = 1.$ We can now identify $b$ with a (complete) isometry $b \colon \mathbb C \to A$
  and we know that the composite
    \begin{equation*}
      \mathbb{C} \cong \mathbb{C}_o \hagertimes \mathbb{C}_c
      \xrightarrow{b \otimes b} A_o \hagertimes A_c
      \xrightarrow{A_o \otimes i} A_o \hagertimes A_o
      \xrightarrow{\gamma} (A \hagertimes A)_o
      \xrightarrow{\mu} A_o
    \end{equation*}
  is also a complete isometry. But this composite may be identified with the element $b^*b$ and therefore $\norm{b^*b} = 1.$
  The last equation now implies $\norm{a^*a} = \norm{a}^2,$ as required. \qed
\end{proof}

\vNtens*

\begin{proof}
  We prove both statements simultaneously, saving the C$^*$-identities for last.
  To prove the first statement, we start by ascertaining that
  $(A \itimes B, \eta, \mu)$ is indeed a monoid. To do so we keep
  \cref{rem:parr-duoidal} in mind.
  \newline

  \noindent We can verify the left and right unital laws:
  \begin{equation*}
    \adjustbox{scale=0.6}{
      \begin{tikzcd}
        \mathbb{C} \hagertimes (A \injtens B) \arrow[rr, "\cong"] \arrow[rrrrddd, "\lambda_{A \injtens B}"', bend right] &  & (\mathbb{C}\injtens \mathbb{C}) \hagertimes (A \injtens B) \arrow[rr, "(\eta_A\otimes \eta_B) \otimes \text{id}"] \arrow[d, "v"']           &  & (A \injtens B) \hagertimes (A \injtens B) \arrow[d, "v"]             &  & (A \injtens B) \hagertimes (\mathbb{C} \injtens \mathbb{C}) \arrow[ll, "\text{id}\otimes(\eta_A\otimes \eta_B)"'] \arrow[d, "v"]      &  & (A \injtens B) \hagertimes \mathbb{C} \arrow[ll, "\cong"'] \arrow[llllddd, "\rho_{A\injtens B}", bend left] \\
                                                                                                              &  & (\mathbb{C} \hagertimes A)\injtens (\mathbb{C}\hagertimes B) \arrow[rr, "(\eta_A \otimes A) \otimes (\eta_B \otimes B)"'] \arrow[rrdd, "\lambda_A \otimes \lambda_B"'] &  & (A\hagertimes A) \injtens (B \hagertimes B) \arrow[dd, "\mu_A \otimes \mu_B"] &  & (A \hagertimes \mathbb{C}) \injtens (B \hagertimes \mathbb{C}) \arrow[ll, "(A \otimes \eta_A) \otimes (B \otimes \eta_B)"] \arrow[lldd, "\rho_A \otimes \rho_B"] &  &                                                                                                 \\
                                                                                                              &  &                                                                                                                                                                    &  &                                                                           &  &                                                                                                                                                              &  &                                                                                                 \\
                                                                                                              &  &                                                                                                                                                                    &  & (A \injtens B)                                                            &  &                                                                                                                                                              &  &
        \end{tikzcd}
    }
  \end{equation*}
  The two inner squares commute by naturality of $v$, the lower triangles by
  assumption and outer squares by duoidality axioms (see \eqref{eq:unital-tens}).

  \noindent The multiplication is associative:
  \begin{equation*}
    \adjustbox{scale=0.6}{
      \begin{tikzcd}
        & ((A\injtens B) \hagertimes (A \injtens B)) \hagertimes (A \injtens B) \arrow[rr, "\alpha"] \arrow[d, "v_{ABAB} \otimes (A\injtens B)"']                               &              & (A\injtens B) \hagertimes ((A \injtens B) \hagertimes (A \injtens B)) \arrow[d, "(A \injtens B) \otimes v_{ABAB}"]                                 &                                                                             \\
        & ((A\hagertimes A) \injtens (B \hagertimes B)) \hagertimes (A \injtens B) \arrow[ld, "(\mu_A \otimes \mu_B) \otimes (A \injtens B)"'] \arrow[d, "v"]                     &              & (A \injtens B) \hagertimes ((A \hagertimes A) \injtens (B \hagertimes B)) \arrow[d, "v"'] \arrow[rd, "(A \injtens B) \otimes (\mu_A \otimes \mu_B)"] &                                                                             \\
          (A \injtens B) \hagertimes (A \injtens B) \arrow[d, "v"']                      & ((A\hagertimes A) \hagertimes A) \injtens ((B \hagertimes B) \hagertimes B) \arrow[rr, "\alpha \otimes \alpha"] \arrow[ld, "(\mu_A \otimes A) \otimes (\mu_B \otimes B)"] &              & (A\hagertimes (A \hagertimes A)) \injtens (B \hagertimes (B \hagertimes B)) \arrow[rd, "(A \otimes \mu_A) \otimes (B \otimes \mu_B)"']                 & (A \injtens B) \hagertimes (A \injtens B) \arrow[d, "v"]                      \\
          (A \hagertimes A) \injtens (B \hagertimes B) \arrow[rrd, "\mu_A \otimes \mu_B"'] &                                                                                                                                                                   &              &                                                                                                                                                & (A \hagertimes A) \injtens (B \hagertimes B) \arrow[lld, "\mu_A \otimes \mu_B"] \\
        &                                                                                                                                                                   & A \injtens B &                                                                                                                                                &
      \end{tikzcd}
    }
  \end{equation*}
The upper square commutes because $\fdOS$ is a BV-category (see \eqref{eq:assoc-w-seq}).
The two diamonds on the
sides commute by naturality of $v$ and the pentagon at the bottom by assumption on
$A$ and $B$ being von Neuman algebras. \\

\noindent Now for the vN-algebra specific coherences. We have that $i$ satisfies
the involution condition:
  \begin{equation*}
    \adjustbox{scale=0.55}{
      \begin{tikzcd}
        & A \injtens B \arrow[equal]{d} \arrow[equal]{rrrrr} \arrow[equal]{ldd}                       &                                                &                               &                               &                                                                        & A \injtens B                             &                                 \\
        & A_{cc} \injtens B_{cc} \arrow[r, "(i_A)_c \otimes (i_B)_c"] \arrow[equal]{d} & (A_o)_c \injtens (B_o)_c \arrow[equal]{rrr} \arrow[equal]{d} &                               &                               & A_{co} \injtens B_{co} \arrow[r, "(i_A)_o \otimes (i_B)_o"] \arrow[equal]{d} & A_{oo} \injtens B_{oo} \arrow[equal]{u}         &                                 \\
      (A \injtens B)_{cc} \arrow[equal]{r} & (A_c\injtens B_c)_c \arrow[r, "(i_A \otimes i_B)_c"']                  & (A_o \injtens B_o)_c \arrow[equal]{r}                 & ((A\injtens B)_o)_c \arrow[equal]{r} & (A \injtens B)_{co} \arrow[equal]{r} & (A_c \injtens B_c)_o \arrow[r, "(i_A \otimes i_B)_o"']                 & (A_o \injtens B_o)_o \arrow[equal]{u} \arrow[equal]{r} & (A \injtens B)_{oo} \arrow[equal]{luu}
      \end{tikzcd}
    }
  \end{equation*}

\noindent The involution $i$ is also reverse multiplicative:
  \begin{equation*}
    \adjustbox{scale=0.75}{
      \begin{tikzcd}
        (A \injtens B)_c \hagertimes (A \injtens B)_c \arrow[equal]{d} \arrow[equal]{r}                                                                & ((A \injtens B)\hagertimes (A \injtens B))_c \arrow[r, "v"]                                                                                    & ((A\hagertimes A) \injtens (B \hagertimes B))_c \arrow[equal]{d} \arrow[r, "(\mu_A \otimes \mu_B)_c"]              & (A \injtens B)_c \arrow[equal]{d}                              \\
        (A_c\injtens B_c) \hagertimes (A_c \injtens B_c) \arrow[r, "v"] \arrow[d, "(i_A\otimes i_B) \otimes (i_A \otimes i_B)"'] & (A_c\hagertimes A_c) \injtens (B_c \hagertimes B_c) \arrow[d, "(i_A\otimes i_A) \otimes (i_B \otimes i_B)"] \arrow[equal]{r} & (A \hagertimes A)_c \injtens (B \hagertimes B)_c \arrow[r, "(\mu_A)_c \otimes (\mu_B)_c "] & A_c\injtens B_c \arrow[d, "i_A \otimes i_B"] \\
        (A_o \injtens B_o) \hagertimes (A_o \injtens B_o) \arrow[equal]{d} \arrow[r, "v"']                                           & (A_o \hagertimes A_o) \injtens (B_o \hagertimes B_o)  \arrow[r, "\gamma \otimes \gamma"']                                         & (A \hagertimes A)_o \injtens (B \hagertimes B)_o \arrow[r, "\mu_{A_o} \otimes \mu_{B_o}"'] \arrow[equal]{d}  & A_o \injtens B_o \arrow[equal]{d}                         \\
        (A \injtens B)_o \hagertimes (A \injtens B)_o \arrow[r, "\gamma"']                                                               & ((A \injtens B) \hagertimes (A \injtens B))_o \arrow[r, "v"']                                                                                       & ((A\hagertimes A) \injtens (B \hagertimes B))_o \arrow[r, "(\mu_A \otimes \mu_B)_o"']                          & (A \injtens B)_o
        \end{tikzcd}
    }
  \end{equation*}

  The coherence diagrams for $A\linftyplus B$ are analogous as the binary product is
  automatically duoidal to any monoidal structure (see \cite[Example 6.19]{AM2010}).
  To give an indication of the proof, note that the interchange law for
  this duoidal structure is given by the universal property of the product:

  \begin{equation*}
    v' \eqdef \langle \pi_A \otimes \pi_C, \pi_B \otimes \pi_D \rangle :
    (A \linftyplus B) \hagertimes (C \linftyplus D) \to (A \hagertimes C) \linftyplus (B \hagertimes D)
  \end{equation*}
  and that the unit and multiplication for
  $A\linftyplus B$ can be factored in the following way:
  \begin{equation*}
    \adjustbox{scale=0.8}{
    \begin{tikzcd}
      \mathbb{C} \arrow[rrdd, "{\langle \eta_A, \eta_B \rangle}"'] \arrow[rr, "\langle \mathrm{id} {,} \mathrm{id} \rangle"] &  & \mathbb{C} \linftyplus \mathbb{C} \arrow[dd, "\eta_A \oplus \eta_B"] & & (A \linftyplus B)\hagertimes (A \linftyplus B) \arrow[rrrdd, "{\langle \mu_A \circ (\pi_A \otimes \pi_A), \mu_B \circ (\pi_B \otimes \pi_B) \rangle}"'] \arrow[rrr, "v'"] & & & (A \hagertimes A) \linftyplus (B \hagertimes B) \arrow[dd, "\mu_A \oplus \mu_B"] \\
                                                                                                                          &  &                                                                       & &                                                                                                                                                                                                                                 & & &                                                                                 \\
                                                                                                                          &  & A \linftyplus B                                                       & &                                                                                                                                                                                                                                 & & & A \linftyplus B
      \end{tikzcd}}
  \end{equation*}

  For the C$^*$-identity of the product, note that
  given f.d operator spaces $A$ and $B$ and a c.c map
  $u = \langle a, b\rangle : \mathbb{C} \to A \linftyplus B$ then $u$ is a
  c.i morphism if and only if $a$ or $b$ is a c.i map.
  Let $u =  \langle a, b\rangle : \mathbb{C} \to A \linftyplus B$ be a c.i. map,
  we know by naturality and the universal property of the product that
  the following diagram commutes:
   \begin{equation*}
    \adjustbox{scale=0.7}{
      \begin{tikzcd}
        \mathbb{C} \arrow[d, "\cong"'] \arrow[rr, "\langle \mathrm{id} {,} \mathrm{id} \rangle"] & & \mathbb{C} \linftyplus \mathbb{C} \arrow[d, "\cong"] \\
        \mathbb{C}_o \hagertimes \mathbb{C}_c \arrow[d, "{\langle a {,} b\rangle \otimes \langle a {,} b\rangle}"'] \arrow[rr, "\langle \mathrm{id} {,} \mathrm{id} \rangle"] &  & (\mathbb{C}_o \hagertimes \mathbb{C}_c) \linftyplus (\mathbb{C}_o \hagertimes \mathbb{C}_c) \arrow[dd, "(a\otimes a) \oplus (b \otimes b)"] \\
        (A \linftyplus B)_o \hagertimes (A \linftyplus B)_c \arrow[equal]{d}                                                           &  &                                                                                                                                             \\
        (A_o \linftyplus B_o) \hagertimes (A_c \linftyplus B_c) \arrow[d, "\mathrm{id} \otimes (i_A \oplus i_B)"'] \arrow[rr, "v'"]                    &  & (A_o \hagertimes A_c) \linftyplus (B_o \hagertimes B_c) \arrow[d, "(\mathrm{id} \otimes i_A) \oplus (\mathrm{id}\otimes i_B)"]              \\
        (A_o \linftyplus B_o) \hagertimes (A_o \linftyplus B_o) \arrow[equal]{d} \arrow[rr, "v'"]                                       &  & (A_o \hagertimes A_o) \linftyplus (B_o \hagertimes B_o) \arrow[d, "\gamma \oplus \gamma"]                                                  \\
        ((A \linftyplus B)_o \hagertimes (A \linftyplus B))_o \arrow[d, "\gamma"']                                              &  &   (A \hagertimes A)_o \linftyplus (B \hagertimes B)_o \arrow[equal]{d}                                                                                                                                         \\
        ((A \linftyplus B) \hagertimes (A \linftyplus B))_o \arrow[d, "v'"'] \arrow[rr, "v'"]                                     &  & ((A \hagertimes A) \linftyplus (B \hagertimes B))_o \arrow[d, "\mu_A \oplus \mu_B"]                                                         \\
        ((A \hagertimes A) \linftyplus (B \hagertimes B))_o \arrow[rr, "\mu_A \oplus \mu_B"']   &  & (A \linftyplus B)_o                                                                                                                \\
        \end{tikzcd}}
   \end{equation*}
  \noindent
  Hence it suffices for either
  \begin{equation*}
    \mathbb{C} \cong \mathbb{C}_o \hagertimes \mathbb{C}_c
    \xrightarrow{a \otimes a} A_o \hagertimes A_c
    \xrightarrow{A_o \otimes i} A_o \hagertimes A_o
    \xrightarrow{\gamma} (A \hagertimes A)_o
    \xrightarrow{\mu} A_o
  \end{equation*}
  or
  \begin{equation*}
    \mathbb{C} \cong \mathbb{C}_o \hagertimes \mathbb{C}_c
    \xrightarrow{b \otimes b} B_o \hagertimes B_c
    \xrightarrow{B_o \otimes i} B_o \hagertimes B_o
    \xrightarrow{\gamma} (B \hagertimes B)_o
    \xrightarrow{\mu} B_o
  \end{equation*}
  to be a c.i map. We know that either $a$ or $b$ is c.i, hence by assumption
  we have that the C$^*$-identity holds for $A \linftyplus B$.

  For $C^*$-identity of the tensor, recall that any f.d vN-algebra is isomorphic
  to a direct sum of the following form: $\linftyplus_{1 \leq i \leq n} M_{k_i}$.
  Now since $\itimes$ distributes over $\linftyplus$, it is sufficient to prove
  that the $C^*$-identity holds on $M_n \itimes M_m \cong M_{nm}$.
  It is straightforward to check that the c.i.i. $M_n \itimes M_m \cong M_{nm}$
  is multiplicative, unital and preserves the involution.
  Suppose that $u : \mathbb{C} \to M_n \itimes M_m$ is a c.i. map. Then
  $r \defeq \left( \mathbb{C} \xrightarrow{u} M_n \itimes M_m \cong M_{nm} \right)$ is a c.i. map too.
  The following diagram commutes:

  \begin{equation*}
    \adjustbox{scale=0.7}{
    \begin{tikzcd}
      \mathbb{C} \arrow[d, "\cong"'] \arrow[equal] {rr}                                                                                                 &  &    \mathbb{C} \arrow[d, "\cong"]                                     \\
      \mathbb{C}_o \hagertimes \mathbb{C}_c \arrow[d, "u \otimes u"']  \arrow[equal]{rr}                                                                      &  &  \mathbb{C}_o \hagertimes \mathbb{C}_c \arrow[d, "r \otimes r"]               \\
      (M_n \itimes M_m)_o \hagertimes (M_n \itimes M_m)_c \arrow[rr, "\cong"] \arrow[d, "\mathrm{id}\otimes (i \otimes i)"'] &  & (M_{nm})_o \hagertimes (M_{nm})_c \arrow[d, "\mathrm{id} \otimes i"] \\
      (M_n \itimes M_m)_o \hagertimes (M_n \itimes M_m)_o \arrow[rr, "\cong"] \arrow[d, "\gamma"']                           &  & (M_{nm})_o \hagertimes (M_{nm})_o \arrow[d, "\gamma"]                          \\
      ((M_n \itimes M_m) \hagertimes (M_n \itimes M_m))_o \arrow[rr, "\cong"] \arrow[d, "\mu \otimes \mu"']                  &  & ((M_{nm}) \hagertimes (M_{nm}))_o \arrow[d, "\mu"]                   \\
      (M_n \itimes M_m)_o \arrow[rr, "\cong"]                                                                                &  & (M_{nm})_o
      \end{tikzcd}}
  \end{equation*}

We know that $M_{nm}$ is a vN-algebra, hence the right side of the
diagram is a c.i map. Therefore the left side must be c.i too, because the horizontal arrows are c.i.i.
Therefore $M_n \itimes M_m$ also satisfies the C$^*$-identity.

\qed
\end{proof}

\vNsmc*

\begin{proof}
  We know from \cref{prop:tensor-alg}
  that the objects in $\vNAlg$ are closed under the injective tensor.
  Suppose that $f: A \to C, g : B \to D$ are morphisms in $\vNalg$.
  The following diagrams:

  \begin{equation*}
    \adjustbox{scale=0.8}{
      \begin{tikzcd}
        (A \itimes B) \hagertimes (A \itimes B) \arrow[d, "v"'] \arrow[rr, "(f \otimes g) \otimes (f \otimes g)"]                        &  & (C \itimes D) \hagertimes (C \itimes D) \arrow[d, "v"]                       &  &                                        & \mathbb{C} \arrow[d, "\cong "]                                                                         &             \\
        (A \hagertimes A) \itimes (B \hagertimes B) \arrow[d, "\mu_A \otimes \mu_B"'] \arrow[rr, "(f \otimes f) \otimes (g \otimes g)"'] &  & (C \hagertimes C) \itimes (D \hagertimes D) \arrow[d, "\mu_C \otimes \mu_D"] &  &                                        & \mathbb{C} \itimes \mathbb{C} \arrow[ld, "\eta_A \otimes \eta_B"'] \arrow[rd, "\eta_C \otimes \eta_D"] &             \\
        A \itimes B \arrow[rr, "f \otimes g"']                                                                                           &  & C \itimes D                                                                  &  & A \itimes B \arrow[rr, "f \otimes g"'] &                                                                                                        & C \itimes D
        \end{tikzcd}
    }
  \end{equation*}

  \begin{equation*}
    \adjustbox{scale=0.8}{
      \begin{tikzcd}
        (A \itimes B)_c \arrow[equal, d] \arrow[rr, "(f\otimes g)_c"]                      &  & (C \itimes D)_c \arrow[equal, d]                    \\
        A_c \itimes B_c \arrow[d, "i_A \otimes i_B"'] \arrow[rr, "f_c \otimes g_c"] &  & C_c \itimes D_c \arrow[d, "i_C \otimes i_D"] \\
        A_o \itimes B_o \arrow[equal]{d} \arrow[rr, "f_o \otimes g_o"]                     &  & C_o \itimes D_o \arrow[equal]{d}                    \\
        (A \itimes B)_o \arrow[rr, " (f \otimes g)_o"]                              &  & (C \itimes D)_o
        \end{tikzcd}
    }
  \end{equation*}
  commute, so the functor
  $\itimes : \vNAlg \times \vNAlg \to \vNAlg$ is well-defined.
  Next, we need to verify that the associator and unitors of $\itimes$ are
  indeed morphisms in $\vNAlg$.

  We begin with the associator. Let $A$, $B$ and $C$ be f.d vN-algebras and $\alpha$
  be the associator of $\itimes$.
  To show that $\alpha$ is multiplicative, consider the following diagram.

  \begin{equation*}
    \adjustbox{scale=0.8}{
    \begin{tikzcd}
      ((A\itimes B) \itimes C) \hagertimes ((A\itimes B)\itimes C) \arrow[rr, "\alpha \otimes \alpha"] \arrow[d, "(v \otimes \mathrm{id}) \circ v"'] &  & (A \itimes (B \itimes C)) \hagertimes (A \itimes (B \itimes C)) \arrow[d, "(\mathrm{id} \otimes v) \circ v"]         \\
      ((A \hagertimes A) \itimes  (B \hagertimes B)) \itimes (C\hagertimes C) \arrow[d, "(\mu_A \otimes \mu_B)\otimes \mu_C"'] \arrow[rr, "\alpha"']       &  & (A \hagertimes A) \itimes  ((B \hagertimes B) \itimes (C \hagertimes C)) \arrow[d, "\mu_A \otimes (\mu_B\otimes \mu_C)"]                                                                 \\
      (A \itimes B) \itimes C \arrow[rr, "\alpha"']                                                                                                        &  & A \itimes (B \itimes C)
      \end{tikzcd}}
  \end{equation*}
  The upper square commutes as $\fdOS$ is a BV-category and the lower one by
  naturality of $\alpha$. The two diagrams below commute by naturality of $\alpha$ (modulo rewriting some equalities).
  \begin{equation*}
    \adjustbox{scale=0.8}{
      \begin{tikzcd}
        & \mathbb{C} \arrow[ld, "\cong"'] \arrow[rd, "\cong"] &                                                                                                        &  & ((A \itimes B) \itimes C)_c \arrow[equal]{d} \arrow[rr, "\alpha_c"]                             &  & (A \itimes (B \itimes C))_c \arrow[equal]{d}                              \\
              (\mathbb{C} \itimes \mathbb{C}) \itimes \mathbb{C} \arrow[d, "(\eta_A \otimes \eta_B )\otimes \eta_C"'] \arrow[rr, "\alpha"] &                                                     & \mathbb{C} \itimes (\mathbb{C} \itimes \mathbb{C}) \arrow[d, "\eta_A \otimes (\eta_B \otimes \eta_C)"] &  & (A_c \itimes B_c) \itimes C_c \arrow[d, "(i\otimes i ) \otimes i"'] \arrow[rr, "\alpha"] &  & A_c \itimes (B_c \itimes C_c) \arrow[d, "i \otimes (i \otimes i)"] \\
              (A \itimes B) \itimes C \arrow[rr, "\alpha"']                                                                                &                                                     & A \itimes (B \itimes C)                                                                                &  & (A_o \itimes B_o) \itimes C_o \arrow[rr, "\alpha"'] \arrow[equal]{d}                           &  & A_o \itimes (B_o \itimes C_o) \arrow[equal]{d}                            \\
        &                                                     &                                                                                                        &  & ((A \itimes B) \itimes C)_o \arrow[rr, "\alpha_o"']                                      &  & (A \itimes (B \itimes C))_o
\end{tikzcd}}
  \end{equation*}

\noindent
  This shows that the associator $\alpha$ is also a natural transformation in $\vNAlg.$
  We proceed similarily for the left unitor $\lambda$.
  \begin{equation*}
    \adjustbox{scale=0.8}{
    \begin{tikzcd}
      (\mathbb{C} \itimes A) \hagertimes (\mathbb{C} \itimes A) \arrow[d, "v"'] \arrow[rr, "\lambda \otimes \lambda"] &                                                            & A \hagertimes A \arrow[d]        \\
      (\mathbb{C} \hagertimes \mathbb{C}) \itimes (A \hagertimes A) \arrow[d, "\mu_\mathbb{C} \otimes \mu_A"'] \arrow[r, "\mu_\mathbb{C} \otimes \mathrm{id}"']              &  \mathbb{C} \itimes (A \hagertimes A) \arrow[r, "\lambda"] & A \hagertimes A \arrow[d, "\mu_A"] \\
      \mathbb{C} \itimes A \arrow[rr, "\lambda"']                                                                     &                                                            & A
      \end{tikzcd}}
  \end{equation*}
  The upper square commutes by the BV-assumption, and the lower one by naturality of $\lambda$.
  Note that $\mu_\mathbb{C} = \lambda_\mathbb{C}$. We proceed similarily for unitality and
  involutivity:

  \begin{equation*}
    \adjustbox{scale=0.8}{
      \begin{tikzcd}
        & \mathbb{C} \arrow[ld, "\cong"'] \arrow[rd, "\cong"] &                                &  & (\mathbb{C} \itimes A)_c \arrow[equal]{d} \arrow[rrd, "(\lambda_{A})_c"] &  &                      \\
\mathbb{C} \itimes \mathbb{C} \arrow[d, "\mathbb{C} \otimes  \eta_A"'] \arrow[rr, "\lambda_\mathbb{C}"] &                                                     & \mathbb{C} \arrow[d, "\eta_A"] &  & \mathbb{C}_c \itimes A_c \arrow[d, "i \otimes i_A"']   \arrow[r, "i \otimes A_c"']           &  \mathbb{C} \itimes A_c \arrow[r , "\lambda_{A_c}"'] & A_c \arrow[d, "i_A"] \\
\mathbb{C} \itimes A \arrow[rr, "\lambda_A"']                                                           &                                                     & A                              &  & \mathbb{C} \itimes A_o \arrow[rr, "\lambda_{A_o}"'] \arrow[equal]{d}     &  & A_o                  \\
        &                                                     &                                &  & (\mathbb{C} \itimes A)_o \arrow[rru, "(\lambda_A)_o"']            &  &
\end{tikzcd}
    }
  \end{equation*}
  This shows that the left unitor $\lambda$ is also a natural transformation in $\vNAlg.$
  We can show that the right unitor $\rho$ is also a natural transformation in $\vNAlg$ using similar arguments.
  Clearly the symmetries are in $\vNAlg$ as they are natural and commute with $v$.

  Finally, in order to see that $\vNAlg$ is symmetric monoidal, we note that the inclusion functor
  $\vNAlg \hookrightarrow \FdOS$ is faithful and preserves the tensor $\itimes$ strictly, so the necessary
  coherence diagrams also commute in $\vNAlg$.
  \qed
\end{proof}

\begin{lemma}
  \label{Lemma:star-hom-cp}
  Let $A,B \in \Ob(\vNAlg)$ and let $f:A \to B$ be a morphism in $\vNAlg$. Then $f$ is a CPU map
  in the sense of \cref{sec:vNAlg}.
\end{lemma}

\begin{proof}
  Let $f \colon A \to B$ be a morphism of $\vNAlg$. This map is unital by assumption, so we only
  have to show it is completely positive.
  Suppose that $p:\mathbb{C} \to M_n \injtens A$ is a positive
  element of $M_n \injtens A$, and let $a \colon \mathbb{C} \to M_n \injtens A$ be such that
  \begin{equation*}
    \adjustbox{scale=0.7}{
    \begin{tikzcd}
      \mathbb{C}  \arrow[equal]{dd} \arrow[r, "\cong"]   & \mathbb{C}_o \hagertimes \mathbb{C}_c \arrow[r, "a \otimes a"] & (M_n \itimes A)_o \hagertimes (M_n \itimes A)_c \arrow[r, "\mathrm{id} \otimes i"] & (M_n \itimes A)_o \hagertimes (M_n \itimes A)_o  \arrow[r, "\gamma"] & ((M_n \itimes A) \hagertimes (M_n \itimes A))_o \arrow[d, "v"] \\
      & & & & ((M_n \hagertimes M_n) \itimes (A \hagertimes A))_o \arrow[d, "\mu"] \\
      \mathbb{C}_o \arrow[rrrr, "p"'] &                                                                    &                                               &                                          & M_n \itimes A_o
    \end{tikzcd}}
  \end{equation*}
  Recall that $M_n$ is a vN-algebra, so by \cref{prop:vN-monoidal}, the map
  $M_n \otimes f : M_n \itimes A \to M_n \itimes B$ is multiplicative and
  preserves the involution, hence the following diagram commutes in $\fdOS:$
  \begin{equation*}
    \adjustbox{scale=0.55}{
      \begin{tikzcd}
                                                                       & (M_n \injtens B)_o \hagertimes (M_n \injtens B)_c  \arrow[rr, "\mathrm{id} \otimes (i_{M_n} \otimes i_B)"]                                                          &                                     & (M_n \injtens B)_o \hagertimes (M_n \injtens B)_o \arrow[r, "\gamma"]                                                        & ((M_n \injtens B) \hagertimes (M_n \injtens B))_o \arrow[equal]{rrdd}                                                                          &                            &                                                                                          \\
                                                                       &                                                                                                                                                                   & &                                                                                                                            &                                                                                                                                                                               & &                                                                                          \\
        \mathbb{C}_o \hagertimes \mathbb{C}_c \arrow[r, "a \otimes a"] & (M_n \injtens A)_o \hagertimes (M_n \injtens A)_c  \arrow[rr, "\mathrm{id} \otimes (i_{M_n} \otimes i_A)"'] \arrow[uu, "(M_n \otimes f) \otimes (M_n \otimes f)"] & & (M_n \injtens A)_o \hagertimes (M_n \injtens A)_o \arrow[r, "\gamma"'] \arrow[uu, "(M_n \otimes f) \otimes (M_n \otimes f)"] & ((M_n \injtens A) \hagertimes (M_n \injtens A))_o \arrow[d, "v"'] \arrow[rr, "(M_n \otimes f) \otimes (M_n \otimes f)"'] \arrow[uu, "(M_n \otimes f) \otimes (M_n \otimes f)"] & & ((M_n \injtens B) \hagertimes (M_n \injtens B))_o \arrow[d, "v"]                         \\
                                                                       &                                                                                                                                                                   & &                                                                                                                            & ((M_n \hagertimes M_n) \injtens (A \hagertimes A))_o \arrow[d, "\mu_{M_n} \otimes \mu_A"']                              & & ((M_n \hagertimes M_n) \injtens (B \hagertimes B))_o \arrow[d, "\mu_{M_n} \otimes \mu_B"] \\
        \mathbb{C}_o \arrow[rrrr, "p"'] \arrow[uu, "\cong"]             &                                                                                                                                                                  & &                                                                                                                            & (M_n \injtens A)_o \arrow[rr, "M_n \otimes f"']                                                                                                                                            & & (M_n \injtens B)_o
      \end{tikzcd}}
  \end{equation*}
  The outer edges of the diagram give us a factorization of
  $(M_n \otimes f) \circ p$ into the product of $(M_n \otimes f) \circ a$ and
  its involution.
  \qed
\end{proof}

\Hsmcc*

\begin{proof}
We know from \cref{prop:tensor-alg}
that the objects in $\HH$ are closed under the injective tensor, and that the
tensor of unital maps is unital. By functoriality of the injective tensor with
respect to c.c. maps and as c.p.u maps are equivalent to c.c.u maps by \cref{prop:alg-cc-cp},
the injective tensor product is indeed a functor $\itimes: \HH \times \HH \to \HH$.
The proof of \cref{prop:vN-monoidal} shows that that the associator, left/right unitors, and symmetry of $\itimes$ are unital maps, and as they are c.c., they are also c.p. maps.

To see that $\HH$ is symmetric monoidal w.r.t $\itimes,$ note that the inclusion functor
$\HH \hookrightarrow \FdOS$ is faithful and preserves $\itimes$ strictly, so the necessary
coherence diagrams also commute in \HH.

Finally, it follows by construction that the inclusion functor strictly preserves the symmetric monoidal
structure and products.
\qed
\end{proof}

\subsection{Duality between von Neumann algebras and von Neumann coalgebras}
Let us recall the definition of a vN-coalgebra.
\begin{definition}
  A f.d. \emph{vN-coalgebra} is an operator space $C$ together with
  \begin{enumerate}
    \item a complete contraction $\delta \colon C \to C \hagertimes C$, called \emph{comultiplication};
    \item a complete contraction $\varepsilon \colon C \to \mathbb C$, called \emph{counit};
    \item a complete contraction $j \colon C_o \to C_c,$ called \emph{involution},
  \end{enumerate}
  such that
  \begin{itemize}
    \item $(C, \varepsilon, \delta)$ is a comonoid object in $(\FdOS, \mathbb C, \hagertimes)$;
    \item $j$ is an \emph{involution}, meaning that the following diagram commutes:
    \begin{equation}
      \tag{involution}
        \begin{tikzcd}
          C \arrow[equal]{rr}                                           &                                                                                   & C     \\
          C_{oo} \arrow[equal]{u} \arrow[r, "j_o"'] & C_{co} = C_{oc} \arrow[r, "j_c"'] & C_{cc} \arrow[equal]{u}
        \end{tikzcd}
    \end{equation}
    \item $j$ is \emph{reverse comultiplicative}, meaning that the following
    diagram commutes:
    \begin{equation}
      \tag{reverse comultiplicative}
      \begin{tikzcd}
        C_o \arrow[d, "\delta"'] \arrow[r, "j"]                 & C_c \arrow[r, "\delta"]                                     & (C\hagertimes C)_c \arrow[equal]{d} \\
        (C\hagertimes C)_o \arrow[r, "{\gamma}"'] & C_o\hagertimes C_o \arrow[r, "j\otimes j"'] & C_c\hagertimes C_c
      \end{tikzcd}
    \end{equation}
    \item for every complete quotient map (i.e. strong epi) $e: C \to \mathbb{C}$,
    the following composite in $\FdOS$ is
    a complete quotient map (i.e. strong epi)
    \begin{equation}
      \tag{co-C$^*$-id}
      C_o \xrightarrow{\delta}
      (C\hagertimes C)_o \xrightarrow{\gamma}
      C_o \hagertimes C_o \xrightarrow{C_o \otimes j}
      C_o\hagertimes C_c \xrightarrow{e \otimes e}
      \mathbb{C}_o \hagertimes \mathbb{C}_c \cong \mathbb{C} .
    \end{equation}
  \end{itemize}
\end{definition}

\begin{remark}
  The definition above is the same as the one in \cref{sec:vNCoalg} only
  slightly expanded to give the coherence conditions names.
\end{remark}

\vNdual*

\begin{proof}
  Given a $A : \vNAlg$, clearly $A^*$ is a comonoid.
  For the remaining structure we can easily check that the following diagrams
  commute.
We have that $j$ satisfies the involution coherence diagram:
\begin{equation*}
  \adjustbox{scale=0.7}{
    \begin{tikzcd}
      A^*                              &                                                   &                                                     &                        &                      &                                          &                                  & A^* \arrow[equal]{lllllll} \arrow[equal]{dd} \arrow[equal]{ld} \\
                                       & (A_{cc})^* \arrow[equal]{lu}                             & (A_{oc})^* \arrow[l, "(i_c)^*"]                     &                        &                      & (A_{co})^* \arrow[equal]{lll}                   & (A_{oo})^*     \arrow[l, "(i_o)^*"]   \arrow[equal]{d}                      &                                           \\
      (A^*)_{cc} \arrow[equal]{uu} \arrow[equal]{ru} & ((A_c)^*)_c \arrow[l, "\cong"] \arrow[u, "\cong"] & ((A_o)^*)_c \arrow[l, "(i^*)_c"] \arrow[u, "\cong"] & ((A^*)_{oc}) \arrow[equal,]{l} & (A^*)_{co} \arrow[equal]{l} & ((A_c)^*)_o \arrow[l, "\cong"] \arrow[equal]{u} & ((A_o)^*)_o \arrow[l, "(i^*)_o"] & (A^*)_{oo} \arrow[equal]{l}
      \end{tikzcd}
  }
\end{equation*}
It is also reverse comultiplicative:
\begin{equation*}
  \begin{tikzcd}
    (A^*)_c \hagertimes (A^*)_c                                       &                                                               & (A^*\hagertimes A^*)_c \arrow[equal]{ll}                      &                                                          & (A^*)_c \arrow[ll, "\delta_c"'] \arrow[lld, "(\mu^*)_c"]                                  \\
                                                                    & (A_c)^* \hagertimes (A_c)^* \arrow[lu, "\cong"]                          & ((A \hagertimes A)^*)_c \arrow[u, "\cong"]                     &                                                          &                                                                            \\
                                                                    & (A_c\hagertimes A_c)^* \arrow[u, "\cong"]                                & (({A\hagertimes A})_c)^* \arrow[equal]{l} \arrow[u, "\cong"] \arrow[u] & (A_c)^* \arrow[l, "(\mu_c)^*"']    \arrow[ruu, "\cong"']                   &                                                                            \\
                                                                    & (A_o \hagertimes A_o)^* \arrow[u, "(i \otimes i)^*"]        & ((A \hagertimes A)_o)^* \arrow[l, "{\gamma^*}"]    & (A_o)^* \arrow[l, "(\mu_o)^*"] \arrow[u, "i^*"'] &                                                                            \\
                                                                    & (A_o)^* \hagertimes (A_o)^* \arrow[u, "\cong"] \arrow[uuu, bend left=80, "i^* \otimes i^*"] & ((A \hagertimes A)^*)_o \arrow[equal]{u}                      &                                                          &                                                                            \\
    (A^*)_o \otimes (A^*)_o \arrow[equal]{ru} \arrow[uuuuu, "j \otimes j"] &                                                               & (A^* \hagertimes A^*)_o \arrow[ll, "\gamma"] \arrow[u, "\cong"] &                                                          & (A^*)_o \arrow[ll, "\delta_o"] \arrow[uuuuu, "j"'] \arrow[equal]{luu} \arrow[llu,"(\mu^*)_o"]
    \end{tikzcd}
\end{equation*}
  C$^*$-identity: Let $a : A^* \to \mathbb{C}$ be a c.q. map.
  By \cref{rem:ci-dual-cq}, $a^*: \mathbb{C} \to A$ is a c.i. map.
  We know that the following diagram commutes:
  \begin{equation*}
    \adjustbox{scale=0.75}{
      \begin{tikzcd}
        (A_o)^* \arrow[equal]{ddd}  \arrow[r, "(\mu_o)^*"]          & ((A \hagertimes A)_o)^* \arrow[r, "{\gamma^*}"] \arrow[equal]{d} & (A_o \hagertimes A_o)^* \arrow[r, "(A_o \otimes i)^*"] \arrow[dd, "\cong"] & (A_o\hagertimes A_c)^* \arrow[rr, "((a^*)_c \otimes (a^*)_o)^*"] \arrow[dd, "\cong"] & & ((\mathbb{C}^*)_o \hagertimes (\mathbb{C}^*)_c)^* \arrow[dd, "\cong"] \arrow[r, "\cong"] & \mathbb{C} \arrow[equal]{dd} \\
                                       & ((A \hagertimes A)^*)_o \arrow[dd, "\cong"]                                            &                                                                                               &                                                                                                        & &                                                          &  \\
                                       &                                                                              & (A_o)^* \hagertimes (A_o)^* \arrow[equal]{d} \arrow[r, "(A_o)^* \otimes i^*"']                                                 & (A_o)^* \hagertimes (A_c)^* \arrow[rr, "((a^*)_o)^* \otimes ((a^*)_c)^*"'] \arrow[d, "\cong"]     & & ((\mathbb{C}^*)_o)^* \hagertimes ((\mathbb{C}^*)_c)^* \arrow[d, "\cong"] \arrow[r, "\cong"] & \mathbb{C} \arrow[equal]{d} \\
        (A^*)_o \arrow[r, "\delta_o"'] & (A^*\hagertimes A^*)_o \arrow[r, "\gamma"']                                             & (A^*)_o \hagertimes (A^*)_o \arrow[r, "(A^*)_o \otimes j"']                                                            & (A^*)_o \hagertimes (A^*)_c \arrow[rr, "a_o \otimes a_c"']                                      & & \mathbb{C}_o \hagertimes \mathbb{C}_c \arrow[r, "\cong"']         & \mathbb{C}
      \end{tikzcd}}
  \end{equation*}
  As $A$ is a vN-algebra, the top morphism is the dual of a c.i. map, i.e. a c.q. map,
  thus the bottom morphism is a c.q. map too.

We can prove that the dual of a vN-coalgebra is a vN-algebra in a similar way.
\qed
\end{proof}

\begin{proposition}
  \label{prop:equiv-star-hom}
  We have the following duality between the morphisms in $\vNAlg$ and
  $\vNCoalg$ (dual taken in $\fdOS$):
  \begin{itemize}
    \item Suppose that $f : A \to B \in \vNAlg$ then $f^*: B^* \to A^* \in \vNCoalg$;
    \item Suppose that $g : C \to D \in \vNCoalg$ then $g^*: D^* \to C^* \in \vNAlg$.
  \end{itemize}
\end{proposition}

\begin{proof}
  We prove the first statement as the second one is proved completely analogously.
  Let $f : A \to B \in \vNAlg$, then the following diagrams commute
  \begin{equation*}
    \adjustbox{scale=0.8}{
  \begin{tikzcd}
      B^* \hagertimes B^* \arrow[r, "f^* \otimes f^*"] & A^* \hagertimes A^*                     & (B^*)_c \arrow[r, "(f^*)_c"]                                 & (A^*)_c                                 & \mathbb{C} \arrow[equal]{r}            & \mathbb{C}    \\
      (B \hagertimes B)^* \arrow[u, "\cong"]  \arrow[r, "(f \otimes f)^*"]         & (A \hagertimes A)^* \arrow[u, "\cong"'] & (B_c)^* \arrow[u, "\cong"] \arrow[r, "(f_c)^*"]              & (A_c)^* \arrow[u, "\cong"']             & \mathbb{C}^* \arrow[u, "\cong"] \arrow[equal]{r}          & \mathbb{C}^* \arrow[u, "\cong"']   \\
      B^* \arrow[u, "(\mu_B)^*"] \arrow[r, "f^*"']     & A^* \arrow[u, "(\mu_A)^*"']             & (B^*)_o = (B_o)^* \arrow[r, "(f^*)_o"'] \arrow[u, "(i_B)^*"] & (A_o)^* = (A^*)_o \arrow[u, "(i_A)^*"'] & B^* \arrow[r, "f^*"'] \arrow[u, "(\eta_B)^*"] & A^* \arrow[u, "(\eta_A)^*"']
      \end{tikzcd}}
  \end{equation*}
  The lower squares all commute by assumption, and upper ones by naturality.
  \qed
\end{proof}

\begin{proposition}
  \label{prop:dual-functors-vNalg}
  The dual on vN-algebras and vN-coalgebras as defined in \cref{prop:duals}
  extends to functors
  \begin{equation*}
    (-)^* \colon \vNcoalg \rightleftarrows \vNalg^\mathrm{op} \colon (-)^*.
  \end{equation*}
\end{proposition}

\begin{proof}
  \cref{prop:equiv-star-hom} assures us that both functors are well-typed.
  \qed
\end{proof}

\begin{lemma}
  \label{lemma:dd-vNalg}
  The functors
  $(-)^* \colon \vNcoalg \rightleftarrows \vNalg^\mathrm{op} \colon (-)^*$
  give rise to:
  \begin{itemize}
    \item a natural isomorphism $d: \mathrm{Id} \cong (-)^{**}: \vNAlg \to \vNAlg$;
    \item a natural isomorphism $d: \mathrm{Id} \cong (-)^{**}: \vNCoalg \to \vNCoalg$;
  \end{itemize}
  and both of these isomorphisms coincide with the double-dual natural isomorphism $d \colon \mathrm{Id} \cong (-)^{**} \colon \FdOS \to \FdOS$
  under the obvious forgetful functors.
\end{lemma}

\begin{proof}
  We prove the first statement as the second one is proved completely analogously.
  Let $A\in \vNAlg$, it is sufficient to prove that the isomorphism
  $d_A : A \to A^{**}$ indeed is in $\vNAlg$, as naturality follows from it being
  natural on $\fdOS$. Note that the $d_A$ commutes with the Haagerup tensor and conjugation
  in the following way:
  \begin{equation}
    \label{eq:dd-alg}
    \begin{tikzcd}
      A \hagertimes A \arrow[d, "d"'] \arrow[r, "d \otimes d"] & A^{**} \hagertimes A^{**} \arrow[d, "\cong"] & A_c \arrow[d, "d"'] \arrow[r, "(d)_c"] & (A^{**})_c \arrow[d, "\cong"] \\
      (A \hagertimes A)^{**}                                                         & (A^* \hagertimes A^*)^* \arrow[l, "\cong"]  & (A_c)^{**}           & ((A^*)_c)^* \arrow[l, "\cong"]
    \end{tikzcd}
  \end{equation}
  this is straightforward to check when reducing to elementary tensors.
  Having these diagrams commute now reduces the rest of the proof to
  leveraging the naturality of $d$, see diagrams below.

  We have that $d$ is multiplicative:
  \begin{equation*}
    \begin{tikzcd}
      A \hagertimes A \arrow[rr, "d \otimes d"] \arrow[ddd, "\mu_A"'] \arrow[rrdd, "d"'] &  &  A^{**} \hagertimes A^{**} \arrow[d, "\cong"]                                               \\
                                                                                                              &                        & (A^* \hagertimes A^*)^* \arrow[d, "\cong"]        \\
                                                                                                              &                                  & (A \hagertimes A)^{**} \arrow[d, "(\mu_A)^{**}"]  \\
      A \arrow[rr, "d"']                                                                                     &                                         & A^{**} &
    \end{tikzcd}
  \end{equation*}

  involutive:
  \begin{equation*}
    \begin{tikzcd}
      A_c  \arrow[ddd, "i"'] \arrow[rr, "(d)_c"] \arrow[rrdd, "d"'] & & (A^{**})_c \arrow[d, "\cong"]     \\
                                                                             & & ((A^*)_c)^* \arrow[d, "\cong"] \\
                                                                             &  &  (A_c)^{**} \arrow[d, "i^{**}"]    &   \\
      A_o \arrow[rr, "(d)_o"']                                             & & (A_o)^{**} = (A^{**})_o                    &
      \end{tikzcd}
  \end{equation*}

  and unital:
  \begin{equation*}
    \begin{tikzcd}
      \mathbb{C} \arrow[dd, "\eta_A"'] \arrow[rr, "d"] & & \mathbb{C}^{**} \arrow[dd, "(\eta_A)^{**}"] \\
      & & \\
      A \arrow[rr, "d"']                                         & & A^{**}
    \end{tikzcd}
  \end{equation*}

  For the second statement, note that the following diagrams commute:
  \begin{equation}
    \label{eq:dd-coalg}
    \begin{tikzcd}
      C \hagertimes C \arrow[d, "d \otimes d"'] \arrow[r, "d"] & (C \hagertimes C)^{**} \arrow[d, "\cong"]  & C_c \arrow[d, "(d)_c"'] \arrow[r, "d"] & (C_c)^{**} \arrow[d, "\cong"]  \\
      C^{**} \hagertimes C^{**}                                                      & (C^* \hagertimes C^*)^* \arrow[l, "\cong"] & (C^{**})_c                               & ((C^*)_c)^* \arrow[l, "\cong"]
      \end{tikzcd}
  \end{equation}
  The rest of the proof is analogous.
  \qed
\end{proof}

\begin{theorem}
  \label{equiv-alg-coalg}
  We have an equivalence of categories $(-)^* \colon \vNcoalg \simeq \vNalg^\mathrm{op} \colon (-)^*$,
  where the action on objects of $(-)^*$ is defined as in \cref{prop:duals} and on morphisms
  in the usual way, i.e. as in $\FdOS.$
\end{theorem}

\begin{proof}
  Follows by combining \cref{prop:dual-functors-vNalg} and \cref{lemma:dd-vNalg}.
  \qed
\end{proof}

\subsection{von Neumann coalgebras}

Most of this section will be dedicated to prove the following proposition:

\vNcoalgtens*

In order to prove this we want to leverage the equivalence between
$\vNcoalg \cong \vNalg^\text{op}$ proved in the previous section, thus we need the
following two lemmas.

\begin{lemma}
  \label{lemma:iso-proj-inj}
  The c.i. natural isomorphism:
  \begin{equation}
    \label{eq:isotens}
    \theta : C \projtens D \cong (C^* \injtens D^*)^*
  \end{equation}
  in $\fdOS$, preserves the comonoid structure and the involution as defined above
  for $C \projtens D$.
\end{lemma}

\begin{proof}
  We can verify that the following diagrams commute. \\

  \noindent The c.i. natural isomorphism $\theta$ is comultiplicative:
  \begin{equation*}
    \adjustbox{scale=0.8}{
    \begin{tikzcd}
      C \projtens D \arrow[rr, "\theta"] \arrow[d, "\delta_C \otimes \delta_D"']           &  & (C^* \injtens D^*)^* \arrow[d, "((\delta_C)^* \otimes (\delta_D)^*)^*"]  &                                                                             \\
      (C \hagertimes C) \projtens (D \hagertimes D) \arrow[rr, "\theta"] \arrow[ddd, "w"'] &  & ((C \hagertimes C)^* \injtens (D \hagertimes D)^*)^* \arrow[rd, "\cong"] &                                                                             \\
                                                                                           &  &                                                                          & ((C^* \hagertimes C^*) \injtens (D^* \hagertimes D^*))^* \arrow[d, "v^*"]   \\
                                                                                           &  &                                                                          & ((C^* \injtens D^*) \hagertimes  (C^* \injtens D^*))^* \arrow[ld, "\cong"] \\
      (C \projtens D) \hagertimes (C \projtens D) \arrow[rr, "\theta \otimes \theta"']     &  & (C^* \injtens D^*)^* \hagertimes  (C^* \injtens D^*)^*                   &
      \end{tikzcd}}
  \end{equation*}
  The upper square commutes by naturality, and the lower one by the duality between
  the two interchange maps $v$ and $w$ in a BV-category with negation. \\

  \noindent It is counital:
  \begin{equation*}
    \begin{tikzcd}
      C \projtens D \arrow[rr, "\theta"] \arrow[d, "\varepsilon_C \otimes \varepsilon_D"'] &            & (C^* \injtens D^*)^* \arrow[d, "((\varepsilon_C)^* \otimes (\varepsilon_D)^*)^*"] \\
      \mathbb{C} \projtens \mathbb{C} \arrow[rr, "\theta"] \arrow[rd, "\cong"']            &            & (\mathbb{C}^* \injtens \mathbb{C}^*)^* \arrow[ld, "\cong"]                        \\
                                                                                           & \mathbb{C} &
      \end{tikzcd}
  \end{equation*}
  \noindent And it preserves the involution:
  \begin{equation*}
    \begin{tikzcd}
      (C \projtens D)_o \arrow[rr, "(\theta)_o"] \arrow[equal]{d}                     &  & ((C^* \injtens D^*)^*)_o \arrow[equal]{d}                                    \\
      C_o \projtens D_o \arrow[rr, "\theta"] \arrow[d, "j_C \otimes j_D"'] &  & ((C_o)^* \injtens (D_o)^*)^* \arrow[d, "((j_C)^* \otimes (j_D)^*)^*"] \\
      C_c \projtens D_c \arrow[rr, "\theta"] \arrow[equal]{ddd}                   &  & ((C_c)^* \injtens (D_c)^*)^* \arrow[d, "\cong"]                       \\
                                                                          &  & ((C^*)_c \injtens (D^*)_c)^* \arrow[equal]{d}                                \\
                                                                          &  & ((C^* \injtens D^*)_c)^* \arrow[d, "\cong"]                           \\
      (C \projtens D)_c \arrow[rr, "(\theta)_c"']                              &  & ((C^* \injtens D^*)^*)_c
      \end{tikzcd}
  \end{equation*}
  We can check that the lower square commutes by reducing to elementary tensors and
  leveraging that the double dual commutes with
  conjugation (see \eqref{eq:dd-coalg}). \\
  \qed
\end{proof}

\begin{lemma}
  \label{lemma:iso-one-infty}
  The c.i. natural isomorphism:
  \begin{equation}
    \label{eq:iso}
    C \loneplus D \cong (C^* \linftyplus D^*)^*
  \end{equation}
  in $\fdOS$, preserves the comonoid structure and the involution as defined above
  for $C \loneplus D$.
\end{lemma}

\begin{proof}
  The proof of this is analogous to the proof of \cref{lemma:iso-proj-inj}.
  \qed
\end{proof}

\begin{proof}[\cref{prop:tensor-coalg}]
  By \cref{prop:duals} we have that $C^*$ and $D^*$ are both vN-algebras, and by
  \cref{prop:tensor-alg} we know $C^* \injtens D^*$ is
  a vN-algebra. Applying \cref{prop:duals} gives us that $(C^* \injtens D^*)^*$
  is a vN-coalgebra. Now, by \cref{lemma:iso-proj-inj}
  $(C \projtens D, \varepsilon, j, \delta)$ is a
  comonoid $(C \projtens D, \varepsilon, \delta)$, where $j$ is an involution
  that is reverse comultiplicative. It remains to show that the co-C$^*$-identity
  is satisfied. To do so, suppose that $e: C \projtens D \to \mathbb{C}$ is a
  c.q. map, then $r \eqdef \left( (C^* \injtens D^*)^* \xrightarrow{\theta^{-1}} C \projtens D \xrightarrow{e} \mathbb{C} \right)$ is
  also a c.q. map.
  Now the following diagram commutes:

  \begin{equation*}
    \adjustbox{scale=0.8}{
      \begin{tikzcd}
        (C \projtens D)_o \arrow[d, "\delta_C \otimes \delta_D"'] \arrow[r, "\theta"]                                    & ((C^* \injtens D^*)^*)_o \arrow[dd, "\delta"]                                     \\
        ((C \hagertimes C) \projtens (D \projtens D))_o \arrow[d, "w"']                                                  &                                                                                   \\
        ((C \projtens D) \hagertimes (C \projtens D))_o \arrow[r, "\theta\otimes \theta"] \arrow[d, "\gamma"']      & ((C^* \injtens D^*)^* \hagertimes (C^* \injtens D^*)^*)_o \arrow[d, "\gamma"] \\
        (C \projtens D)_o \hagertimes (C \projtens D)_o \arrow[d, "\mathrm{id} \otimes (j_C \otimes j_D)"'] \arrow[r, "\theta \otimes \theta"] & ((C^* \injtens D^*)^*)_o \hagertimes ((C^* \injtens D^*)^*)_o \arrow[d, "\mathrm{id} \otimes j"]      \\
        (C \projtens D)_o \hagertimes (C \projtens D)_c \arrow[d, "e \otimes e"'] \arrow[r, "\theta \otimes \theta"]     & ((C^* \injtens D^*)^*)_o \hagertimes ((C^* \injtens D^*)^*)_c \arrow[d, "r \otimes r"]           \\
        \mathbb{C} \projtens \mathbb{C} \arrow[equal]{r} \arrow[d, "\cong"']                                                    & \mathbb{C} \projtens \mathbb{C} \arrow[d, "\cong"]                                         \\
        \mathbb{C} \arrow[equal]{r}                                                                                             & \mathbb{C}
        \end{tikzcd}}
  \end{equation*}
  As $(C^* \injtens D^*)^*$ is a vN-coalgebra the right side is c.q, hence the
  left one is too.
  The result for the coproduct follows similarily, using \cref{lemma:iso-one-infty}.
  \qed
\end{proof}

\begin{lemma}
  \label{lemma:star-hom-cp-coalg}
  Let $C, D \in \vNCoalg$ and let $g: C \to D$ be a morphism in $\vNCoalg.$ Then $g$ is a
  CPTP map in the sense of \cref{sec:vNCoalg}.
\end{lemma}

\begin{proof}
  This proof is analogous to the proof of \cref{Lemma:star-hom-cp}.
  \qed
\end{proof}

\subsection{Duality between $\HH$ and $\SSS$}

\begin{lemma}
  \label{lemma:equiv-pos}
  Let $A \in \vNAlg$ and $p : \mathbb{C} \to A$ be a c.b map, then
  the following are equivalent:
  \begin{enumerate}
    \item $p$ is a \emph{positive element} of $A$;
    \item $p^*$ is a \emph{positive linear functional} on $A^*$.
  \end{enumerate}
\end{lemma}

\begin{proof}
  We prove that 1. implies 2. (the other direction is proved completely analogously).
  Let $p : \mathbb{C} \to A$ be a \emph{positive element} of $A$, this means
  that we have some $a \colon \mathbb{C} \to A$ such that the following diagram

  \begin{equation*}
  \begin{tikzcd}
    \mathbb{C}  \arrow[equal]{d} \arrow[r, "\cong"]   & \mathbb{C}_o \hagertimes \mathbb{C}_c \arrow[r, "a \otimes a"] & A_o \hagertimes A_c \arrow[r, "A_o \otimes i"] & A_o \hagertimes A_o  \arrow[r, "\gamma"] & (A \hagertimes A)_o \arrow[d, "\mu"] \\
    \mathbb{C}_o \arrow[rrrr, "p"'] &                                                                    &                                               &                                          & A_o
  \end{tikzcd}
  \end{equation*}
  commutes.
  The following diagram shows that $p^*$ is factorized as required.
  \begin{equation*}
    \begin{tikzcd}
      (A^* \hagertimes A^*)_o \arrow[r, "\gamma"]              & (A^*)_o \hagertimes (A^*)_o \arrow[rr, "(A^*)_o \otimes j"] &  & (A^*)_o \hagertimes (A^*)_c \arrow[rr, "(a^*)_o \otimes (a^*)_c"] & & (\mathbb{C}^*)_o \hagertimes (\mathbb{C}^*)_c \arrow[d, "\cong"] \\
      ((A \hagertimes A)^*)_o \arrow[u, "\cong"]               & (A_o)^* \hagertimes (A_o)^* \arrow[equal]{u} \arrow[rr, "(A_o)^* \otimes i^*"]                      &  & (A_o) \hagertimes (A_c)^* \arrow[u, "\cong"] \arrow[rr, "(a_o)^* \otimes (a_c)^*"] & & (\mathbb{C}_o)^* \hagertimes (\mathbb{C}_c)^* \arrow[d, "\cong"] \\
      ((A \hagertimes A)_o)^* \arrow[r, "\gamma^*"] \arrow[equal]{u} & (A_o \hagertimes A_o)^* \arrow[rr, "(A_o \otimes i)^*"] \arrow[u, "\cong"]      &  & (A_o \hagertimes A_c)^* \arrow[u, "\cong"] \arrow[rr, "(a_o \otimes a_c)^*"]  & & (\mathbb{C}_o \hagertimes \mathbb{C}_c)^* \arrow[d, "\cong"]     \\
      (A^*)_o \arrow[rrrrr, "(p^*)_o"'] \arrow[u, "\mu^*"]          &                                                    &  &                                    & & \mathbb{C}_o
      \end{tikzcd}
  \end{equation*}
  \qed
\end{proof}

\begin{proposition}
  \label{prop:equiv-cp}
  We have the following duality between c.p. maps on vN-algebras and
  c.p. maps on vN-coalgebras:
  \begin{itemize}
    \item Suppose that $A, B \in \vNAlg$ and that $f \colon A \to B$ is a completely positive map.
    Then $f^* : B^* \to A^*$ is a c.p map between vN-coalgebras.
    \item Suppose that $C, D \in \vNCoalg$ and that $g \colon C \to D$ is a completely positive map.
    Then $g^*: D^* \to C^*$ is a c.p map between vN-algebras.
  \end{itemize}
\end{proposition}

\begin{proof}
  We prove the first statement as the second one is proved completely analogously.
  Let $q: T_n \projtens A^* \to \mathbb{C}$ be a positive linear functional.
  By \cref{lemma:equiv-pos}, $q$ is the dual of a positive element
  $p: \mathbb{C} \to (T_n \projtens A^*)^* \cong M_n \injtens A$.
  By assumption, $(M_n \otimes f) \circ p: \mathbb{C} \to M_n \injtens B$ is positive as well,
  i.e. it can be factorised using an element $b: \mathbb{C} \to M_n \injtens B$
  and its involution.
  Let $r$ denote the linear functional $r \defeq \left(T_n \ptimes B^* \cong (M_n \injtens B)^* \xrightarrow{b^*} \mathbb{C} \right)$,
  then the following diagram commutes:
  \begin{equation*}
    \adjustbox{scale=0.55}{
      \begin{tikzcd}
        ((T_n \projtens B^*) \hagertimes (T_n \projtens B^*))_o \arrow[r, "\gamma"]                      & (T_n \projtens B^*)_o \hagertimes (T_n \projtens B^*)_o \arrow[rr, "\mathrm{id} \otimes j_{T_n \projtens B^*}"]                                                                                   &                                                        & (T_n \projtens B^*)_o \hagertimes (T_n \projtens B^*)_c \arrow[rr, "r_o \otimes r_c"]                                                 &  & \mathbb{C}_o \hagertimes \mathbb{C}_c \arrow[d, "\cong"]        \\
        ((M_n \injtens B)^* \hagertimes (M_n \injtens B)^*)_o \arrow[r, "\gamma"] \arrow[u, "\cong"]     & ((M_n \injtens B)^*)_o \hagertimes ((M_n \injtens B)^*)_o \arrow[rr, "\mathrm{id} \otimes j_{(M_n \injtens B)^*}"] \arrow[rrd, "\mathrm{id} \otimes (i_{M_n} \otimes i_B)^*"'] \arrow[u, "\cong"] &                                                        & ((M_n \injtens B)^*)_o \hagertimes ((M_n \injtens B)^*)_c \arrow[rr, "(b^*)_o \otimes (b^*)_c"] \arrow[u, "\cong"] &  & (\mathbb{C}^*)_o \hagertimes (\mathbb{C}^*)_c \arrow[d, "\cong"'] \\
                                                                                                         &                                                                                                                                                                                                   &                                                        & ((M_n \injtens B)_o)^* \hagertimes ((M_n \injtens B)_c)^* \arrow[rr, "(b_o)^* \otimes (b_c)^*"] \arrow[u, "\cong"] &  & (\mathbb{C}_o)^* \hagertimes (\mathbb{C}_c)^* \arrow[d, "\cong"'] \\
        (((M_n \injtens B) \hagertimes (M_n \otimes B))_o)^* \arrow[r, "\gamma^*"] \arrow[uu, "\cong"]     & ((M_n \injtens B)_o \hagertimes (M_n \otimes B)_o)^* \arrow[rr, "(\mathrm{id} \otimes (i_{M_n} \otimes i_B))^*"] \arrow[uu, "\cong"]                                                                  &                                                        & ((M_n \injtens B)_o \hagertimes (M_n \otimes B)_c)^* \arrow[rr, "(b_o \otimes b_c)^*"] \arrow[u, "\cong"]          &  & (\mathbb{C}_o \hagertimes \mathbb{C}_c)^* \arrow[dd, "\cong"]   \\
        (((M_n \hagertimes M_n) \injtens (B \hagertimes B))_o)^* \arrow[u, "v^*"]                        &                                                                                                                                                                                                   &                                                        &                                                                                                                    &  &                                                               \\
        ((M_n \injtens B)_o)^* \arrow[rr, "(M_n \otimes f)^*"'] \arrow[u, "(\mu_{M_n} \otimes \mu_B)^*"] &                                                                                                                                                                                                   & ((M_n \injtens A)_o)^* \arrow[rrr, "p^*"']             &                                                                                                                    &  & (\mathbb{C}_o)^* \arrow[d, "\cong"]                                                             \\
        (T_n \projtens B^*)_o \arrow[u, "\cong"] \arrow[rr, "T_n \otimes f^*"']                          &                                                                                                                                                                                                   & (T_n \projtens A^*)_o \arrow[rrr, "q"'] \arrow[u, "\cong"] &                                                                                                                    &  & \mathbb{C}_o
        \end{tikzcd}}
  \end{equation*}
  Following the diagram above the linear functional
  $q \circ (T_n \otimes f^*): T_n \ptimes B^* \to \mathbb{C}$ is factorized
  by the comultiplication into the linear functional $r$ and its involution.
  \qed
\end{proof}

\begin{proposition}
  The dual on vN-algebras and vN-coalgebras as defined in \cref{prop:duals}
  extends to functors
  \begin{equation*}
    (-)^* \colon \SSS \rightleftarrows \HH^\mathrm{op} \colon (-)^*.
  \end{equation*}
\end{proposition}

\begin{proof}
  We know from \cref{prop:equiv-star-hom} that
  the dual sends unital maps to counital/trace-preserving ones and vice-versa.
  \cref{prop:equiv-cp} now assures us that the functors are well-typed.
  \qed
\end{proof}

\begin{theorem}
  \label{thm:equiv-H-S}
  We have an equivalence of categories $(-)^* \colon \SSS \simeq \HH^\mathrm{op} \colon (-)^*$
  where the action on objects of $(-)^*$ is defined as in \cref{prop:duals} and on morphisms
  in the usual way, i.e. as in $\FdOS.$
\end{theorem}

\begin{proof}
  As any vN-algebra morphism is also CPU (\cref{Lemma:star-hom-cp}) and
  any vN-coalgebra morphism is also CPTP (\cref{lemma:star-hom-cp-coalg}), the
  double dual natural isomorphisms (\cref{lemma:dd-vNalg}) extend to natural isomorphisms:
  \begin{itemize}
    \item $d: \mathrm{Id} \cong (-)^{**}: \HH \to \HH$;
    \item $d: \mathrm{Id} \cong (-)^{**}: \SSS \to \SSS$.
  \end{itemize}
  \qed
\end{proof}

\vNcoalgsmcc*

\begin{proof}
Recall that an equivalence of categories where one of the categories has a symmetric monoidal
structure, induces a symmetric monoidal structure on the other category.
We know that $\vNcoalg \simeq \vNalg^\mathrm{op}$ by \cref{equiv-alg-coalg},
and that $\vNAlg$ is monoidal by \cref{prop:vN-monoidal}.
\cref{lemma:iso-proj-inj} now assures us that the induced tensor is
naturally isomorphic to the one given by $\projtens$ in \cref{prop:tensor-coalg}.
We can make a similar argument for the coproducts.
\qed
\end{proof}

\Ssmcc*

\begin{proof}
  This follows by a similar argument to the proof of \cref{prop:vNcoalgsmcc},
  using the equivalence $\SSS \simeq \HH^\mathrm{op}$ from \cref{thm:equiv-H-S} and the
  monoidal structure on $\HH$ (see \cref{prop:Hsmcc}) instead. The coproducts
  follows by the same argument too.

  It now follows by construction that the inclusion functor preserves the monoidal
  structure and coproducts.
  \qed
\end{proof}

\begin{proposition}
  The equivalences $(-)^* \colon \SSS \simeq \HH^\mathrm{op} \colon (-)^*$
  and $(-)^* \colon \vNcoalg \simeq \vNalg^\mathrm{op} \colon (-)^*$
  are strong monoidal.
\end{proposition}

\begin{proof}
  Follows by \cref{lemma:iso-proj-inj}.
  \qed
\end{proof}

\TnCoalgebra*
\begin{proof}
  For the first part of the statement,
  recall that the canonical vN-algebra structure on $M_n$, has:
  \begin{enumerate}
    \item unit given by $\eta \defeq 1_n: \mathbb{C} \to M_n :: \lambda \mapsto \lambda 1_n$
    \item involution given by the conjugate transpose $i \defeq (-)^* : M_n^c \to M_n^o$
    \item multiplication given by matrix multiplication
    $\mu : M_n \hagertimes M_n \to M_n :: a \otimes b \mapsto ab$
  \end{enumerate}
  To avoid confusion we will use $i$ to denote the conjugate transpose.
  Recall also that we have the isomorphism $T_n \cong (M_n)^*:: a \mapsto \mathrm{tr}(a (-))$.
  By \cref{prop:duals}, $(M_n)^*$ has a vN-coalgebra structure. Taking this structure modulo
  the isomorphism $T_n \cong (M_n)^*$ gives us
  \begin{enumerate}
    \item counit:
    \begin{equation*}
      T_n  \cong  (M_n)^*  \xrightarrow{\eta^*}  \mathbb{C}^*  \cong \mathbb{C} \\
    \end{equation*}
    \begin{equation*}
      \begin{array}{rcl}
        a & \xmapsto{\cong} & \mathrm{tr}(a(-)) \\
          & \xmapsto{\cong \circ \eta^*} & \mathrm{tr}(a1_n) = \mathrm{tr}(a)
      \end{array}
    \end{equation*}
    \item involution:
      \begin{equation*}
        T_n^o  \cong  ((M_n)^*)_o = (M_n^o)^* \xrightarrow{i^*} (M_n^c)^* \cong ((M_n)^*)_c \cong T_n^c
      \end{equation*}
      \begin{equation*}
        \begin{array}{rcl}
          a & \xmapsto{\cong} & \mathrm{tr}(a (-)) \\
            & \xmapsto{(i(-))^*} & \mathrm{tr}(a i(-)) \\
            & \xmapsto{\cong} & \overline{\mathrm{tr}(a i(-))} = \mathrm{tr}(i(a i(-))) = \mathrm{tr}((-) i(a))  = \mathrm{tr}(i(a) (-)) \\
            & \xmapsto{\cong} & i(a)
        \end{array}
      \end{equation*}
    \item comultiplication
    \begin{equation*}
      T_n  \cong  (M_n)^*  \xrightarrow{\mu^*}  (M_n \hagertimes M_n)^*  \cong (M_n)^* \hagertimes (M_n)^*  \cong T_n \hagertimes T_n
    \end{equation*}
      \begin{equation*}
      \begin{array}{rcl}
        e_{ij} & \xmapsto{\cong} & \mathrm{tr}(e_{ij}(-)) \\
               & \xmapsto{\mu^*} & (a \otimes b \mapsto \mathrm{tr}(e_{ij} ab) = \sum_{k} a_{jk}b_{ki}= \sum_{k} \mathrm{tr}(e_{kj}a)\mathrm{tr}(e_{ik}b)) \\
               & \xmapsto{\cong} & \sum_{k} \mathrm{tr}(e_{kj}(-)) \otimes \mathrm{tr}(e_{ik}(-)) \\
               & \xmapsto{\cong} & \sum_{k} (e_{kj} \otimes e_{ik})
      \end{array}
    \end{equation*}
  \end{enumerate}

  For the second part of the statement, let $C \in \vNCoalg$. By
  \cref{prop:duals}, $C^* \in \vNAlg$ and
  we get the following chain of natural isomorphisms in $\vNCoalg$
  \begin{equation*}
    \begin{array}{rclr}
      C & \cong & C^{**} & \qquad \text{(by \cref{lemma:dd-vNalg})} \\
      & \cong & (\linftyplus_{i} M_{n_i})^* & \qquad \text{(by \cref{ex:fd-vN-algebra})} \\
      & \cong & \loneplus_{i} (M_{n_i})^* & \qquad \text{(de Morgan duality)} \\
      & \cong & \loneplus_{i} T_{n_i} & \qquad \text{(by argument above)}
    \end{array}
  \end{equation*}
  \qed
\end{proof}

%%%%%%%%%%%%%%%%%%%%%%%%%%%%%%%%%%%%%%%%%%%%%%%%%%%%%%%%%%%%%%%%%%%%%%%%%%%%%%
\newpage
\section{von Neumann coalgebras and (complete) positivity}
\label{app:vn-coalg-cp}
%%%%%%%%%%%%%%%%%%%%%%%%%%%%%%%%%%%%%%%%%%%%%%%%%%%%%%%%%%%%%%%%%%%%%%%%%%%%%%

The purpose of this appendix is to explain that the definitions of positive
linear functionals and (completely) positive maps between vN-coalgebras from
\secref{sec:vNCoalg} coincide with the usual concrete definitions of these
notions on the spaces $\loneplus_i T_{k_i}$, which give all the f.d.
vN-coalgebras modulo isomorphism. We refer to the vN-coalgebras of the
form $\loneplus_i T_{k_i}$ as \emph{concrete} vN-coalgebras.

Recall that an element $p \in \MM_n$ is called \emph{positive} iff $\langle h, ph \rangle \geq 0$
for every $h \in \mathbb C^n$ where we have written $\langle \cdot , \cdot \rangle \colon \mathbb C^n \times \mathbb C^n \to \mathbb C$
for the standard inner product on $\mathbb C^n.$
Moreover, $p \in \MM_n$ is positive iff there exists $a \in \MM_n$ such that $p = a^* a.$
An element $(p_1, \ldots, p_n) \in \oplus_{1 \leq i \leq n} \MM_{k_i}$ is positive
iff each $p_i \in \MM_{k_i}$ is positive.

\textbf{Notation.} Throughout the remainder of this section, we simply write $V$ and $W$ to refer to two (potentially different) vector spaces of the form $\oplus_{1 \leq i \leq n} \MM_{k_i}$.

Note that the vector spaces $V$ and $W$ can be equipped with the structure of f.d. vN-(co)algebras, as we have already explained in the main body
of the document. Moreover, all f.d. vN-(co)algebras have underlying vector spaces of this form, modulo isomorphism.
We say that a linear map $f \colon V \to W$ is \emph{positive} (in the concrete sense) if $f$ maps positive elements of $V$ to positive elements of $W.$
A \emph{positive linear functional} (in the concrete sense) is a positive map $f \colon V \to \mathbb C$. From the theory of C*-algebras, we have
the following proposition.

\begin{proposition}
  \label{prop:positive-maps}
  Let $f \colon V \to W$ be a linear map. The following are equivalent:
  \begin{itemize}
    \item $f$ is positive;
    \item $f$ reflects positive linear functionals, i.e. for each positive linear functional $p \colon W \to \mathbb C$, the linear functional $p \circ f \colon V \to \mathbb C$ is also positive.
  \end{itemize}
\end{proposition}
\begin{proof}
  $(\Rightarrow)$ This is obvious.

  $(\Leftarrow)$ This follows easily from \cite[Theorem 4.3.4]{kadison-ringrose}. In fact, the statement holds more generally for infinite-dimensional C*-algebras.
\end{proof}

Our first objective is to show that the positivity of a linear functional in the coalgebraic sense is equivalent to positivity in the concrete sense.
We begin with a lemma.

\begin{lemma}
  Let $\varphi \colon T_n \to \mathbb C$ be a linear functional on the indicated vN-coalgebra. The following are equivalent:
  \begin{itemize}
    \item $\varphi$ is positive in the coalgebraic sense, i.e. there exists a linear functional $\psi \colon T_n \to \mathbb C$, such that the following diagram:
      \begin{equation}
      \begin{tikzcd}
        (T_n\hagertimes T_n)^o \arrow[r, "\gamma"]    & T_n^o\hagertimes T_n^o \arrow[r, "T_n^o\otimes j"] & T_n^o \hagertimes T_n^c \arrow[r, "\psi \otimes \psi"] & \mathbb{C}_o \hagertimes \mathbb{C}_c \arrow[r, "\cong"] & \mathbb{C} \arrow[equal]{d} \\
        T^o_n \arrow[u, "\delta"] \arrow[rrrr, "\varphi"'] &                                              &                                                  &                                                         & \mathbb{C}_o
      \end{tikzcd}
      \end{equation}
      commutes.
    \item $\varphi$ is positive in the concrete sense.
  \end{itemize}
\end{lemma}
\begin{proof}
  $(\Rightarrow)$
  Let $\psi \colon T_n \to \mathbb C$ be a linear functional. It follows $\psi = \tr(a\cdot)$ for some $a \in T_n.$
  A diagram chase along the upper leg of the diagram shows that $e_{ij} \in T_n$ is mapped to
  \begin{align*}
    e_{ij} &\xmapsto{\delta} \sum_k e_{kj} \otimes e_{ik} \\
           &\xmapsto{\gamma} \sum_k e_{ik} \otimes e_{kj} \\
           &\xmapsto{T_n^o \otimes j} \sum_k e_{ik} \otimes e_{jk} \\
           &\xmapsto{\psi \otimes \psi} \sum_k \tr(e_{ik}a) \otimes \tr(e_{jk}a) \\
           &\xmapsto{\cong} \sum_k \tr(e_{ik}a) \overline{\tr(e_{jk}a)} \\
           &= \sum_k a_{ki} \overline{a_{kj}}
  \end{align*}
  By taking the linear extension of the above composite and using the commutativity of the above diagram, it follows that
  \[ \varphi(b) = \sum_{i,j,k} a_{ki} \overline{a_{kj}} b_{ij} \]
  for any $b = [b_{ij}] \in T_n.$ Straightforward verification now shows that
  \[ \varphi(b) = \sum_{i,j,k} a_{ki} \overline{a_{kj}} b_{ij}  = \tr(a^* a b) \]
  and therefore $\varphi = \tr((a^* a)\ \cdot)$ is a positive linear functional in the concrete sense.

  $(\Leftarrow)$ Let $\varphi \colon T_n \to \mathbb C$ be a positive linear functional in the concrete sense.
  It follows that $\varphi = \tr(p\ \cdot)$ for some positive element $p \in T_n$ in the concrete sense.
  This is equivalent to $\varphi = \tr((a^* a)\ \cdot)$ for some element $a \in T_n$.
  Now we can define $\psi \colon T_n \to \mathbb C$ by $\psi \eqdef \tr(a\ \cdot)$ and the same reasoning
  as above shows that the required diagram commutes. \qed
\end{proof}

Therefore, for linear functionals $\varphi \colon T_n \to \mathbb C$, the
abstract (coalgebraic) and concrete notions of positivity coincide, so we do
not distinguish them any further.

Next, we extend this equivalence for linear functionals $\varphi \colon \loneplus_{1\leq i \leq n }T_{k_i} \to \mathbb C$.
Every such functional is necessarily of the form $\varphi = [\varphi_1, \ldots, \varphi_n] $ where we have $\varphi_i \colon T_{k_i} \to \mathbb C$
and $\varphi(t_1, \ldots, t_n) = \sum_i \varphi_i(t_i).$

\begin{lemma}
  Let $\varphi \colon \loneplus_{1\leq i \leq n }T_{k_i} \to \mathbb C$ be a linear functional.
  Then $\varphi$ is positive in the coalgebraic sense iff each $\varphi_i \colon T_{k_i} \to \mathbb C$
  (from the above decomposition) is positive.
\end{lemma}
\begin{proof}
  For simplicity, let us assume that we have only two entries in the direct sum, i.e. we have a linear functional
  $\varphi \colon T_n \loneplus T_m \to \mathbb C.$ The general case is analogous.

  $(\Rightarrow)$ Let $\varphi$ be positive in the coalgebraic sense. By performing a diagram chase in \eqref{eq:positive-functional} on
  elements of the form $(t,0) \in T_n \oplus T_m$ we see that $\varphi_1$ is positive. The argument for $\varphi_2$ is symmetric.

  $(\Leftarrow)$ If both $\varphi_1$ and $\varphi_2$ are positive, then there exist $t_1 \in T_n$ and $t_2 \in T_m$ such that
  $\varphi_1 = \tr(t_1^* t_1\ \cdot)$ and $\varphi_2 = \tr(t_2^* t_2\ \cdot).$ We can now define a linear functional
  $a \colon T_n \loneplus T_m \to \mathbb C$ by $a(x_1,x_2) =\tr(t_1 x_1) + \tr(t_2 x_2). $
  To verify that diagram \eqref{eq:positive-functional} commutes, it suffices to perform the diagram chase on elements of the
  form $(e_{ij}, 0)$ and $(0, e_{lp})$. In both of these cases, the diagram chase easily reduces to the diagram chase in the
  proof of the preceding lemma and the commutativity holds using the same arguments.
\end{proof}

\begin{proposition}
  \label{prop:positive-functionals}
  Let $\varphi \colon \loneplus_{1\leq i \leq n }T_{k_i} \to \mathbb C$ be a linear functional. Then,
  $\varphi$ is positive in the coalgebraic sense iff $\varphi$ is positive in the concrete sense.
\end{proposition}
\begin{proof}
  By the previous lemma, it suffices to show that $\varphi$ is positive in the
  concrete sense iff each $\varphi_i \colon T_{k_i} \to \mathbb C$ is positive.
  This is obviously true.
  \qed
\end{proof}

Therefore, for concrete vN-coalgebras, the positivity of linear functionals in
the abstract (coalgebraic) and concrete senses coincide.

\begin{theorem}
  Let $f \colon C_1 \to C_2$ be a linear map between two concrete f.d. vN-coalgebras. Then:
  \begin{enumerate}
    \item $f$ is positive in the abstract (coalgebraic) sense iff $f$ is positive in the concrete sense;
    \item $f$ is completely positive in the abstract (coalgebraic) sense iff $f$ is completely positive in the concrete sense.
  \end{enumerate}
\end{theorem}
\begin{proof}
  The first statement follows by combining Propositions \ref{prop:positive-functionals} and \ref{prop:positive-maps}. The second
  statement now immediately follows from the first one, because the linear map
  $T_n \ptimes \varphi \colon T_n \ptimes C_1 \to T_n \ptimes C_2$, from the coalgebraic definition of complete positivity, is precisely the same linear map
  as $\MM_n \otimes \varphi \colon \MM_n \otimes C_1 \to \MM_n \otimes C_2$ in the concrete definition of complete positivity.
  \qed
\end{proof}

It now easily follows that an element $p \in C$ of a concrete f.d. vN-coalgebra
is positive in the coalgebraic sense (see \secref{sec:qcs}) iff $p$ is
positive in the concrete sense.

To summarise, all the notions of (complete) positivity that we have introduced
in the coalgebraic sense coincide with the standard concrete definitions of
(complete) positivity when we consider concrete f.d. vN-coalgebras.

Our next objective is to show that the set of density
operators on a vN-coalgebra $C$ may be characterised in two equivalent ways:
\begin{equation*}
  \{ c \in C \ |\ \norm{c} = 1 = \varepsilon(c) \} = \{ \rho \in C \ |\ \varepsilon(\rho) = 1 \text{ and } \rho \geq 0 \} = P_C .
\end{equation*}

Towards this end, we prove the following proposition. Its dual analogue for
vN-algebras/C*-algebras is well-known, so we use a simple duality argument to
establish our next proposition.

\begin{proposition}
  \label{prop:density-equivalent}
  Let $C$ be a f.d. vN-coalgebra and let $\rho \in C$ be an element with $\varepsilon(\rho) = 1.$ The
  following are equivalent:
  \begin{enumerate}
    \item $\rho \geq 0$, i.e. $\rho$ is positive;
    \item $\norm{\rho} = 1$;
    \item $\norm{\rho} \leq 1$.
  \end{enumerate}
\end{proposition}
\begin{proof}
  Consider the linear map $p \colon \mathbb C \to C$ defined by $p(1) \eqdef \rho.$ Since $\varepsilon(\rho) = 1$, the map
  $p$ is counital (i.e. trace-preserving).

  (1) $\Longrightarrow$ (2).
  The map $p$ is positive (by definition). The proof of Proposition \ref{prop:equiv-cp} shows that
  $C^* \xrightarrow{p^*} \mathbb C^* \cong \mathbb C$ is positive (simple special case)
  and it is clearly unital.
  From the theory of C*-algebras, we know that every positive unital functional $\varphi$ has norm $\norm{\varphi} = 1 = \norm{\varphi}_{\mathrm{cb}}$
  \cite[Theorem 4.3.2]{kadison-ringrose}
  (and is also completely positive).
  The o.s. dual $(-)^*$ preserves the cb-norm of $p$, so $ 1 = \norm{p^*}_{\mathrm{cb}} = \norm{p}_{\mathrm{cb}} = \norm{p} = \norm{\rho}.$

  \noindent (2) $\Longrightarrow$ (3). Obvious.

  \noindent (3) $\Longrightarrow$ (1).
  The map $p$ is a (necessarily complete) contraction and also trace-preserving.
  By Proposition \ref{prop:coalg-cc-cp} it follows that $p$ is completely positive and
  therefore also positive, so $p(1) = \rho$ is positive (by definition).
  \qed
\end{proof}

\begin{corollary}
  \label{cor:density}
  Let $C$ be a f.d. vN-coalgebra. Then the set of density operators $P_C$ may be defined in two equivalent ways:
\begin{equation}
  \label{eq:density}
  P_C = \{ \rho \in C \ |\ \varepsilon(\rho) = 1 \text{ and } \rho \geq 0 \} = \{ \rho \in C \ |\ \norm{\rho} = 1 = \varepsilon(\rho) \}
\end{equation}
\end{corollary}

%%%%%%%%%%%%%%%%%%%%%%%%%%%%%%%%%%%%%%%%%%%%%%%%%%%%%%%%%%%%%%%%%%%%%%%%%%%%%%
\newpage

%%%%%%%%%%%%%%%%%%%%%%%%%%%%%%%%%%%%%%%%%%%%%%%%%%%%%%%%%%%%%%%%%%%%%%%%%%%%%%
\section{Background on double gluing and ommitted proofs from \cref{sec:qcs}}
In this section we revisit the hom-functor double gluing construction for classical
linear logic from \cite{HS}.

\begin{definition}[{Glued category, \cite[Sec. 4.3.2]{HS}}]
  Given a $*$-autonomous category $(\mathbf{C}, I, \otimes, [ - , - ])$ the (hom-functor)
  double glued category $G(\mathbf{C})$ has
  \begin{itemize}
    \item objects: $(c, S, V)$,
  with $c$ an object in $\mathbf{C}$ and $S \subseteq \mathbf{C}(I,c)$ and
  $R \subseteq \mathbf{C}(c, \bot)$
    \item morphisms:
    $f : (c, S, V) \to (d, R, W)$,
    where $f : c \to d$ in $\mathbf{C}$ and
    for all $s \in S$ we have $f \circ s \in R$ and for all $w \in W$ we have
    $ w \circ f \in V$
  \end{itemize}
\end{definition}
\begin{remark}
  When we do double gluing along the hom-functor, we can see the second
argument in an object as a set
of ``proofs'' of the formula $A$ and the third argument as a set
of ``proofs'' of the formula $\lnot A$ or ``counter-proofs''.
\end{remark}
\begin{proposition}[{Glued category, \cite[Prop. 41]{HS}}]
  The category $G(\mathbf{C})$ is $*$-autonomous, with
  \begin{itemize}
    \item Negation: $(c, S, V)^* \eqdef (c^*, V, S)$
    \item Tensor: $(c, S, V) \otimes (d, R, W) \eqdef (c \otimes d, S\otimes R, G(\mathbf{C})((c, S, V), (d, R, W)^*))$
    where $S\otimes R \eqdef \{s\otimes u \ | \ s\in S, \ u \in U\}$,
    with unit: $(I, \{\text{id}_I\}, \mathbf{C}(I,\bot))$
  \end{itemize}
\end{proposition}
This $*$-autonomous structure is not very meaningful, as we have no interaction
between the ``proofs'' and ``counter-proofs'', hence the dual is trivial.
Note also that the monoidal structure adds a lot of extra data to the third
argument. We thus need to carve out a subcategory of ``meaningful'' objects.
To impose a restriction on objects we want the ''proofs" to interact with the
''counter-proofs", we thus use the notion of \emph{focussed orthogonality}
from \cite{HS}.

\begin{definition}[focussed orthogonality]
  Let $\mathbf{C}$ be a $*$-autonomous category, then a \emph{focussed orthogonality}
  $\independent$ on $\mathbf{C}$ is a family of relations
  $\{\independent_c \in \mathbf{C}(I,c)\times \mathbf{C}(c, \bot) \}$ induced by
  $\independent \subseteq C(I, \bot)$ meaning that
  \begin{equation*}
    u \independent_c f \quad \iff \quad f \circ u \in \independent
  \end{equation*}
\end{definition}

We say that a focussed orthogonality induced by
$\independent \subseteq \mathbf{C}(I, \bot)$ has \emph{focus} on $\independent$.

\begin{proposition}[{focussed precise, \cite[Sec. 5.3]{HS}}]
  \label{precise}
  Let $\independent$ be a focussed orthogonality on $\mathbf{C}$, then $\independent$
  is \emph{precise}, meaning it satisfies
  \begin{itemize}
    \item For all $u: I \to c$, $v: I \to d$ and $f:c\otimes d \to \bot$
      \begin{equation*}
        u \independent \langle v | f \rangle \quad \text{and} \quad v \independent \langle u | f \rangle
        \quad \iff \quad u \otimes v \independent f
      \end{equation*}
      where
      \begin{equation*}
        \langle v | f \rangle \defeq \left( c \cong c \otimes I \xrightarrow{v} c \otimes d \xrightarrow{f} \bot \right)
      \end{equation*}
      \begin{equation*}
        \langle u | f \rangle \defeq \left( d \cong I \otimes d \xrightarrow{u} c \otimes d \xrightarrow{f} \bot \right)
      \end{equation*}
    \item For all $u: I \to c$, $v: d \to \bot$ and $f: c \to d$
      \begin{equation*}
        f \circ u \independent v \quad \text{and} \quad u \independent v \circ f \quad \iff
        \quad f \independent u \multimap v
      \end{equation*}
  \end{itemize}
\end{proposition}
We introduce some notation, let $S \subseteq \mathbf{C}(I, c)$, $R \subseteq \mathbf{C}(I, d)$ and
$V \subseteq \mathbf{C}(c, \bot)$, then
\begin{equation*}
  S \otimes R \eqdef \{s \otimes r \ | \ s \in S \text{ and } r \in R \} \subseteq
  \mathbf{C}(I, c \otimes d)
\end{equation*}
\begin{equation*}
  R \multimap V \eqdef \{r \multimap v \ | \ r \in R \text{ and } v \in V \} \subseteq
  \mathbf{C}(d \multimap c , \bot)
\end{equation*}
The \emph{polar set} of a subset $S \subseteq \mathbf{C}(I,c)$ can now be defined as
$S^\circ \eqdef \{f \in C(c, \bot) \ | \ \forall s\in S \ s \independent f\}$.
The set $S^\circ$ can be seen as the maximal subset where all elements are
orthogonal to $S$. It has the following properties:
\begin{lemma}[{\cite[Sec. 5.1]{HS}}]
  \label{galois}
  For an $R \subseteq S\subseteq \mathbf{C}(I, c)$
  we have that $S\subseteq S^{\circ\circ}$, $S^{\circ\circ\circ} = S^{\circ}$
  and $S^\circ \subseteq R^\circ$
\end{lemma}
\begin{definition}[{focussed orthogonal category, \cite[Sec. 5.1]{HS}}]
  Let $\mathbf{C}$ be a symmetric monoidal closed category and
  $\independent \subseteq \mathbf{C}(I,c)$, the \emph{focussed orthogonal category}
  denoted $\mathcal{O}_{\footnotesize{\independent}}(\mathbf{C})$ has
  \begin{itemize}
    \item objects: $(c, S, V)$, where $c\in \text{ob}(\mathbf{C})$, $S\subseteq \mathbf{C}(I, c)$,
    $V\subseteq \mathbf{C}(c,\bot)$ such that $S=V^\circ$ and $V=S^\circ$
    \item morphisms: $f:(c, S, V) \to (d, R, W)$, where
    $f:c \to d \in \text{mor}(\mathbf{C})$, such that $f\circ s \in R$ and $w \circ f \in V$
  \end{itemize}
\end{definition}
The condition that $S=V^\circ$ and $V=S^\circ$ is equivalent to
$S=S^{\circ\circ}$, we can thus see the objects of
$\mathcal{O}_{\independent}(\mathbf{C})$ as pairs $(c, S^{\circ\circ})$ instead of triples.
\begin{proposition}[{stable}]
  \label{stable}
  A focussed orthogonality $\independent$ on $\mathbf{C}$ is
  \emph{stable}, that is, for all $S\subseteq \mathbf{C}(I,c)$, $R\subseteq \mathbf{C}(I, d)$ and
  $W\subseteq \mathbf{C}(d,\bot)$
  \begin{itemize}
    \item $(S^{\circ\circ}\otimes R^{\circ\circ})^{\circ}=
    (S\otimes R^{\circ\circ})^{\circ}=
    (S^{\circ\circ}\otimes R)^{\circ}$
    \item $(R^{\circ\circ} \multimap W^{\circ\circ})^\circ =
    (R \multimap W^{\circ\circ})^\circ =
    (R^{\circ\circ} \multimap W)^\circ$
  \end{itemize}
\end{proposition}

\begin{proof}
  \begin{itemize}
    \item For a focussed orthogonality we have the following equivalence
    \begin{equation}
      s \independent \langle r | f \rangle \quad \iff \quad r \independent \langle s | f \rangle
    \end{equation}
    as the two composites below are equal
    \begin{equation*}
      I \xrightarrow{s} S \otimes I \xrightarrow{S \otimes r} S \otimes R \xrightarrow{f} \bot
    \end{equation*}
    \begin{equation*}
      I \xrightarrow{r} R \otimes I \xrightarrow{s \otimes R} S \otimes R \xrightarrow{f} \bot
    \end{equation*}

    First equality: By \cref{galois} it is sufficient to prove that $(S \otimes R^{\circ\circ})^{\circ}
    \subseteq (S^{\circ\circ} \otimes R^{\circ\circ})^{\circ}$. Suppose
    $f\in (S \otimes R^{\circ\circ})^{\circ}$, by preciseness
    (\cref{precise}) we have that $\langle r | f \rangle \in S^\circ$ for all
    $r \in R^{\circ\circ}$. Consequently we have
    $s \independent \langle r | f \rangle$ for all $s\in S^{\circ\circ}$ and
    $r \in R^{\circ\circ}$ , hence by the equivalence above
    $f \in (S^{\circ\circ} \otimes R^{\circ\circ})^{\circ}$.
    The second equality can be proved in a similar way.

    \item For a focussed orthogonality
    \begin{equation*}
      f \circ r \independent w \quad \iff \quad r \independent w \circ f
    \end{equation*}
    as $w \circ (f \circ u) = (w \circ f) \circ u \in \independent$.

    First equality: By \cref{galois} it is sufficient to prove
    $(R \multimap W^{\circ\circ})^\circ \subseteq
    (R^{\circ\circ} \multimap W^{\circ\circ})^\circ$. Suppose
    $f \in (R \multimap W^{\circ\circ})^\circ$, then by preciseness
    $w \circ f \in R^\circ$ for all $w\in W^{\circ\circ}$. Consequently
    we have $r \independent w \circ f$ for all $r\in R^{\circ\circ}$ and
    $w\in W^{\circ\circ}$, hence by preciseness
    $f \in (R^{\circ\circ} \multimap W^{\circ\circ})^\circ$.
    The second equality can be proved in a similar way.
  \end{itemize}
  \qed
\end{proof}

\begin{theorem}[{\cite[Prop. 62]{HS}}]
  \label{thrm:orth-star-autonomous}
  Let $\mathbf{C}$ be a $*$-autonomous category and $\independent \subseteq \mathbf{C}(I,\bot)$,
  then the induced focussed orthogonal category $\mathcal{O}_{\independent}(\mathbf{C})$
  is $*$-autonomous, with
  \begin{itemize}
    \item Tensor: $(c,S) \projtens (d, R) \eqdef (A\projtens B, (S\otimes U)^{\circ\circ})$
    with unit: $(I, \{\text{id}_I\}^{\circ\circ})$, where
    $S\otimes U \eqdef \{s\otimes u \ | \ s\in S, \ u \in U\}$
    \item Negation: $(c, S)^* \eqdef (c^*, S^\circ)$
  \end{itemize}
\end{theorem}

\begin{proposition}[{focussed morphism, \cite[Sec. 5.2]{HS}}]
  \label{focussed}
  Let $\mathbf{C}$ be a $*$-autonomous category with an orthogonality
  $\independent \subseteq \mathbf{C}(I,\bot)$ and $S \subseteq \mathbf{C}(I, c)$, $V \subseteq \mathbf{C}(d, \bot)$,
  then any $f: c \to d$ is \emph{focussed (w.r.t $S$ and $V$)}, meaning that
  \begin{equation*}
    f \circ s \independent v \quad \iff \quad s \independent v \circ f \quad
    \forall s \in S \ \forall v\in V
  \end{equation*}
\end{proposition}

\begin{proposition}[{(co)products}]
  \label{coproduct-orth}
  Let $\mathbf{C}$ be a $*$-autonomous category with an orthogonality
  $\independent \subseteq \mathbf{C}(I,\bot)$
  \begin{itemize}
    \item Suppose $\mathbf{C}$ is cartesian, then $\mathcal{O}_{\independent}(\mathbf{C})$ is
    cartesian. The binary products are explicitly defined as
    \begin{equation*}
      (c, S) \times (d, R) \eqdef (c\times d, S\times R)
    \end{equation*}
    where $S \times R \eqdef \{\langle s, r \rangle: I \to c \times d \ | \ s\in S, r \in R \}$
    and the terminal object is $(1, \mathbf{C}(I, 1))$.
    \item Suppose $\mathbf{C}$ is cocartesian, then $\mathcal{O}_{\independent}(\mathbf{C})$ is
    cocartesian. The binary coproducts are explicitly defined as
    \begin{equation*}
      (c, S) + (d, R) \eqdef (c + d, (S^\circ + R^\circ)^\circ)
    \end{equation*}
    where $S^\circ + R^\circ \eqdef \{ [ f, g ] \colon c + d \to \bot \ | \ f \in S^\circ, g \in R^\circ \} $
    and the initial object is $(0, \mathbf{C}(0, \bot)^\circ)$.
  \end{itemize}
\end{proposition}

\begin{proof}
  By \cite[Prop. 63]{HS} it is sufficient to show that
  projections (injections) are focussed. We know by \cref{focussed} that all
  morphisms in a focussed orthogonal category are focussed, hence projections
  (injections) are too.
  \qed
\end{proof}

Now we have that $\QCoh$ is the \emph{focussed orthogonal category} of
$G(\fdOS)$ with focus
$\{1\}\subseteq \fdOS(\mathbb{C}, \mathbb{C})=\Ball (\mathbb{C})$.

\qMall*
\begin{proof}
    Any focussed orthogonality is \emph{precise} by \cref{precise}
    and \emph{stable} by \cref{stable}. Since $\fdOS$ is $*$-autonomous, it follows that $\QCoh$
    is $*$-autonomous as well with the given structure by
    \cref{thrm:orth-star-autonomous}.
    We also know that $\fdOS$ has finite products, thus $\QCoh$ has
    finite products (and consequently also coproducts) as given by \cref{coproduct-orth}.
  \qed
\end{proof}
\newpage
%%%%%%%%%%%%%%%%%%%%%%%%%%%%%%%%%%%%%%%%%%%%%%%%%%%%%%%%%%%%%%%%%%%%%%%%%%%%%%

\section{Bipolar subsets in $\QQ$}
\label{app:bipolar}
%%%%%%%%%%%%%%%%%%%%%%%%%%%%%%%%%%%%%%%%%%%%%%%%%%%%%%%%%%%%%%%%%%%%%%%%%%%%%%

Our notion of polarity is a strong condition which implies that non-trivial
polar sets must be subsets of the unit \emph{sphere}.

\begin{proposition}
  Given an operator space $X$ and a non-empty subset $S \subseteq \Ball(X)$, if
  $S^\circ\neq \varnothing,$ then $\forall s \in S.\ \norm{s} = 1$ and
  $\forall s' \in S^\circ.\ \norm{s'}=1.$
\end{proposition}
\begin{proof}
  Let $s \in S$ and $s' \in S^\circ$. By definition, $s'(s) = 1$ and therefore
  $ 1 = \norm{s'(s)} \leq \norm{s'} \norm{s} \leq 1$,
  because $s$ and $s'$ are in the unit balls of $X$ and $X^*.$
  This obviously implies that $\norm{s} = \norm{s'} = 1.$
\end{proof}

We recall some facts from functional analysis. For a matrix $a \in \MM_n$, we
write $\norm{a}$ for its operator norm, i.e. its norm as an element $a \in M_n$.
We write $\trnorm{a}$ for its trace norm, which is simply its norm as an element $a \in T_n.$

\begin{lemma}
  \label{lem:trace-operator-norm}
  Let $a, b \in \MM_n$ be two matrices and $u \in \MM_n$ a unitary matrix. Then:
  \begin{enumerate}
    \item $\norm{u} = 1.$
    \item $\trnorm{u} = n.$
    \item $\abs{\trace(ab)} \leq \trnorm{a} \norm{b}$.
  \end{enumerate}
\end{lemma}
\begin{proof}
  The first two statements follow easily from the definitions. The last statement is a special case of \cite[p. 267]{conway-course}.
\end{proof}

\begin{lemma}
  \label{lem:strict-convexity}
  Let $X$ be a strictly convex normed space and let $x_1, \ldots, x_n \in X$ with $\norm{x_i} = 1$ be $n$ points on the unit sphere.
  Assume further that $\norm{\sum_{k=1}^n x_k} = n$. Then $x_1 = x_2 = \cdots = x_n.$
\end{lemma}
\begin{proof}
  A strictly convex normed space is characterised by the property \cite[p. 263]{strict-convexity}:
  \begin{center}
    if $\norm{x + y} = \norm{x} + \norm{y}$ and $x \neq 0$, then $y = cx$ for some $c \geq 0.$
  \end{center}
  The lemma clearly holds for $n=0$ or $n=1$. The rest follows via a simple inductive argument.
  If the lemma is true for $n=k \geq 1$, then to prove it for $n=k+1$, we reason as follows.
  Let $y = \sum_{i=1}^{k} x_i$. Then
  \[ k+1 = \norm{y + x_{k+1}} \leq \norm{y} + \norm{x_{k+1}} \leq  k + 1. \]
  Therefore $\norm{y + x_{k+1}} = \norm{y} + \norm{x_{k+1}}$ and also $\norm{y} = k$.
  From the latter and the induction hypothesis, it follows $x_1 = \cdots = x_k.$
  From the former and strict convexity, it follows $x_{k+1} = cy$ for some $c \geq 0.$
  It follows $c=\frac{\norm{x_{k+1}}}{\norm{y}} = \frac{1}{k}$ and therefore
  $x_1 = \cdots = x_{k+1}.$ \qed
\end{proof}

\begin{lemma}
  \label{lem:iso-lift}
  Let $i \colon X \to Y$ be an isomorphism in $\fdOS$. If $(X,U)$ is an object of $\QQ,$ then
  $(Y, i[U])$ is also an object of $\QQ,$ where $i[U] \eqdef \{i(u) \ |\ u \in U \}$
  is the image of $U$ under $i.$
\end{lemma}
\begin{proof}
  This can be proven by straightforward verification.
\end{proof}

%%%%%%%%%%%%%%%%%%%%%%%%%%%%%%%%%%%%%%
\subsection{Bipolar sets relevant for von Neumann algebras}
%%%%%%%%%%%%%%%%%%%%%%%%%%%%%%%%%%%%%%

We now proceed to identify important bipolar sets that we can associate to f.d. vN-algebras in $\QQ.$

\begin{lemma}
  \label{lem:unitary-polar}
  Let $X$ be an operator space whose underlying vector space is $\MM_n$ and such that
  for every $a \in \MM_n$ we have $\norm{a}_{M_n} \leq \norm{a}_X$. Assume that $(X,S) \in \Ob(\QQ),$
  i.e. $S$ is bipolar.
  If $S^\circ \subseteq \Ball(X^*)$ contains an element of the form $\frac{1}{n} \tr(u^*\ \cdot),$
  where $u \in \MM_n$ is a unitary matrix, then $S \subseteq \{u\}$.
\end{lemma}
\begin{proof}
  Let $g \in S = S^\bipolar$ be an arbitrary element.
  By definition of $S^\bipolar$ we have that $\tr( \frac{1}{n} u^*g) = 1$
  and therefore  $\tr( u^*g) = n.$ Let $\{e_k\}_{k=1}^n$ be an orthonormal basis of $\mathbb C^n.$ Then
  \begin{align*}
    n
    &= \tr( u^*g) & \\
    &= \sum_{k=1}^n \langle u^*g e_k, e_k \rangle & (\text{definition of } \trace{}) \\
    &= \sum_{k=1}^n \langle g e_k, u e_k \rangle & (\text{adjoint})
  \end{align*}
  and therefore  $n = \abs{\sum_{k=1}^n \langle g e_k, u e_k \rangle}$. The Cauchy-Schwarz inequality gives for every $1 \leq k \leq n:$
  \begin{equation}
    \label{eq:cs-argument}
    \abs{\langle g e_k, u e_k \rangle} \leq \norm{g e_k} \norm{u e_k} = \norm{g e_k} \leq \norm{g}_{M_n} \norm{e_k} \leq \norm{g}_X \leq 1
  \end{equation}
  Writing $c_k \eqdef \langle g e_k, u e_k \rangle \in \mathbb C$, we therefore have
  that $0 \leq \abs{c_k} \leq 1$ and
  \[ n = \abs{\sum_{k=1}^n c_k} \leq \sum_{k=1}^n \abs{c_k} \leq n \]
  and therefore $\abs{c_k} = 1$ for every $1 \leq k \leq n.$ It follows that $c_1 = c_2 = \cdots = c_n = 1$ by Lemma \ref{lem:strict-convexity}.
  Returning to \eqref{eq:cs-argument}, we have that
  \[ 1 = \abs{\langle g e_k, u e_k \rangle} =  \norm{g e_k} \norm{u e_k} \]
  and by applying the Cauchy-Schwarz argument again (for the equality) we get that
  \[ g e_k = \frac{\langle g e_k, u e_k \rangle}{\norm{u e_k}^2} u e_k = u e_k \]
  for every $1 \leq k \leq n.$ Therefore $g = u.$
  \qed
\end{proof}

\begin{lemma}
  \label{lem:unitary-Mn}
  Let $u \in M_n$ be a unitary element with $n \in \mathbb N.$ Then
  the pair $(M_n, \{u\})$ is an object of $\QQ$.
\end{lemma}
\begin{proof}
  If $n = 0$, then $M_0 \cong 0$ and the unique (zero) element $0 \in M_0$ is unitary, because it is also the identity element.
  The pair $(M_0, \{ 0 \})$ is the terminal object of $\QQ$ by \cref{thm:q-mall}.

  If $n >0$, then let $ U \eqdef \{ u \}.$
  Consider the pair $ (M_n^*, U^\polar ).$
  By definition, we have that
  \[ U^\polar = \{ f \colon M_n \to \mathbb C \ |\ f(u) = 1 \} \subseteq \fdOS(M_n, \mathbb C) . \]
  It follows that $\trace(\frac{1}{n}u^* \cdot) \in U^\polar$, where $u^*$ is the conjugate transpose (and inverse) of $u.$
  Indeed, we have that $\trace(\frac{1}{n}u^* u) = \frac{1}{n} \trace(1_{n}) = 1$ and to see that $\trace(\frac{1}{n}u^* \cdot)$
  is a (necessarily complete) contraction, we verify that if $b \in M_n$, then
  \begin{align*}
       \abs{\trace \left (\frac{1}{n}u^* b \right)}
    &= \frac{1}{n} \abs{\trace \left (u^* b \right)} & \\
    &\leq \frac{1}{n} \trnorm{u^*} \norm{b} & (\text{Lemma \ref{lem:trace-operator-norm}})\\
    &\leq \norm{b}  & (\text{Lemma \ref{lem:trace-operator-norm}})
  \end{align*}
  We always have that $U \subseteq U^\bipolar$ and \cref{lem:unitary-polar} gives the other direction.
  \qed
\end{proof}

\begin{restatable}{proposition}{vnAlgUnitary}
  \label{prop:vn-algebra-unitary-object}
  Let $M$ be a f.d. vN-algebra and $u \in M$ a unitary element, i.e. $uu^* =
  1_M = u^*u.$ Then $(M, \{u\}) \in \Ob(\QQ).$
\end{restatable}
\begin{proof}
  From \cref{ex:fd-vN-algebra},
  we have a vN-algebra isomorphism $M \cong \bigoplus_{1 \leq i \leq n}^{\infty}M_{k_i}$ for some $n \in \mathbb N$ and $k_i \geq 1$. This isomorphism preserves unitaries and therefore
  by Lemma \ref{lem:iso-lift}, it suffices to show that
  \[ \left( \bigoplus_{1 \leq i \leq n}^{\infty}M_{k_i}, \{u\} \right) \]
  is an object of $\QQ$ for any unitary $u \in \bigoplus_{1 \leq i \leq n}^{\infty}M_{k_i}$.
  The unitary $u$ is of the form
  \[ u = (u_1, \ldots , u_n) \]
  with $u_i \in M_{k_i}$ and such that each $u_i$ is a unitary in $M_{k_i}$. By Lemma \ref{lem:unitary-Mn} each pair $(M_{k_i}, \{u_i\})$ is an object of $\QQ$
  and by \cref{thm:q-mall} we have that
  \[ \left( \bigoplus_{1 \leq i \leq n}^{\infty}M_{k_i}, \{(u_1, \ldots , u_n )\} \right)  \in \Ob(\QQ) \]
  is the categorical product of the objects $(M_{k_i}, \{u_i\})$.
  \qed
\end{proof}

It is now easy to prove that we have a fully faithful inclusion $H \colon \HH \fullsubright \QQ$.

\QvNAlg*
\begin{proof}
  By the previous proposition, we see that the functor $H$ is well-defined on objects.
  Using Proposition \ref{prop:alg-cc-cp}, it follows that $H$ is well-defined on morphisms.
  Faithfulness is obvious and fullness follows by Proposition \ref{prop:alg-cc-cp} (again)
  because any morphism $f \colon (A, \{1_A\}) \to (B, \{1_B \})$
  in $\QQ$ is clearly a unital complete contraction and therefore completely positive.

  Theorem \ref{thm:q-mall} shows immediately that $H$ strictly preserves products.
  It remains to show that $H$ strictly preserves $\itimes$. We show this first for
  vN-algebras of the form $M_n$. This amounts to proving that
  $S \eqdef ( \{ 1_n \}^\circ \otimes \{ 1_m \}^\circ )^\circ = \{ 1_n \otimes 1_m \}$ when viewed as subsets of $M_n \itimes M_m$.
  To see that $ 1_n \otimes 1_m  \in S$, let
  \[ \varphi \otimes \psi \in (\{ 1_n \}^\circ \otimes \{ 1_m \}^\circ) . \]
  Then,
  by definition $\varphi(1_n) = 1 = \psi(1_m)$ and so $(\varphi \otimes \psi)(1_n \otimes 1_m) = 1,$
  as required. To see that $S \subseteq \{ 1_n \otimes 1_m \},$ observe that we have linear functionals
  $\frac{1}{n}\tr \in \{ 1_n \}^\circ$ and $\frac{1}{m}\tr \in \{ 1_m \}^\circ$, which we know from the proof
  of Lemma \ref{lem:unitary-Mn}. The functional $\frac{1}{n}\tr \otimes \frac{1}{m}\tr$ may now be identified
  with $\frac{1}{nm}\tr$ via the vN-algebra isomorphism $M_n \itimes M_m \cong M_{nm}$ and the proof of Lemma \ref{lem:unitary-Mn}
  shows that the only element $g \in \Ball(M_{nm})$ that satisfies $\tr(g) = nm$ is $g = 1_{nm}$ which corresponds to
  $1_n \otimes 1_m$ via the isomorphism $M_{nm} \cong M_n \itimes M_m$, as required.

  Therefore $H(M_n) \itimes H(M_m) = H(M_n \itimes M_m).$ If $A = M_{n} \linftyplus M_k$, then
  \begin{align*}
    H(A \itimes M_m) &\xrightarrow{H(\cong)} H( (M_{n} \itimes M_m )\linftyplus (M_k \itimes M_m)) \\
                     &= H( M_{n} \itimes M_m ) \linftyplus H(M_k \itimes M_m) \\
                     &= (H( M_{n}) \itimes H( M_m ) ) \linftyplus (H(M_k) \itimes H(M_m)) \\
                     &\xrightarrow{\cong^{-1}} (H(M_n) \linftyplus H(M_k)) \itimes H(M_m) \\
                     &= H(A) \itimes H(M_m)
  \end{align*}
  and since the action of $H$ on morphisms is trivial and the MALL structure of $\QQ$ coincides with that of $\FdOS$ (under the forgetful functor),
  it follows that the above
  composite is the identity map. This argument can be obviously generalised and symmetrised and it follows $H(A) \itimes H(B) = H(A \itimes B)$
  for all concrete vN-algebras. The result for all (abstract) f.d. vN-algebras follows by using Example \ref{ex:fd-vN-algebra} and similar arguments.
  \qed
\end{proof}

%%%%%%%%%%%%%%%%%%%%%%%%%%%%%%%%%%%%%%
\subsection{Bipolar sets relevant for von Neumann coalgebras}
%%%%%%%%%%%%%%%%%%%%%%%%%%%%%%%%%%%%%%

We now turn our attention to the inclusion of f.d. vN-coalgebras in $\QQ.$

\begin{proposition}
  \label{prop:vn-coalgebra-trace-object}
  Let $C$ be a finite-dimensional von Neumann coalgebra with counit $\varepsilon \colon C \to \mathbb C$ and let
  $S_C \eqdef \{ c \in C \ |\ \norm{c} = 1 = \varepsilon(c) \}$. Then $(C, S_C) \in \Ob(\QQ).$
  Moreover, the double dual isomorphism $d \colon C \cong C^{**} : \mathrm{Id} \to \vNCoalg$
  lifts to an isomorphism of $\QQ$ in the following way
  \[ d \colon (C, S_C) \cong (C^{**}, \{\varepsilon\}^{\circ}) = (C^*, \{\varepsilon\})^*, \]
  with $(C^*, \{\varepsilon\}) = H(C^*)$ the canonical inclusion of the vN-algebra $C^*$.
\end{proposition}
\begin{proof}
  By Theorem \ref{thm:alg-coalg-dual} we know that $C^*$ is a vN-algebra with unit $\varepsilon.$ Proposition \ref{prop:vn-algebra-unitary-object} shows that $(C^*, \{ \varepsilon \} ) \in \Ob(\QQ)$
  and Theorem \ref{thm:q-mall} shows that $(C^{**}, \{\varepsilon\}^\polar)$ is also an object of $\QQ.$
  By definition
  \[ \{\varepsilon\}^\polar = \left \{ f \in \Ball(C^{**}) \ |\ f(\varepsilon) = 1 \right \} . \]
  But for $f \in \{\varepsilon\}^\polar$, we have $1 = \norm{f(\varepsilon)} \leq \norm{f} \norm{\varepsilon} \leq \norm{f} \leq 1$
  and therefore $f(\varepsilon) = 1 = \norm{f}_{C^{**}}.$
  Let
  \[ V \eqdef \left \{ f \in C^{**} \ |\ f(\varepsilon) = 1 = \norm{f}_{C^{**} } \right \} = \{\varepsilon\}^\circ . \]
  We have a completely isometric isomorphism $d \colon C \cong C^{**} :: c \mapsto (\varphi \mapsto \varphi(c))$
  and by Lemma \ref{lem:iso-lift}, it follows that
  $(C, d^{-1}\left[ V \right])$ is an object of $\QQ.$ It is easy to see that
  \[ d^{-1}[V] = \left \{ c \in C \ |\ \varepsilon(c) = 1 = \norm{c}  \right \} = S_C \]
  which completes the proof. \qed
\end{proof}

\begin{remark}
  Note that when $\dim(C) = 0$, the above proposition shows that $(0, \varnothing) \in \Ob(\QQ)$. In fact, this is the initial object of $\QQ$ (\cref{thm:q-mall}).
\end{remark}

\begin{corollary}
  \label{cor:coalg-well-defined}
  If $C$ is a f.d. vN-coalgebra, then $(C, P_C) \in \Ob(\QQ).$
\end{corollary}
\begin{proof}
  Combine Proposition \ref{prop:vn-coalgebra-trace-object} and \cref{cor:density}.
  \qed
\end{proof}

We can also show that the inclusion functors $H$ and $S$ behave well with respect to the duality.

\begin{proposition}
  \label{prop:hs-duals-commute}
  Let $A$ be a f.d. vN-algebra and $C$ a f.d. vN-coalgebra. Then $H(A)^* = S(A^*)$ and $S(C)^* = H(C^*).$
\end{proposition}
\begin{proof}
  We verify
  \begin{align*}
    H(A)^* &= (A, \{ 1_A\})^* & &  \\
           &= (A^*, \{ 1_A\}^\circ )  & &  \\
           &= (A^*, \{ f \in \Ball(A^*) \ | \ f(1_A) = 1 \} )  & &  \\
           &= (A^*, \{ f \in A \ |\ f(1_A) = 1 = \norm{f} \} )  & &  \\
           &= (A^*, P_{A^*} )  & & \text{(Corollary \ref{cor:density})} \\
           &= S(A^*)  & &
  \end{align*}
  and for the other statement we have
  \begin{align*}
    S(C)^* &= (C, P_C)^* & &  \\
           &= (C^*, P_C^\circ) & &  \\
           &= (C^*, \{ \epsilon \})  & & \text{(Proposition \ref{prop:vn-coalgebra-trace-object})} \\
           &= H(C^*) & &
  \end{align*}
  which completes the proof. \qed
\end{proof}

In order to prove that we have a fully faithful inclusion $S \colon \SSS \fullsubright \QQ$, we
prove an extra proposition first.

\begin{proposition}
  Let $f \colon (C, P_C) \to (D, P_D)$ be a morphism in $\QQ$, where $C = \loneplus_{1 \leq i \leq n} T_{k_i}$
  is a concrete vN-coalgebra and $D$ a f.d. vN-coalgebra. Then $f$ is counital, i.e. trace-preserving.
\end{proposition}
\begin{proof}
  Let $c \in C$ be a positive element with $c \neq 0$. From the results in Appendix \ref{app:vn-coalg-cp}, we know that
  $c = (c_1, \ldots, c_n)$ where each $c_i \in T_{k_i}$ is a positive matrix. It is easy to see that
  \begin{align*}
    \norm{c} = \sum_i \trnorm{c_i} = \sum_i \tr(c_i) = \varepsilon_C(c)
  \end{align*}
  using the fact that the trace norm coincides with the trace on positive matrices. If we define $c' \eqdef \frac{c}{\norm{c}}$,
  then we see that $c' \in P_C$ and therefore $f(c') \in P_D$. Therefore
  $\varepsilon_D(f(c')) = 1 = \varepsilon_C(c')$ and by multiplying by $\norm{c}$
  we obtain $\varepsilon_D(f(c)) = \varepsilon_C(c)$. Therefore $f$ is trace-preserving
  on non-zero positive elements. But it is obviously trace-preserving for $c=0$, so $f$ is
  trace-preserving on all positive elements of $C.$ For the general case, recall that every element in (the underlying vector space of) $C$
  is a finite linear combination of positive elements of $C$ (special case of \cite[Corollary 4.2.4]{kadison-ringrose}), so if
  $c = \sum_j \lambda_j p_j$, where each $p_j \in C$ is positive, then
  \[ \varepsilon_D(f(c)) = \sum_j \lambda_j \varepsilon_D(f(p_j)) = \sum_j \lambda_j \varepsilon_C(p_j) = \varepsilon_C(c), \]
  and therefore $f$ is indeed counital, i.e. trace-preserving. \qed
\end{proof}

\begin{corollary}
  \label{cor:vn-coalg-full}
  Let $f \colon (C, P_C) \to (D, P_D)$ be a morphism in $\QQ$, where $C$ and $D$ are f.d. vN-coalgebras.
  Then $f$ is counital, i.e. trace-preserving.
\end{corollary}
\begin{proof}
  Combine the previous proposition with \cref{ex:coalg} which states that $C$ is isomorphic as a vN-coalgebra to a concrete one.
\end{proof}

\QvNCoalg*
\begin{proof}
  Corollary \ref{cor:coalg-well-defined} shows that the functor $S$ is well-defined on objects.
  From Proposition \ref{prop:coalg-cc-cp} it follows that the morphisms in $\SSS$ may be equivalently described as the completely contractive
  trace-preserving maps. To prove that $S$ is well-defined on morphisms, it remains to check
  that for a completely contractive trace-preserving map $f \colon C \to D$, we have $f[C] \subseteq D.$
  If $c \in P_c$, then $\varepsilon_D(f(c)) = \varepsilon_C(c) = 1,$ using trace-preservation and therefore
  \[ 1 = \norm{\varepsilon_D(f(c))} \leq \norm{f(c)} \leq \norm{c} = 1  \]
  so that $\norm{f(c)} = 1 = \varepsilon_D(f(c)),$ as required.
  Faithfulness is obvious and fullness follows by \cref{cor:vn-coalg-full}.

  We can now prove that $S$ strictly preserves $\ptimes$
  by verifying
  \begin{align*}
    S(C) \ptimes S(D) &\xrightarrow{d} (S(C) \ptimes S(D))^{**} & &  \\
         &\xrightarrow{\cong} (S(C)^* \itimes S(D)^*)^{*} & & \\
         &= (H(C^*) \itimes H(D^*))^{*} & & \text{(Proposition \ref{prop:hs-duals-commute})} \\
         &= H(C^* \itimes D^*)^{*} & & \text{(\cref{thm:vn-alg-include})} \\
         &= S((C^* \itimes D^*)^{*}) & &  \text{(Proposition \ref{prop:hs-duals-commute})} \\
         &\xrightarrow{S(\cong^{-1})}  S((C \ptimes D)^{**}) & & \\
         &\xrightarrow{S(d^{-1})}  S(C \ptimes D) & &
  \end{align*}
  and recognising that the above composite is the identity. The proof for $\loneplus$ is fully analogous.
  \qed
\end{proof}

%%%%%%%%%%%%%%%%%%%%%%%%%%%%%%%%%%%%%%
\subsection{Bipolar sets relevant for pure state computation}
\label{app:qsw}
%%%%%%%%%%%%%%%%%%%%%%%%%%%%%%%%%%%%%%

Finally, we prove the claims that we made in the final paragraph of the paper where we discuss pure state computation.
Note that in the following proposition we are \emph{not} using the $\itimes$ tensor, so $M_n \ptimes M_n$ need not be a vN-algebra, nor a vN-coalgebra.

\begin{proposition}
  Let $u,v \in M_n$ be two unitary matrices. Then we have the following object equalities in $\QQ$
  \[ (M_n, \{u\}) \ptimes (M_n, \{v\}) = (M_n \ptimes M_n, \{u \otimes v\}^\bipolar) = (M_n \ptimes M_n, \{u \otimes v\}) . \]
\end{proposition}
\begin{proof}
  The first equality holds by definition, so we only have to show $\{u \otimes v\}^\bipolar = \{u \otimes v\}.$
  The inclusion $\{u \otimes v\}^\bipolar \supseteq \{u \otimes v\}$ is obvious, so we prove the other one.
  For simplicity, we identify $\MM_n \otimes \MM_n$ with $\MM_{n^2}$ and write $X$ for the operator space with the induced matrix norm from $M_n \ptimes M_n$, so that we have
  a c.i.i. $X \cong M_n \ptimes M_n$ given by the aforementioned isomorphism.
  Note that for any $a \in M_n \itimes M_n$ we have that $\norm{a}_{M_n \itimes M_n} \leq \norm{a}_{M_n \ptimes M_n}$ and therefore
  for any $a \in \MM_{n^2}$ we have
  $\norm{a}_{M_{n^2}} \leq \norm{a}_X,$ because $M_n \itimes M_n \cong M_{n^2}.$
  Next, observe that the linear functional $\frac{1}{n^2}\tr((u^* \otimes v^*)\ \cdot) \colon X \to \mathbb C$ satisfies
  \[ \frac{1}{n^2}\tr((u^* \otimes v^*) (u \otimes v) ) = 1 \]
  and this is a (necessarily complete) contraction, because for any $x \in X$ we have
  \begin{align*}
    \left| \frac{1}{n^2} \tr((u^* \otimes v^*) x) \right| &= \frac{1}{n^2} \left| \tr((u^* \otimes v^*) x) \right| \\
      &\leq \frac{1}{n^2} \norm{u^* \otimes v^*}_{\tr} \norm{x}_{M_{n^2}} \\
      &= \norm{x}_{M_{n^2}} \\
      &\leq \norm{x}_{X} .
  \end{align*}
  Therefore $\frac{1}{n^2}\tr((u^* \otimes v^*)\ \cdot) \in \{ u \otimes v\}^\circ$ and we complete the proof via Lemma \ref{lem:unitary-polar} which shows
  that $\{u \otimes v\}^\bipolar \subseteq \{u \otimes v\}.$
  \qed
\end{proof}

\begin{proposition}
  Let $u,v \in M_n$ be two unitary matrices. Then
  \[ (M_n \hagertimes M_n, \{u \otimes v\}) \in \Ob(\QQ) . \]
\end{proposition}
\begin{proof}
  The proof is analogous to the previous proposition, because for any $a \in \MM_n \otimes \MM_n$
  we have that $\norm{a}_{M_n \itimes M_n} \leq \norm{a}_{M_n \hagertimes M_n} .$
  \qed
\end{proof}

\begin{proposition}
  Let $u, v \in \MM_n$ be two unitary matrices. Then
  \[ \left( M_{2n}, \{ \ket 0 \bra 0 \otimes (uv)  + \ket 1 \bra 1 \otimes (vu) \} \right) \in \Ob(\QQ). \]
\end{proposition}
\begin{proof}
  Straightforward verification shows that $\ket 0 \bra 0 \otimes (uv)  + \ket 1 \bra 1 \otimes (vu)$ is a unitary matrix in $\MM_{2n}$, so the proposition follows
  by \cref{lem:unitary-Mn}.
  \qed
\end{proof}

\begin{proposition}
  Let $u, v \in \MM_n$ be two unitary matrices. Then
  \[ \qsw \colon (M_n \ptimes M_n, \{u \otimes v\}) \to \left(M_{2n}, \{ \ket 0 \bra 0 \otimes (uv)  + \ket 1 \bra 1 \otimes (vu) \} \right) \]
  is a morphism of $\QQ$, but
  \[ \qsw \colon (M_n \hagertimes M_n, \{u \otimes v \}) \to \left(M_{2n}, \{ \ket 0 \bra 0 \otimes (uv)  + \ket 1 \bra 1 \otimes (vu) \} \right) \]
  is \emph{not} a morphism of $\QQ.$
\end{proposition}
\begin{proof}
  All the pairs are valid objects of $\QQ$ and $\qsw \colon M_n \ptimes M_n \to M_{2n} $ is a complete contraction \cite[\S 4]{os-lics}.
  Obviously, $\qsw(u \otimes v ) = \ket 0 \bra 0 \otimes (uv)  + \ket 1 \bra 1 \otimes (vu)$, so
  \[ \qsw \colon (M_n \ptimes M_n, \{u \otimes v\}) \to \left(M_{2n}, \{ \ket 0 \bra 0 \otimes (uv)  + \ket 1 \bra 1 \otimes (vu) \} \right) \]
  is a valid morphism of $\QQ.$
  However, $\qsw \colon M_n \hagertimes M_n \to M_{2n}$ is \emph{not} a complete contraction \cite[\S 4]{os-lics}, so we do not have a morphism
  \[ \qsw \colon (M_n \hagertimes M_n, \{u \otimes v \}) \to \left(M_{2n}, \{ \ket 0 \bra 0 \otimes (uv)  + \ket 1 \bra 1 \otimes (vu) \} \right) . \]
  \qed
\end{proof}

\end{document}